%% file: main/main.tex
\documentclass{article}

\usepackage{todonotes}
\usepackage{constants}
\newconstantfamily{epsilons}{
symbol=\varepsilon
}
\newconstantfamily{k}{
symbol=k
}
\input{preamble-fabian.tex}

\usepackage{siunitx}
\newcommand{\numLambdaLower}{0.01610}
\newcommand{\numLambdaUpper}{0.6405}
\newcommand{\numKtilde}{60}
\newcommand{\numCOneTilde}{4400}
\newcommand{\numCTwoTilde}{2.98e7}
\newcommand{\numMetricComparison}{2.177}
\newcommand{\numRicciBound}{4.45e-13}
\newcommand{\numResidualCZero}{5.0e-20}

\newcommand{\numConnectionLength}{%
  \fpeval{sqrt(3)*(pi/4-atan(3^(-1/4)))
          +1/4+pi/8+ln(2)/4}}

\newcommand{\numUpperBoundDiam}{%
  \fpeval{ceil(
    sqrt(\numMetricComparison/(pi*sqrt(2)))*
    (2*\numConnectionLength+4*pi/(3*sqrt(3))),3)}}

\newcommand{\numUpperBoundQuotientDiam}{%
  \fpeval{ceil(
    2*\numConnectionLength*
    sqrt(\numMetricComparison/(pi*sqrt(2))),3)}}

\newcommand{\numLambdaOne}{%
  \fpeval{floor(
    (pi/\numUpperBoundQuotientDiam
     -3*\numRicciBound*\numUpperBoundQuotientDiam/(4*pi))^2,2)}}

\author{Michael R.~Douglas\thanks{%
    Center of Mathematical Sciences and Applications,
    Harvard University, Cambridge, MA, USA.}
    \qquad
  Javier G\'omez-Serrano\thanks{%
    Department of Mathematics,
    Brown University, Providence, RI, USA.}
    \qquad
  Fabian Lehmann\thanks{%
    Institut f\"ur Mathematik,
    Humboldt-Universit\"at zu Berlin, Berlin, Germany.}
    \\[0.6em]
  Daniel Platt\thanks{%
    Department of Mathematics,
    Imperial College London, London, United Kingdom.}
    \qquad
  Yidi Qi\thanks{%
    John A.~Paulson School of Engineering and Applied Sciences,
    Harvard University, Cambridge, MA, USA.}
    }
\date{\today}
\title{A constructive computer-assisted existence proof\texorpdfstring{\\}{ } for a Calabi--Yau metric}

\usepackage{booktabs}

\begin{document}

\begingroup
\let\|\Vert
\maketitle
\endgroup

\begin{abstract}
    We give a rigorous, constructive proof that there exists a Ricci-flat K\"ahler metric on the Fermat quartic hypersurface, independently of Yau's famous theorem,
    along with an explicit and accurate analytic approximation to the metric.
    The proof combines numerical techniques to obtain the approximately Ricci flat metric, 
    a fixed point theorem to prove existence, and verified numerics to get the bounds used by the fixed point theorem.
    The method is general and we discuss the computational issues in applying it to other manifolds.
\end{abstract}

\tableofcontents

\section{Introduction}

Calabi-Yau manifolds are complex manifolds that are central objects in superstring theory and pure mathematics.
In differential geometry, they constitute the most fertile supply of Ricci-flat metrics thanks to Yau's proof of the Calabi conjecture in \cite{Yau1977,Yau1978}.
In superstring theory, they appear as compactifying dimensions, used to model the laws of physics.
In algebraic geometry, they play a distinguished role in birational geometry and the minimal model program.

For applications in geometry and physics it is often necessary to work with the Ricci-flat metric, but Yau's theorem alone does not give much beyond its existence, which makes it notoriously difficult to work with.
This shortcoming has been addressed in different ways:

First, the SYZ conjecture, motivated by T-duality in string theory \cite[Section~1]{StromingerYauZaslow1996} makes some predictions for Calabi-Yau metrics near certain singular limits.
Thanks to recent breakthrough work in \cite{Li2022,Li2026}, some of these predictions have now been proven.
Second, near other degenerate limits methods from geometric analysis allow to give fairly detailed descriptions of Calabi-Yau metrics.
Examples are too many to list exhaustively, but we mention the classical Kummer construction in \cite{Topiwala1987I,Topiwala1987II}, conical asymptotics \cite{HeinSun2017}, and gluing constructions of collapsing metrics \cite{GrossWilson2000,Foscolo2019,SunZhang2019}.
Third, there is now a vast literature on numerically solving for the Ricci-flat metric.
References are too numerous to list, but we mention the pioneering work \cite{HeadrickWiseman2005}, the fundamental work suggesting to use algebraic metrics in numerical optimisation \cite{Donaldson2009}, and the recent work \cite{DouglasLakshminarasimhanQi2022} as a representative of a large body of work applying machine learning for this approximation task.
In special cases it is also possible to write down series that converge to the Calabi-Yau metric, see \cite{FuLiuXiaoYau2026Kummer,FuLiuXiaoYau2026Schoen,KachruTripathyZimet2020}.

The problem with these approaches is that they either describe metrics near degenerate limits, 
or they give numerical candidates for which it is not known how close they are to the true Calabi-Yau metric.
We are not aware of a quantitative description of genuine Calabi-Yau metrics away from a degenerate limit. 
This is a particularly severe limitation in Physics, where numerous important calculations have been carried out with respect to numerical metrics far away from degenerate limits, but there is currently no error bound at all for how close the computed data are to the data for the true Calabi-Yau metric.
Examples include computations of physical Yukawa couplings \cite{ConstantinEtAl2025} and of Kaluza--Klein spectra testing the swampland distance conjecture \cite{AshmoreRuehle2021}.

We address these issues by:
\begin{enumerate}
\item
Explaining an algorithm for computing an ultra high-accuracy approximate Calabi-Yau metric $\omega$ and applying it to compute an approximate metric on the Fermat quartic whose error is far smaller than that of other candidates if evaluated using the non-rigorous Monte-Carlo estimation that is standard in the literature.
This method is readily extended to other complete intersection Calabi-Yau manifolds.

\item
We certify the smallness of our error using rigorous computing.
While certified error bounds are standard in numerically verified proofs, certifying the tiny errors in our work required some engineering implementing rigorous bounds using Taylor models using graphics processing unit (GPU) computing.
The Taylor-model arithmetic is described in the technical report \cite{technical-report-reference}; the steps for the Fermat quartic are stated in \cref{section:rigorous-error-bounds}, and the code and further details are in our repository \cite{code-repository}.
To our knowledge, these are the first rigorous error bounds for numerical Calabi-Yau metrics.

\item
A small error bound alone does not imply that the true Calabi-Yau metric is nearby.
We then proceed to use elliptic analysis to prove the existence of a Calabi-Yau metric.
This proof method automatically gives a bound for the difference between the approximate and the true Calabi-Yau metric.
Existence is nowhere assumed, so this is an existence proof that is independent of Yau's original existence proof.
\end{enumerate}

Let 
\[
M=\{(X_0:X_1:X_2:X_3) \in \mathbb{CP}^3: X_0^4+X_1^4+X_2^4+X_3^4=0\}
\]
denote the Fermat quartic.
Our main theorem is then:

\begin{theorem}
\label{theorem:main}
Let $\omega$ be the explicitly computed Kähler form determined by
\eqref{eqn:certified-potential}, and let $g$ be its associated metric.
There exists $u \in C^\infty(M)$ such that $\omega+i \partial \bar \partial u$ is Calabi-Yau, and
\[
\|{u}_{\mathscr{X}} \leq 1.551\cdot10^{-4}.
\]
In particular
\[
\|{i\partial\bar\partial u}_{C^0(M,g)}
\leq 1.553\cdot10^{-4}.
\]
\end{theorem}

Here, $\|{\cdot}_{\mathscr{X}}$ denotes a norm that is a weighted sum of the $C^{2,\alpha}$-norm and the $H^2$-norm defined in \cref{section:functionspaces}.
Control of the Kähler potential in the $\mathscr{X}$-norm gives a bound for the difference of the Kähler forms in the ordinary $C^0$-norm, as also stated in the theorem.

The article is structured as follows:
\begin{itemize}
\item
In \cref{section:high-accuracy} we describe our method to compute ultra high-accuracy approximate Calabi-Yau metrics.
During optimisation, accuracy is measured using non-rigorous Monte-Carlo integration.
We then rigorously bound the approximation error using verified computing.

\item
In \cref{section:bounding-quantities} we rigorously bound the geometric quantities needed for the later functional analysis.
These are:
bounds for the metric in coordinates and derivative estimates for transition functions of an explicit atlas on $M$;
an upper bound for the diameter of $M$;
a lower bound for the Ricci tensor;
a lower bound for the spectral gap of the Riemannian Laplace operator acting on symmetry-invariant functions.
All these bounds are for the approximate Calabi-Yau metric $\omega$.
Geometric analysis proofs routinely use that these quantities are bounded by a constant, but we exhibit \emph{fully explicit} numerical bounds for all of them.

\item
In \cref{section:perturbation-theory}, we use pen-and-paper geometric analysis techniques to prove that there exists a genuine Calabi-Yau metric near our approximation.
We begin in \cref{section:analytic-set-up} by defining an atlas of small ball charts depending on the explicit atlas defined before and bound their geometry.
In \cref{section:functionspaces} we define the composite $\mathscr{X}$-norm and $\mathscr{Y}$-norm and state the residual bound in these norms.
The most important functional analytic bound is an explicit operator norm bound for the inverse of the Laplacian, which we derive in \cref{subsection:estimate-inverse-laplacian}.
\end{itemize}

\paragraph{Acknowledgments.}
The authors thank Simon Donaldson, Simon Salamon, and Freid Tong for helpful conversations. JGS was partially supported by the National Science Foundation
through Grants DMS-2247537, DMS-2434314, and DMS-2554957, and by a Simons Fellowship.
FL was funded by the Deutsche Forschungsgemeinschaft (DFG, German Research Foundation) – Projektnummer  573641031
and partially
funded by the DFG under Germany´s Excellence Strategy – The Berlin Mathematics 
Research Center MATH+ (EXC-2046/1, project ID: 390685689).
DP was supported by the Eric and
Wendy Schmidt AI in Science Postdoctoral Fellowship, a Schmidt Futures program.

\section{Approximate Calabi-Yau metrics}
\label{section:high-accuracy}

Denote the Fubini-Study Kähler form restricted to the Fermat quartic and normalised to volume one by $\omega_0$ and the induced metric by $g_0$.
In this section we explain how to compute a Kähler form $\omega$ with induced metric $g$ that is close to being Ricci-flat and we will refer to this metric throughout as the \emph{approximate Calabi-Yau metric}, even though a priori it is not clear that $\omega$ is anywhere near the true Calabi-Yau metric, which is known to exist on abstract grounds.

\subsection{Computing ultra high-accuracy approximate Calabi-Yau metrics}
\label{subsection:computing-metrics}

In this subsection we explain how we compute approximate Calabi-Yau metrics and how we estimate their error.
We first describe the general method, and then the metric that we use in the rest of the paper.
The residual numbers in this subsection are non-rigorous Monte-Carlo estimates; rigorous bounds follow in \cref{subsection:certified-errors-bounds}.
\Cref{section:numerical-details} contains further details of the method and the smaller residuals that it reaches.

Our Kähler potentials are algebraic potentials in the sense of \cite{Donaldson2009} with explicit corrections.
Fits of the type of Headrick and Nassar \cite{headrick_energy_2009}, which fit the coefficients of one Hermitian form $p$ so that $\log p$ minimises the residual, stop near degree $28$, where the products of sections are almost linearly dependent on $M$.
We therefore represent the potential in an orthonormal basis on $M$, which we build by a recurrence of Vandermonde-with-Arnoldi type \cite{BrubeckNakatsukasaTrefethen2021} and evaluate in double-double arithmetic \cite{Dekker1971}.
Such a basis spans a linear space, so we write the potential as a fixed base of low degree plus a correction of higher degree, which approximates the difference between a potential of the Calabi--Yau metric and the potential of the base.
To first order the residual depends linearly on the correction, so one linear least-squares problem gives it.
This removes the limit of the fits of Headrick and Nassar, and a correction with many consecutive degrees then reaches a small residual at a lower top degree, which reduces the cost of the rigorous bounds.
The perturbation argument of \cref{section:perturbation-theory} needs much less accuracy than this method can give, so in this paper we use corrections up to degree $24$, which give a residual $f$ with $\|{f}_{L^1(M)}\approx\num{1.7e-22}$.

\textbf{The potential.}
All our potentials $\Phi$ are functions on $\CC^4\setminus\{0\}$ such that $\Phi-\log p_0$ is a smooth function on $M$ for a Hermitian form $p_0$ of bidegree $(k_0,k_0)$ with $p_0>0$ on $M$.
Then $i\partial\bar\partial\Phi$ is a closed real $(1,1)$-form on $M$, cohomologous to $k_0\,i\partial\bar\partial\log|X|^2$.
If it is positive, we let $\omega$ be the constant multiple of $i\partial\bar\partial\Phi$ with $\vol(M,\mu_\omega)=1$, so that $\omega$ is cohomologous to $\omega_0$.

We impose the symmetries of $M$: every potential is invariant under the group $G$ of order $1536$ from \cref{subsection:diameter-bound} and under complex conjugation of the coordinates.
Then a Hermitian form $N(X)=\sum_{A,B}N_{AB}X^A\bar X^B$ of bidegree $(k,k)$ has a real symmetric coefficient matrix, and $N_{AB}=0$ unless $A\equiv B \pmod 4$ in each component.
So the matrix has $64$ blocks, one for each residue $A\bmod 4$, and the permutations in $G$ map blocks to blocks.
This leaves $135$ free parameters at degree $k=16$ and $542$ at degree $k=24$ (\cref{subsection:invariant-forms}).

Our potentials have the form
\begin{align}
\label{eqn:one-stage-potential}
    \Phi_1=\log p_0+\sum_{L\in\mathcal L_1}N_L\,p_0^{-L/k_0}
\end{align}
for a finite set $\mathcal L_1$ of degrees.
Here $p_0(X)=\sum_{A,B}H_{AB}X^A\bar X^B$ with a $G$-invariant real symmetric matrix $H$, and each $N_L$ is a real $G$-invariant Hermitian form of bidegree $(L,L)$.
The base $p_0$ is a fit of the type of Headrick and Nassar: we fit $H$ so that $\log p_0$ minimises the residual.
In our experiments such fits do not improve beyond a mean residual of about $\num{3e-14}$, which they reach at degree $28$.
At these degrees the products of sections are almost linearly dependent on $M$: the map from the coefficients of $H$ to the values of $p_0$ has a condition number up to about $10^{20}$.
We therefore keep $k_0$ low and correct $\log p_0$ by the sum in \eqref{eqn:one-stage-potential}.

\textbf{The levels.}
Each correction term must be homogeneous of degree $0$ in $X$.
A common choice is $N_L|X|^{-2L}$; we take $N_Lp_0^{-L/k_0}$ instead, and for $L>0$ we call the space $V_L$ of all these functions the \emph{level} of degree $L$.
The correction $\Phi_{\mathrm{CY}}-\log p_0$ from the base to a potential $\Phi_{\mathrm{CY}}$ of the Calabi--Yau metric is, to first order and up to a constant, a function in $V_L$, by the Tian--Yau--Zelditch expansion \cite{Tian1990,Zelditch1998}; so each level only needs to supply an almost constant function.
With the weight $|X|^{-2L}$, the same correction carries the factor $(|X|^2p_0^{-1/k_0})^L$, whose variation over $M$ grows exponentially with $L$ (\cref{subsection:level-weight}).
The same argument gives the weight $e^{-(L/k_0)\Phi_1}$ for a further correction of $\Phi_1$ (\cref{subsection:two-stage-results}).

Two further properties of the levels explain the choice of $\mathcal L_1$:
\begin{itemize}
\item
The levels are nested: $V_L\subset V_{L+k_0}$, because $Np_0^{-L/k_0}=(Np_0)\,p_0^{-(L+k_0)/k_0}$.
\item
Levels whose degrees differ modulo $k_0$ contain different fractional powers of $p_0$.
At the top degree $24$, eight consecutive levels give a residual about $100$ times smaller than levels with a degree step of $4$.
\end{itemize}
The monomials $X^A\bar X^B$ with $|A|=|B|=L$ give a very ill-conditioned basis of $V_L$ on $M$.
We therefore build an orthonormal basis of each level on a fixed sample of points with a recurrence of Vandermonde-with-Arnoldi type \cite{BrubeckNakatsukasaTrefethen2021}.
Each new basis function is an earlier one times a $G$-invariant function of degree $0$, with its projections onto the earlier functions removed.
The stored coefficients of the recurrence then evaluate the basis at any point of $M$, without large coefficients.

\textbf{The residual in a chart.}
As in \cref{section:perturbation-theory}, $\Omega$ is the holomorphic volume form of $M$, $\hat\mu$ is the multiple of $\mu_\Omega$ with volume $1$, and $\eta=\mu_\omega/\hat\mu$.
To evaluate $\eta$ at a point, we set the homogeneous coordinate of largest modulus to $1$, solve the equation of $M$ for the coordinate $X_m$ of largest modulus among the other three, and use the remaining two coordinates $t=(t_1,t_2)$ as holomorphic coordinates.
In this chart $\Omega$ is a constant multiple of $dt_1\wedge dt_2/X_m^3$, and
\begin{align}
\label{eqn:eta-chart}
    \eta=\frac{|X_m|^6}{\bar v}\det\left(\frac{\partial^2\Phi}{\partial t_i\partial\bar t_j}\right),
    \qquad
    \bar v=\frac{32\pi^4k_0^2}{\Gamma(1/4)^8}.
\end{align}
The constant $\bar v$ is exact.
It is the ratio of $\int_M\det(\partial_{t_i}\partial_{\bar t_j}\Phi)\,d^4t=2\pi^2k_0^2$, which only depends on the cohomology class, to the period $\int_M|X_m|^{-6}\,d^4t=\Gamma(1/4)^8/(16\pi^2)$.
Here $d^4t$ is the Lebesgue measure in the chart, and both integrands are volume forms on $M$.
So $\eta$ needs no normalisation from a sample, which would hide a constant error of $\mu_\omega$.
We report the residual $f=1-1/\eta$ from \cref{eqn:residual-def}; our programs compute $\eta-1$, which agrees with $f$ to second order (\cref{subsection:certified-errors-bounds}).

We use Monte-Carlo integration over $M$ in the fit and for the estimates of the residual in this subsection and in \cref{section:numerical-details}; the rigorous bounds of \cref{subsection:certified-errors-bounds} do not use it.
We sample points of $M$ as the intersections of $M$ with random lines in $\P^3$, as in \cite{DouglasLakshminarasimhanQi2022}, and each point $x_i$ has a weight $w_i$ with $\sum_iw_i=1$, such that $\sum_iw_ih(x_i)$ estimates $\int_Mh\,\hat\mu$.
We report
\[
    \operatorname{mean}|f|=\sum_iw_i|f(x_i)|,
    \qquad
    \operatorname{rms}(f)=\Bigl(\sum_iw_if(x_i)^2\Bigr)^{1/2},
\]
which estimate $\|{f}_{L^1(M)}$ and $\|{f}_{L^2(M)}$, and the sample maximum of $|f|$, which is only a lower bound for $\sup_M|f|$; we compute these numbers on hold-out points that no fit uses.
The ratio $r(x)=\min_a|X_a|/\max_a|X_a|$ describes the position of a point $x$; it is small near the curves $M\cap\{X_a=0\}$.

\textbf{Arithmetic.}
The residuals that we aim at are far below the unit roundoff $2^{-53}\approx\num{1.1e-16}$ of double precision.
We therefore evaluate every step except the linear solves in double-double arithmetic \cite{Dekker1971}, which represents a number as an unevaluated sum of two double-precision numbers.
This includes the projection of each sample point onto $M$ by a Newton step, the base, the determinant, the weights, the recurrences and the Gram--Schmidt processes.
The computations run on GPUs.
The recurrence amplifies the rounding errors of double-double arithmetic, so the error of a basis function grows with its degree.
We measure this error at sample points against the same recurrence in mpmath \cite{mpmath} at $50$ significant digits, with the same stored coefficients and the same point (\cref{subsection:level-bases}).

\textbf{The fit.}
We fit the potential in two steps:
\begin{enumerate}
\item
Fit the matrix $H$ of the base $p_0$ with the Levenberg--Marquardt method, which minimises the sample mean of $(\eta-1)^2$.
This is the only nonlinear fit.
\item
Keep $p_0$ fixed and fit the forms $N_L$, $L\in\mathcal L_1$, in one linear least-squares problem.
\end{enumerate}
We call step~2 a \emph{stage}.
With $p_0$ fixed, $\eta$ depends linearly on the coefficients of the stage, up to terms of second order in the change of the metric.
So the stage is one linear least-squares problem for the values of $\eta-1$ at $200{,}000$ fit points.
We solve it in double precision with scaled columns.
Then a few steps of iterative refinement, with $\eta-1$ at the fit points computed in double-double arithmetic, remove the second-order terms and the rounding error of the solve.

The levels of a stage share one orthonormal basis.
A pivoted Gram--Schmidt process takes the level functions from the highest degree down.
A function enters only if its new part $\nu$, the norm of its part orthogonal to the functions before it relative to its own norm, is larger than a cutoff $\tau$.
A function with the new part $\nu$ carries the error of its level function times $1/\nu$.
If $\nu$ is below that error, the new part is rounding noise, which depends on the arithmetic and on the chart of the point.
We therefore set $\tau$ near the measured error of the level functions at the top degree.
During the construction, a second chart and computations in mpmath \cite{mpmath} test whether a fit uses this noise (\cref{subsection:two-stage-results}); for the metric of this paper, the rigorous bounds of \cref{subsection:certified-errors-bounds} make these tests unnecessary.
A second stage fits the correction of $\Phi_1$ with the weights $e^{-(L/k_0)\Phi_1}$ in the same way, and in principle more stages can follow; \cref{subsection:two-stage-results} gives a metric with two stages and $\operatorname{mean}|f|\approx\num{2.7e-28}$.

\textbf{The approximate Calabi-Yau metric of this paper.}
From now on, $\omega$ denotes the metric of the one-stage potential
\begin{align}
\label{eqn:certified-potential}
    \Phi=\log p_0+\sum_{L=13}^{24}N_L\,p_0^{-L/k_0},
    \qquad
    k_0=16,
\end{align}
that is, \eqref{eqn:one-stage-potential} with $k_0=16$ and $\mathcal L_1=\{13,14,\dots,24\}$.
Further stages and higher degrees give a much smaller residual, but this metric is accurate enough for the perturbation argument of \cref{section:perturbation-theory}.
The cost of the rigorous bounds of \cref{subsection:certified-errors-bounds} grows quickly with the degree, and our Taylor-model code treats one stage only.

For the rigorous computations, we write $p_0$ with blocks that are exact copies of each other within each $S_4$-orbit, so that $p_0$ is exactly $G$-invariant.
We also convert the sum of each level into one explicit $G$-invariant Hermitian form $N_L$, by a least-squares fit in ball arithmetic \cite{python-flint} on about twice as many points of $M$ as $N_L$ has free parameters; at other test points its error is at most $\num{2e-70}$ of the largest value of $|N_L|$.
We take the stored coefficients of $H$ and of the $N_L$ as exact numbers, so $\omega$ is an explicit metric.
\Cref{tab:certified-parts,tab:certified-residual} give the parts of \eqref{eqn:certified-potential} and its residual.

\begin{table}[ht]
\centering
\begin{tabular}{@{}lcccc@{}}
\toprule
part & degrees & cutoff $\tau$ & coefficients & condition number\\
\midrule
base $p_0$ & $16$ & -- & $135$ & --\\
stage & $24,23,\dots,13$ & $10^{-23}$ & $2{,}482$ & $50$\\
\bottomrule
\end{tabular}
\caption{The parts of the potential \eqref{eqn:certified-potential}. The coefficients of the stage are those of its orthonormal basis, and the condition number is that of its least-squares matrix with scaled columns.}
\label{tab:certified-parts}
\end{table}

\begin{table}[ht]
\centering
\begin{tabular}{@{}lccc@{}}
\toprule
potential & $\operatorname{mean}|f|$ & $\operatorname{rms}(f)$ & sample maximum of $|f|$\\
\midrule
base $p_0$ & \num{3.54e-9} & \num{6.68e-9} & \num{2.31e-7}\\
base and stage & \num{1.73e-22} & \num{2.38e-22} & \num{3.29e-21}\\
\bottomrule
\end{tabular}
\caption{The residual of the potential \eqref{eqn:certified-potential} and of its base on $400{,}000$ hold-out points, in double-double arithmetic.}
\label{tab:certified-residual}
\end{table}

The solve in double precision gives $\operatorname{mean}|f|=\num{1.2e-17}$ on the fit points, about the square of the residual of the base, as expected for the second-order terms.
The first refinement step takes it to $\num{1.72e-22}$, and further steps change it by less than $0.01\%$.
The error is almost uniform over $M$: $\operatorname{mean}|f|$ is about $\num{1.6e-22}$ to $\num{1.9e-22}$ in each bin of $r$ above $0.05$, and $\num{2.5e-22}$ below $0.05$.
\Cref{subsection:certified-errors-bounds} proves $\sup_M|f|\leq\num{\numResidualCZero}$, and \cref{subsection:ricci-bound} bounds the Ricci curvature of $\omega$.

\subsection{Certified error bounds for the approximate Calabi-Yau metric}
\label{subsection:certified-errors-bounds}

The estimates of \cref{subsection:computing-metrics} evaluate the residual $f$ only at finitely many points, and the sample maximum of $|f|$ is only a lower bound for $\sup_M|f|$.
The perturbation argument of \cref{section:perturbation-theory} needs bounds on all of $M$: by Proposition~\ref{prop-residual-bound}, it needs bounds for $\|{f}_{C^0(M)}$ and $\sup_M|\nabla f|_g$.
In this subsection we prove such bounds for the approximate Calabi-Yau metric $\omega$ of \eqref{eqn:certified-potential}.
The proof is computer-assisted.
We describe the method and the computation here.
\Cref{section:rigorous-error-bounds} states each step of the proof as a lemma, and the technical report \cite{technical-report-reference} gives the arithmetic of the Taylor models.
The code and further details are in our repository \cite{code-repository}.
The same method also gives the bounds for $\omega$ in \cref{section:bounding-quantities}.

\begin{proposition}
\label{proposition:certified-residual}
Let $\omega$ be the metric of the potential \eqref{eqn:certified-potential}, with the stored coefficients of $p_0$ and of the forms $N_L$ taken as exact numbers.
\begin{enumerate}
\item
$p_0>0$ on $M$, and $i\partial\bar\partial\Phi$ is positive on $M$.
So $\omega$ is a Kähler form cohomologous to $\omega_0$, and the residual $f$ from \cref{eqn:residual-def} is a smooth function on $M$.
\item
$\|{f}_{C^0(M)}\leq\num{\numResidualCZero}$.
\item
$\sup_M|\nabla f|_g\leq\num{3.37e-16}$.
\end{enumerate}
\end{proposition}

\begin{proof}
The computation described in the rest of this subsection and in \cref{section:rigorous-error-bounds} proves (1)--(3).
\end{proof}

Together with $\vol(M,\mu_\omega)=1$ and Proposition~\ref{prop-residual-bound}, these bounds control the residual norm introduced in \cref{section:functionspaces}.

Our programs bound $\eta-1$ and $\nabla\eta$, where $\eta=\mu_\omega/\hat\mu$ as in \cref{eqn:eta-chart}.
Since $f=1-1/\eta=(\eta-1)/\eta$ and $\nabla f=\nabla\eta/\eta^2$, the bounds $|\eta-1|\leq B$ and $|\nabla\eta|_g\leq B'$ with $B<1$ give $|f|\leq B/(1-B)$ and $|\nabla f|_g\leq B'/(1-B)^2$.
For our bounds $B=\num{4.9985e-20}$ and $B'=\num{3.3654e-16}$, these differ from $B$ and $B'$ by a relative amount below $10^{-19}$.

\textbf{Taylor models.}
Interval or ball arithmetic \cite{Johansson2017} alone does not give useful bounds for $f$.
It encloses the result of each operation for all inputs in the input balls independently, so it cannot follow cancellations between correlated quantities, and $f$ is small only because of such cancellations.
For example, in a chart the forms of the potential become polynomials that are up to $\num{2e5}$ times smaller than the sum of the absolute values of their terms.
A Taylor model of a function on a small box is a polynomial together with a rigorous bound for the difference between the function and the polynomial on the box \cite{MakinoBerz2003}.
Arithmetic on Taylor models acts on the coefficients of the polynomials, so it keeps the cancellations, and it bounds only the terms that it drops and the rounding errors.

\textbf{The method.}
The proof has four steps.
\begin{enumerate}
\item
\emph{Reduction by symmetry.}
The group $G$ and complex conjugation act on $M$ by isometries of $g$ that preserve $\hat\mu$, because $\Phi$ is invariant under them.
So $f$ and $|\nabla f|_g$ are invariant under them, and it is enough to bound them at the points $[1:X_1:X_2:X_3]$ of $M$ with $t=(X_2,X_3)$ in the region
\[
    \mathcal F=\bigl\{t\in\CC^2:|X_3|\leq|X_2|\leq|X_1|\leq1,\ 0\leq\arg X_2\leq\pi/4,\ 0\leq\arg X_3<\pi/2\bigr\},
\]
where $X_1^4=-1-X_2^4-X_3^4$.
The angular conditions are imposed only on non-zero coordinates.
A permutation orders the coordinates by modulus, phases put $\arg X_2$ and $\arg X_3$ in $[0,\pi/2)$, and complex conjugation followed by a phase maps $\arg X_2$ to $\pi/2-\arg X_2$.
The four solutions $X_1$ over a point $t$ differ by phases, so they give the same values.
On $\mathcal F$ we have $|X_1|\geq3^{-1/4}$, and $t$ is the chart of \cref{eqn:eta-chart}.
\item
\emph{A finite cover.}
A grid of cubes of side $0.03$ in the real coordinates $(\Re X_2,\Im X_2,\Re X_3,\Im X_3)$ covers $\mathcal F$.
Interval arithmetic shows that all but $113{,}357$ of its $943{,}296$ cubes miss $\mathcal F$; we call the remaining cubes the \emph{cells}.
For each cell we choose a point $p\in M$ over its centre and holomorphic coordinates $w=(w_1,w_2)$ on $M$ near $p$ along the tangent plane at $p$.
Ball arithmetic proves that every point of $M$ over the cell lies in the \emph{box} $\{|w_1|,|w_2|\leq\delta\}$ of the cell.
\item
\emph{A bound on each box.}
On each box we compute Taylor models of $\eta-1$ and of its first derivatives.
Their polynomials have degree $28$ or $32$ in $w$ and $\bar w$, and their remainder bounds include all dropped terms and all rounding errors.
The models give bounds for $|\eta-1|$ and $|\nabla\eta|_g$ on the box.
They also prove that $p_0\neq0$ and that $i\partial\bar\partial\Phi$ is positive on the box.
\item
\emph{An exact check.}
A separate program proves with exact arithmetic that the boxes cover all of $\mathcal F$, that each box contains the points of $M$ over its cell, and that each bound of $|\eta-1|$ is at most $\num{5e-20}$.
\end{enumerate}
So the largest bounds of the boxes bound $|\eta-1|$ and $|\nabla\eta|_g$ on $M$.
Since $M$ is connected, $p_0$ has one sign on $M$, and one evaluation in ball arithmetic shows that $p_0>0$ (\cref{subsection:cover-and-sign}).
This step is necessary, because the coefficient matrix of $p_0$ is not positive definite.

\textbf{Rounding errors.}
The Taylor models run on GPUs in double-double and triple-double arithmetic, which store a number as an unevaluated sum of two or three double-precision numbers.
Each operation also computes a rigorous bound for its own rounding error from its computed result, in the manner of Arb \cite{Johansson2017}: the error of a rounded sum with the result $r$ is at most $u|r|$, and the error of a rounded product is at most $u|r|+2^{-1074}$, where $u=2^{-53}$ \cite{IEEE754}.
The arithmetic on these bounds rounds upwards.
So no bound depends on a count of operations, on the order of the parts of a number, or on the absence of underflow.
The numbers that the GPUs take as input, such as the frames of the boxes and the re-expanded forms of the potential, come from ball arithmetic at $256$ or $512$ bits or from exact rational arithmetic, with rigorous bounds for their rounding.

\textbf{The computation.}
We aimed at the bound $\num{5e-20}$ for every box.
Boxes with $r<0.3$ use degree $28$ and the others degree $32$.
A box whose bound exceeded the target was computed again with degree $32$ and then, if necessary, split into the cubes of half the side that meet $\mathcal F$; no box needed a second split.
\Cref{tab:certified-run} summarises the result.
The largest bound of $|\eta-1|$ is at $r\approx0.23$, and $97\%$ of it is the remainder of the Taylor model.
At the centre of one box $|f|\approx\num{4.1e-21}$, so our bound exceeds $\sup_M|f|$ by at most a factor of about $12$.
The sample maximum $\num{3.29e-21}$ of \cref{tab:certified-residual} lies below this value.

\begin{table}[ht]
\centering
\begin{tabular}{@{}ll@{}}
\toprule
quantity & value\\
\midrule
cells & $113{,}357$\\
boxes of the proof & $134{,}552$, of which $27{,}695$ come from $3{,}538$ split cells\\
largest bound of $|\eta-1|$ & \num{4.9985e-20}\\
median bound of $|\eta-1|$ & \num{1.8e-21}\\
largest $|\eta-1|$ at a box centre & \num{4.06e-21}\\
largest bound of $|\nabla\eta|_g$ & \num{3.3654e-16}\\
cost & $12.5$ GPU-hours on NVIDIA H200 GPUs\\
\bottomrule
\end{tabular}
\caption{The computation that proves Proposition~\ref{proposition:certified-residual}. The cost includes the bounds for the metric of \cref{subsection:certified-metric-bounds}.}
\label{tab:certified-run}
\end{table}

\textbf{Tests.}
The exact check of step 4 is part of the proof.
We also tested the computation along independent paths; these tests are not part of the proof.
\begin{itemize}
\item
$3{,}000$ random points of $\mathcal F$ each lie in a box of the proof.
\item
At random points of $8$ random boxes, the Taylor model of $\eta-1$ agrees with $\eta-1$ computed in ball arithmetic in the chart coordinates $t$, within $\num{1.9e-5}$ of its remainder bound.
\item
For $8$ boxes, each with $4$ radii, every intermediate Taylor model of the computation encloses its function at $512$ points of the box, where an independent program computes the function in ball arithmetic.
\item
At $2$ points of each of the $332$ reference centres of \cref{subsection:box-model}, the re-expanded forms of the potential agree with $p_0$ and the $N_L$ evaluated directly in homogeneous coordinates, within their error bounds.
\end{itemize}

\textbf{Assumptions.}
The proof assumes that the GPUs and the CPUs execute the floating-point operations according to the IEEE~754 standard \cite{IEEE754} with the requested rounding modes, that Arb \cite{Johansson2017,python-flint} and the interval arithmetic of mpmath \cite{mpmath} return correct enclosures, and that our programs implement the method of \cref{section:rigorous-error-bounds} correctly.
The tests above support the last assumption but do not prove it.

\section{Bounding geometric quantities}
\label{section:bounding-quantities}

In this section we rigorously bound the geometric quantities for the metric $g$ constructed in \cref{section:high-accuracy} that will be needed for our later functional analysis.
Throughout, we write $\Delta_g=\operatorname{div}_g\nabla^g$ and compute spectral gaps for $-\Delta_g$.
In \cref{subsection:atlas-bounds}, we begin by defining one explicit atlas for the quartic $M$ and then prove bounds depending on this choice of atlas:
bounds for the metric in local coordinates, and bounds for the transition functions of the atlas.
We then proceed in \cref{subsection:diameter-bound,subsection:ricci-bound,subsection:spectral-gap} to compute an upper bound for the diameter, a lower bound for $\Ric$, and a lower bound for the spectral gap of the Riemannian Laplacian acting on symmetry-invariant functions, respectively.
All bounds are proved and do not rely on non-rigorous numerical evaluations.

\subsection{Bounds for a fixed atlas on $M$}
\label{subsection:atlas-bounds}

In the later functional analysis we will use an atlas using \emph{small ball charts}.
Like in the typical proofs for elliptic theory on manifolds, we rely on these balls being small enough to be close to Euclidean.
Their number is far too large to enumerate them, so we instead define one fully explicit \emph{parent atlas}, consisting of just a small number of charts that are explicitly written down.
In \cref{section:analytic-set-up}, we will then explain how to construct the small ball charts from the parent atlas.
We will be able to estimate how close to Euclidean the small ball charts are, based on some explicit bounds for the parent charts, which themselves need not be close to Euclidean.
The purpose of this subsection is to explain the definition of the parent atlas and compute the explicit bounds needed later.

As before, we write $X_0,X_1,X_2,X_3$ for the homogeneous coordinates on $\mathbb{CP}^3$ and work on the Fermat quartic $M=Z(X_0^4+\dots+X_3^4)$.
We define an atlas as follows:

\textbf{Construction of two charts.}
Let $U_0=\{X \in \mathbb{CP}^3: X_0=1\}$ be a standard affine open.
On this chart,
\begin{align}
\label{eqn:s-def}
X_1=s:=\sqrt[4]{-1-X_2^4-X_3^4}.
\end{align}
The 4th root is unique if we ask $|X_1|>0$ and fix either $\arg(X_1) \in (-40^\circ, 40^\circ)$ or $\arg(X_1) \in (45^\circ-40^\circ, 45^\circ+40^\circ)$.
It's important these intervals are smaller than $360^\circ/4=90^\circ$, because on bigger intervals the 4th root isn't unique.
At $X_1=0$ the fourth root function is continuous but not smooth, so the chart degenerates there.
Because of this, we require on our chart that
$|X_1| > 0.98 \max_{j \neq 0,1} \{|X_j|\}$.
At large $X_j$ for $j \neq 0$, the map sending a point to $X_2,X_3$ degenerates, so we demand $|X_0|>0.99 \max_{j \neq 0} \{ |X_j| \}$.
The final chart definitions are:
\begin{align*}
U^*&:=
\left\{
X \in M\cap U_0 :
\arg(X_1) \in (-40^\circ, 40^\circ),
|X_0|>0.99 \max_{j \neq 0} \{ |X_j| \},
|X_1| > 0.98 \max_{j \neq 0,1} \{|X_j|\}
\right\},
\\
U^{**}&:=
\left\{
X \in M\cap U_0 :
\arg(X_1) \in (45^\circ-40^\circ, 45^\circ+40^\circ),
|X_0|>0.99 \max_{j \neq 0} \{ |X_j| \},
|X_1| > 0.98 \max_{j \neq 0,1} \{|X_j|\}
\right\}
\end{align*}
with the map
\begin{align*}
    \varphi^*: U^* & \rightarrow \mathbb{C}^2
    \\
    [1:X_1:X_2:X_3] & \mapsto (X_2,X_3)
\end{align*}
and analogously for $U^{**}$.
Note that $\arg(X_1)$ is defined, because a point on $M$ must have at least two non-zero coordinates, and because $|X_1| > 0.98 \max_{j \neq 0,1} \{|X_j|\}$ we have $X_1 \neq 0$.

\textbf{Definition of the atlas.}
Define
\begin{align}
    \label{equation:explicit-atlas}
    \mathcal{A}
    :=
    \{
    g \cdot U^*, g \cdot U^{**}:g\in G
    \}
    \text{ where }
    G=((\mu_4)^4/\mu_4)\rtimes S_4.
\end{align}

Here $\mu_4$ denotes the group of fourth roots of unity, the quotient is by diagonal scalars, and $S_4$ permutes the coordinates.

\begin{proposition}
    The charts of $\mathcal{A}$ from \cref{equation:explicit-atlas} cover all of $M$.
\end{proposition}

\begin{proof}
Let $y \in M$.
Let $g \in G$ be a permutation of coordinates so that $|y_0| \geq |y_1| \geq \max_{j \neq 0,1} |y_j|$ and without loss of generality write $y_0=1$.
Let $g' \in G$ be a multiplication of $y_1$ by a fourth root of unity so that $\arg(y_1) \in (-22.5^\circ, 67.5^\circ]$.
(Note this is possible because the length of this interval is $360^\circ/4=90^\circ$.)
Then $\arg(y_1) \in (-40^\circ,40^\circ)$ or $\arg(y_1) \in (5^\circ,85^\circ)$, so $g'g y \in U^*$ or $g'g y \in U^{**}$.
Equivalently, $y \in (g'g)^{-1} U^*$ or $y \in (g'g)^{-1} U^{**}$, so $\mathcal{A}$ covers all of $M$.
\end{proof}

Each chart of the atlas is uniquely identified by the three numbers:
index of biggest modulus coordinate, index of second biggest modulus coordinate, centre of argument angle window.
Define the chart assignment $a:M \rightarrow \mathcal{A}$ by assigning $p \in M$ to the following chart:
rescale the largest-modulus entry of $p$ to be $1$.
Then choose a chart which has the (ordered) two largest modulus entries with the same indices as $p$, and we choose the argument angle window as follows:
we partition $[0,2\pi]$ modulo addition of $2\pi$ into the eight intervals $[-22.5^\circ,22.5^\circ)$, $ [22.5^\circ,67.5^\circ)$,\,\dots, $[-67.5^\circ, -22.5^\circ)$.
We then assign $p$ to $U^*$ and $U^{**}$ and their translates in the obvious way, i.e. if the argument of the second-largest modulus entry of $p$ is in $[-22.5^\circ,22.5^\circ)$, then it is assigned to a translate of $U^*$;
if the argument is in $[22.5^\circ,67.5^\circ)$, then it is assigned to a translate to $U^{**}$ etc.
Where ties occur, choose some fixed assignment to one of the tied possible options.
The exact choice does not matter for the later analysis.

\textbf{Bounding $\rho$ for the chosen atlas.}
In the following, we use the notation $(U_a,\psi_a)_{a\in\mathcal A}$ for the holomorphic charts $\psi_a:U_a\to V_a\subset\mathbb C^m$ of this parent atlas.
We first bound how far in the interior any point $p \in M$ is for its assigned chart $a(p)$.
More precisely, we bound:
\begin{equation}
 \operatorname{dist}\bigl(\psi_{a(p)}(p),\partial V_{a(p)}\bigr)\geq\rho
 \qquad\text{for every }p\in M.
 \label{equation:parent-chart-margin}
\end{equation}

We have:

\begin{proposition}
\label{prop:atlas-rho-bound}
    For the chosen atlas, $\rho>1/200$.
\end{proposition}

\begin{proof}
    First suppose that $p=[1:X_1:X_2:X_3]$ is assigned to $U^*$, so that $|X_2|,|X_3| \leq |X_1| \leq 1$ and $\arg(X_1) \in [-22.5^\circ,22.5^\circ]$.
    We show that the ball of radius $1/200$ around $\varphi^*(p)=(X_2,X_3)$ is contained in $\varphi^*(U^*)$.
    To this end, let $|(Y_2,Y_3)-(X_2,X_3)| \leq 1/200$.
    
    Because $1=|X_1^4+X_2^4+X_3^4| \leq 3 |X_1|^4$ and hence $|X_1|\geq 3^{-1/4}>\frac{3}{4}$, we have that
    \[
    \left| \frac{-1-Y_2^4-Y_3^4}{X_1^4}-1 \right|
    \leq
    \frac{|Y_2^4-X_2^4|+|Y_3^4-X_3^4|}{|X_1|^4}
    \leq
    \frac{4 \sqrt{2}}{200 |X_1|} \left( 1+ \frac{1}{200|X_1|} \right)^3 < \frac{1}{25},
    \]
    where in the first step we used $X_1^4=-1-X_2^4-X_3^4$ and in the second step we used
    $|a^4-b^4|=|a-b| \cdot |a^3+a^2b+ab^2+b^3| \leq 4 \max \{|a|,|b|\}^3 |a-b|$ for $a,b \in \CC$ together with the Cauchy-Schwarz inequality.
    
    We may therefore define $Y_1 := X_1 \left( \frac{-1-Y_2^4-Y_3^4}{X_1^4} \right)^{1/4}$ choosing the fourth root near $1$.
    Write $w:=Y_1/X_1$.
    The chosen fourth root satisfies $\Re w, \Re(w^2)>0$, so $|w+1|,|w^2+1|>1$.
    Thus
    \begin{align}
    \label{eqn:w4-estimate}
    |w-1|
    =
    \frac{|w^2-1|}{|w+1|}
    \leq
    |w^2-1|
    =
    \frac{|w^4-1|}{|w^2+1|}
    \leq
    |w^4-1|
    <\frac{1}{25}.
    \end{align}
    Therefore
    \[
    |(w+1)(w^2+1)|
    =
    |2+(w-1)| |2+(w^2-1)|
    \geq
    (2-|w-1|)(2-|w^2-1|)>\left( \frac{49}{25} \right)^2>\frac{96}{25}.
    \]
    Furthermore
    \[
    \begin{aligned}
    |X_1||w^4-1|
    &=
    \frac{|(Y_2^4-X_2^4)+(Y_3^4-X_3^4)|}{|X_1|^3}\\
    &\leq
    \frac{
    |Y_2-X_2|\,|Y_2^3+Y_2^2X_2+Y_2X_2^2+X_2^3|
    +|Y_3-X_3|\,|Y_3^3+Y_3^2X_3+Y_3X_3^2+X_3^3|
    }{|X_1|^3}\\
    &\leq
    \frac{4(|X_1|+1/200)^3}{|X_1|^3}
    \bigl(|Y_2-X_2|+|Y_3-X_3|\bigr)\\
    &\leq
    \frac{\sqrt2}{50}
    \left(1+\frac{1}{200|X_1|}\right)^3\\
    &<
    \frac3{100},
    \end{aligned}
    \]
    where in the second step we used the triangle inequality;
    in the third step we used
    $\max\{|X_j|,|Y_j|\}\leq |X_1|+1/200$ for $j=2,3$;
    in the fourth step we used the norm equivalence of the $1$-norm and Euclidean norm and the smallness assumption from the beginning,
    $|Y_2-X_2|+|Y_3-X_3|
    \leq
    \sqrt2\,|(Y_2,Y_3)-(X_2,X_3)|
    \leq
    \frac{\sqrt2}{200}$;
    and in the fifth step we used $|X_1|>3/4$.
        Therefore
    \[
    |Y_1-X_1|
    =
    |X_1(w-1)|
    =
    |X_1| \left| \frac{w^4-1}{(w+1)(w^2+1)} \right|
    <
    \frac{3/100}{96/25}=\frac{1}{128}.
    \]
    Hence,
    \begin{align*}
    \max \{|Y_1|,|Y_2|,|Y_3|\}
    &< \frac{129}{128} < \frac{100}{99},
    \\
    |Y_1|-\frac{49}{50}\max \{|Y_2|,|Y_3|\}
    &>
    |X_1|-\frac{1}{128}
    -\frac{49}{50}\left(|X_1|+\frac{1}{200}\right)
    \\
    &=
    \frac{|X_1|}{50}-\frac{1}{128}-\frac{49}{10000}
    \\
    &>
    \frac{3}{200}-\frac{1}{128}-\frac{49}{10000}>0.
    \end{align*}
    Moreover, $|w^4-1|<1/25$ implies
    $\Re(w^4)>24/25$ and $|\Im(w^4)|<1/25$.
    Thus, for the chosen fourth root,
    \[
    |\arg(Y_1/X_1)|
    =
    \frac14\arctan\left(\frac{|\Im(w^4)|}{\Re(w^4)}\right)
    <
    \frac14\arctan(1/24)
    <
    \frac1{96}.
    \]
    Since $22.5^\circ=\pi/8$ and $40^\circ=2\pi/9$, we have
    \[
    \arg(Y_1)
    \in
    \left(-\frac{\pi}{8}-\frac{1}{96},
           \frac{\pi}{8}+\frac{1}{96}\right)
    \subset
    \left(-\frac{2\pi}{9},\frac{2\pi}{9}\right).
    \]
    Hence, $(Y_2,Y_3) \in \varphi^*(U^*)$.
    The same argument applies to $U^{**}$, and taking $G$-translates proves the claim.
\end{proof}

\textbf{Bounding pointwise quantities.}
Write $G_a$ for the matrix of $g$ in the real coordinates underlying $\psi_a$.
We need uniform lower and upper bounds for the eigenvalues of $G_a$:
\begin{equation}
 \lambda I\leq G_a\leq\Lambda I,
 \label{equation:eigenvalues-metric}
\end{equation}
throughout $V_a$, for every $a\in\mathcal A$, where $0<\lambda\leq\Lambda$.

The small ball charts will later be constructed to be isometries at one interior anchor point.
To this end, for $a\in\mathcal A$ and $x\in U_a$ define the normalising matrix to be the uniquely defined symmetric, positive definite root $Q_{a,x}:=G_a(\psi_a(x))^{1/2}$, which then satisfies $Q_{a,x}G_a(\psi_a(x))^{-1}Q_{a,x}^{T}=I$.
Since $g$ is Hermitian, $Q_{a,x}$ is complex linear.

We also require
\begin{equation}
 \left|
 DG_a\bigl(\psi_a(x)\bigr)
 \left[
 Q_{a,x}^{-1}\,\cdot,
 Q_{a,x}^{-1}\,\cdot,
 Q_{a,x}^{-1}\,\cdot
 \right]
 \right|
 \leq\widetilde K
 \label{equation:3rd-der-potential}
\end{equation}
for every $a\in\mathcal A$ and $x\in U_a$.

Write $F_{ab}:=\psi_a\circ\psi_b^{-1}$ for the transition between parent charts.
Furthermore, we will need bounds
\begin{align}
 \left|Q_{a,x}D^2F_{ab}\bigl(\psi_b(x)\bigr)
 \left[Q_{b,x}^{-1}\,\cdot,Q_{b,x}^{-1}\,\cdot\right]\right|
 &\leq\tilde c_1,
 \label{equation:parent-transition-second-derivative-bound}\\
 \left|Q_{a,x}D^3F_{ab}\bigl(\psi_b(x)\bigr)
 \left[
 Q_{b,x}^{-1}\,\cdot,
 Q_{b,x}^{-1}\,\cdot,
 Q_{b,x}^{-1}\,\cdot
 \right]\right|
 &\leq\tilde c_2
 \label{equation:parent-transition-third-derivative-bound}
\end{align}
for every $a,b\in\mathcal A$ and $x\in U_a\cap U_b$.
All derivatives and tensor norms are real.
Both bounds use the metric normalisation at the same point $x$.

\begin{proposition}
    The estimates \cref{equation:eigenvalues-metric,equation:3rd-der-potential,equation:parent-transition-second-derivative-bound,equation:parent-transition-third-derivative-bound} hold with
        \begin{align*}
        \lambda &= \num{\numLambdaLower}
        \\
        \Lambda &= \num{\numLambdaUpper}
        \\
        \tilde K &= \num{\numKtilde}
        \\
        \tilde c_1 &= \num{\numCOneTilde}
        \\
        \tilde c_2 &= \num{\numCTwoTilde}.
    \end{align*}
\end{proposition}

\begin{proof}
    We bound this using computer assistance using the method described in \cref{section:rigorous-error-bounds}.
\end{proof}

\subsection{Diameter bound}
\label{subsection:diameter-bound}

The automorphism group of $(M,\omega_0)$ is the group
\begin{align}
    G := \bigl(\mu_4^4/\mu_{4,\mathrm{diag}}\bigr)\rtimes S_4
\end{align}
of order $|G|=1536$, see \cite[Remark 3.4]{Oguiso2005}.
Here $S_4$ acts by permutations, and $\mu_4^4$ acts by multiplications of the four coordinates by fourth roots of unity, and $\mu_{4,\mathrm{diag}}$ acts by multiplying each coordinate by the same fourth root of unity.
In this subsection, we will obtain \emph{upper bounds} for $\diam(M,\omega)$ and $\diam(M/G,\omega)$, where we take the quotient metric on $M/G$, i.e. $d(p,q)=\min_{g \in G} d(p,gq)$.
We normalise $\omega_0$ by $\int_M\omega_0^2/2=1$ and assume that $G$ preserves $\omega$.
Our approach is to prove diameter bounds for the metric $g_0$ induced by the ambient Fubini-Study metric and then simply multiply these bounds with one global comparison factor.

\begin{proposition}
    \label{proposition:fubini-study-diameter-bound}
    We have
    \[
    \diam(M,\omega_0) \leq \frac{2\Cr{const:connection-to-E-length}+4\pi/(3\sqrt3)}{\sqrt{\pi\sqrt2}},
    \qquad
    \diam(M/G,\omega_0) \leq 2\Cr{const:connection-to-E-length}/\sqrt{\pi\sqrt2},
    \]
    where $\Cr{const:connection-to-E-length}$ is defined in \cref{eqn:connection-to-E-length-def}.
\end{proposition}

\begin{proof}
    Write $\zeta=e^{\pi i/4}$ and define $E:=G\cdot[1:\zeta:0:0]
        =\{[e_i+\xi e_j]:0\leq i<j\leq3,\ \xi^4=-1\}$.
    For any $p\in M$, there exists $g\in G$ such that
    \[
        gp=[1:X_1:X_2:X_3],
        \qquad |X_j|\leq1,
        \qquad \operatorname{Re}(X_1^4)\leq-\frac13.
    \]
    Define the absolutely continuous path
    $\hat{\gamma} = [1:\gamma]=[1:\gamma_1:\gamma_2:\gamma_3]:[0,1] \rightarrow M$ by
    \[
        \gamma_2(t)=(1-t)^{1/4}X_2,\qquad
        \gamma_3(t)=(1-t)^{1/4}X_3,\qquad
        \gamma_1(t)^4=(1-t)X_1^4-t,
        \qquad \gamma_1(0)=X_1.
    \]
    The initial condition determines a unique continuous fourth root,
    since its fourth power has negative real part.
    Thus $\hat{\gamma}$ joins $gp$ to $E$ in $M$.
    
   We set $w=(\gamma_2,\gamma_3)$, $r=|w|$ and $A=1 + |\gamma_1|^2$. 
   We have
   \begin{align*}
       g_{\operatorname{FS}}(\hat{\gamma}',\hat{\gamma}') = \frac{(1+|\gamma|^2) |\gamma'|^2 - |\langle \gamma,\gamma'\rangle|^2}{(1+|\gamma|^2)^2}
   \end{align*}
   If for $z\in \mathbb{C}^3$ we introduce the Hermitian form 
   \begin{align*}
       b_z(\xi) = \frac{(1+|z|^2) |\xi|^2 - |\langle z,\xi\rangle|^2}{(1+|z|^2)^2},
   \end{align*}
   using the triangle inequality for the induced norms gives 
   \begin{align*}
       \sqrt{g_{\operatorname{FS}}(\hat{\gamma}',\hat{\gamma}')} \leq \sqrt{b_{\gamma}(\gamma_1',0,0)} + \sqrt{b_{\gamma}(0,w')}
       \leq 
       \frac{|\gamma_1'|}{1+|\gamma_1|^2} + \frac{\sqrt{A}}{A+r^2} |r'|,
       \end{align*}
   where we have used that $|\langle w, w'\rangle|^2 = |w|^2 |w'|^2$ and $|w'| = |r'|$ since $w$ and $w'$ are parallel.

   Now we claim
   \begin{align*}
        |\gamma_1'|
        =\frac{|1+X_1^4|}{4|\gamma_1|^3},
        \qquad
        |1+X_1^4|\leq\frac2{\sqrt3},
        \qquad
        |\gamma_1|^4\geq\frac{1+2t}{3},
        \qquad
        \frac{\sqrt A}{A+r^2}
        \leq
        \begin{cases}
            \frac12\sqrt{2-r^2},&r\leq1,\\[2pt]
            \frac1{2r},&r\geq1.
        \end{cases}
       \end{align*}
       The equality is a direct calculation.
    The first and second inequalities follow by using $\operatorname{Re}(X_1^4)\leq -\frac{1}{3}$. We have
    \begin{align*}
        |1+X_1^4|^2 = 1 + 2\operatorname{Re}(X_1^4) + |X_1^4|^2 \leq \frac{4}{3}
    \end{align*}
    and
    \begin{align*}
        \operatorname{Re}((1-t)X_1^4-t)\leq -\frac{1}{3}(1-t)-t = - \frac{1+2t}{3}
    \end{align*}
    and thus
    \begin{align*}
        |\gamma_1(t)|^4 \geq |\operatorname{Re}((1-t)X_1^4-t)| \geq \frac{1+2t}{3}.
    \end{align*}
    To explain the third inequality, we first note that $A+r^2 \geq 2$. This follows from $\sum_j \operatorname{Re}(\gamma_j(t)^4) =-1$.  Therefore, if $r\leq 1$ we have $A\geq r^2$ and the function $f(B)=\frac{\sqrt{B}}{B+r^2}$ is monotone decreasing for $B\geq r^2$. Therefore, $f(A) \leq f(2-r^2)$ which gives the case $r\leq 1$. For $r\geq 1$ we use the general inequality $A+r^2 \geq 2 r \sqrt{A}$. 
    Since $r$ decreases from at most $\sqrt{2}$ to $0$, using the above inequalities gives
     \begin{align}
        \sqrt{\pi\sqrt2}\,
        \operatorname{length}_{g_0}(\gamma)
        &=
        \int_0^1
        \sqrt{g_{\operatorname{FS}}(\gamma',\gamma')}\,dt\notag\\
        &\leq
        \int_0^1
        \left(
            \frac{|\dot\gamma_1|}{1+|\gamma_1|^2}
            +
            \frac{\sqrt A}{A+r^2}|r'|
        \right)\,dt\notag\\
        &\leq
        \frac1{2\sqrt3}
        \int_0^1
        \frac{dt}{
            \left(\frac{1+2t}{3}\right)^{3/4}
            \left(1+\sqrt{\frac{1+2t}{3}}\right)
        }
        +
        \frac12\int_0^1\sqrt{2-r^2}\,dr
        +
        \frac12\int_1^{\sqrt2}\frac{dr}{r}\notag\\
        &=
        \sqrt3\int_{3^{-1/4}}^1\frac{du}{1+u^2}
        +\frac14+\frac{\pi}{8}+\frac{\log2}{4}
        \notag\\
        &=
        \sqrt3\left(
            \frac{\pi}{4}-\arctan(3^{-1/4})
        \right)
        +\frac14+\frac{\pi}{8}+\frac{\log2}{4}
        =:\Cl{const:connection-to-E-length}.
        \label{eqn:connection-to-E-length-def}
    \end{align}
    Thus, every point in $M$ can be joined to a point of $E$ through a curve of length at most $\Cr{const:connection-to-E-length}/\sqrt{\pi\sqrt2}$.
    Hence,
    $\operatorname{diam}(M/G,\omega_0)
        \leq 2\Cr{const:connection-to-E-length}/\sqrt{\pi\sqrt2}$.

    It remains to bound $\diam (M,\omega_0)$.
    To this end, we bound the distance between two arbitrary points in $E$.
    
    Let $i,j,k\in\{0,1,2,3\}$ be distinct and let $\xi,\eta$
    be fourth roots of $-1$.
    We first connect $[e_i+\xi e_j]$ and $[e_i+\eta e_k]$.
    Define $\gamma:[0,1]\to M$ by
    \[
        \gamma(t)=[e_i+z_1(t)e_j+z_2(t)e_k],
    \]
    where
    \[
        a(t)=\sqrt{
            \cos^4\left(\frac{\pi t}{2}\right)
            +\sin^4\left(\frac{\pi t}{2}\right)
        },
        \qquad
        z(t)=\frac{
            \left(
                \xi\cos\left(\frac{\pi t}{2}\right),
                \eta\sin\left(\frac{\pi t}{2}\right)
            \right)
        }{\sqrt{a(t)}}.
    \]
    This curve connects the two points. We bound its length by
    \begin{align*}
        \sqrt{\pi\sqrt2}\,
        \operatorname{length}_{g_0}(\gamma)
        &=
        \int_0^1
        \frac{
            \sqrt{
                (1+|z|^2)|z'|^2
                -|\langle z,z'\rangle|^2
            }
        }{1+|z|^2}\,dt\\
        &=
        \int_0^1
        \sqrt{
            \frac{\pi^2}{4(1+a)}
            +\frac{(a')^2}{4a(1+a)^2}
        }\,dt\\
        &=
        \frac{\pi}{2}\int_0^1
        \sqrt{\frac{2a-1}{2a^3}}\,dt\\
        &\leq
        \frac{\pi}{2}\int_0^1\frac4{3\sqrt3}\,dt
        =\frac{2\pi}{3\sqrt3}.
    \end{align*}
    In the third line we have used
    \[
        (a')^2=\frac{\pi^2}{a^2} \frac{1-a^2}{2} (2a^2-1).
    \]
    Indeed, if we set $c(t) = \cos(\frac{\pi t}{2})$ and $s(t)=\sin(\frac{\pi t}{2})$, then
    \begin{align*}
        (a')^2 = \frac{\pi^2}{a^2} c^2s^2 (s^2-c^2)^2
    \end{align*}
    and
    \begin{align*}
        2c^2s^2 = (c^2+s^2)^2-(c^4+s^4)=1-a^2,
        \\
        (s^2-c^2)^2 = (s^2+c^2)^2-4c^2s^2 = 1-2(1-a^2) = 2a^2-1.
    \end{align*}
    In the fourth line we have used
    \[
        \frac{2a-1}{2a^3}\leq\frac{16}{27}
        \qquad\text{for}\qquad
        \frac1{\sqrt2}\leq a\leq1.
    \]
    This connects two points in $E$ if they have exactly one common non-zero coordinate index.
    Points $[e_i+\xi e_j]$ and $[e_k+\zeta e_l] \in E$ for $i,j,k,l$ pairwise distinct and points with the same non-zero indices can be connected via a third point.
    Hence,
    \[
        \diam(M,\omega_0)
        \leq
        \frac{2\Cr{const:connection-to-E-length}+4\pi/(3\sqrt3)}{\sqrt{\pi\sqrt2}}.
    \]
\end{proof}

\begin{remark}
    The general-purpose estimate for Fubini-Study diameters of algebraic varieties from \cite{FengSchumacher1999} gives a diameter bound of $\operatorname{diam}(M,\omega_0)
    \leq 1625\pi/\sqrt{\pi\sqrt2}
    <2422$.
    The bound from \cref{proposition:fubini-study-diameter-bound} is sharper and also the reference does not readily provide bounds for $\diam(M/G)$, other than the trivial $\diam(M/G) \leq \diam(M)$.
\end{remark}

The previous result gives the diameter with respect to the Fubini-Study metric.
To obtain the diameter with respect to the approximate Calabi-Yau metric $\omega$, we use the following certified comparison:

\begin{proposition}
    The approximate Calabi-Yau metric $g$ satisfies 
    \[
    g \leq \beta g_0
    \]
    as quadratic forms,
    where
    $\beta=\num{\numMetricComparison}$.
\end{proposition}

\begin{proof}
    We bound this using computer assistance using the method described in \cref{section:rigorous-error-bounds}.
\end{proof}

Combining this comparison with
\cref{proposition:fubini-study-diameter-bound} gives
\[
    \diam(M,\omega)\leq\num{\numUpperBoundDiam},
    \qquad
    \diam(M/G,\omega)\leq\num{\numUpperBoundQuotientDiam}.
\]

\subsection{Ricci bound}
\label{subsection:ricci-bound}

In \cref{subsection:certified-errors-bounds}, we rigorously bounded the Monge–Ampère error of our approximate solution.
This was the optimisation objective, so it is unsurprising this number was indeed certified to be small.

The Ricci tensor is, roughly speaking, the second derivative of the Monge–Ampère error.
Making it small was not an objective during optimisation.
However, to have any hope to be near the true Calabi-Yau metric, one expects the Ricci tensor of this approximation to be small.
Indeed, this (or at least small scalar curvature) has been empirically observed several times, e.g. \cite[Table 2]{MirjanicMishra2025}, \cite[Figure 1]{ButbaiaEtAl2025}, \cite[Figure 1]{LarforsLukasRuehleSchneider2022}.
In the following proposition, we promote this empirical observation to a certified bound:

\begin{proposition}
    We have
    \begin{align}
    \label{eqn:ricci-bound}
    \Ric_g > -(n-1) \kappa g = -3 \kappa g
    \end{align}
    for $\kappa=\num{\numRicciBound}$.
\end{proposition}

\begin{proof}
    We bound this using computer assistance using the method described in \cref{section:rigorous-error-bounds}.
\end{proof}

\subsection{Spectral gap bound}
\label{subsection:spectral-gap}

Using the upper diameter bound from \cref{subsection:diameter-bound} and the Ricci lower bound from \cref{subsection:ricci-bound}, one can compute a bound for the spectral gap of the Laplace operator acting on invariant functions.
There are many different routes to obtaining such a bound.
We use the following estimate for $G$-invariant functions, which is sufficient for our application and easy to compute rigorously.

\begin{theorem}
    Let $(X,g)$ be a closed connected $n$-dimensional Riemannian manifold,
    $n\geq 2$, and $G \subset \operatorname{Isom}(X,g)$ be a finite subgroup.
    If $\operatorname{Ric}_g \geq -(n-1)\kappa g$ and
    $\operatorname{diam}(X/G) \leq D$, where $\kappa\geq 0$, $D>0$ and
    $(n-1)\kappa D^2<4\pi^2$, then the smallest non-zero eigenvalue
    $\lambda_{1,G}$ of the Riemannian Laplacian
    $-\Delta$ on $X$, restricted to $G$-invariant
    functions, satisfies
    \[
        \lambda_{1,G}
        \geq \left(\frac{\pi}{D}
        -\frac{(n-1)\kappa D}{4\pi}\right)^2.
    \]
\end{theorem}

\begin{proof}
    By \cite[Corollary~1.2]{GGKMS2018}, $X/G$ is an $\mathrm{RCD}^{*}(-(n-1)\kappa,n)$ space.
By \cite[Theorem~5.12]{GGKMS2018}, $\lambda_{1,G}=\lambda_1(X/G)$.
By \cite[Theorem~1.2]{JiangZhang2016}, $\lambda_1(X/G) \geq \mu_{n,\kappa}(D)$, 
where $\mu$ is the first positive
eigenvalue of the Neumann problem on $(-D/2,D/2)$:
\[
    u''+(n-1)\sqrt{\kappa}\tanh(\sqrt{\kappa}\,t)u'
    =-\mu u,
    \qquad u'(\pm D/2)=0.
\]
The computation in \cite[proof of Corollary~4.3]{QianZhangZhu2013},
with $K=-\kappa$, gives
\[
    \mu\geq4s(1-s)\frac{\pi^2}{D^2}-s(n-1)\kappa
    \qquad\text{for all }s\in(0,1).
\]
Taking
\[
    s=\frac12-\frac{(n-1)\kappa D^2}{8\pi^2}\in(0,1)
\]
proves the claim.
\end{proof}

\begin{corollary}
The $G$-invariant spectral gap for the Riemannian Laplacian on $(M,\omega)$ satisfies
\begin{align}
\label{eqn:spectral-gap}
    \lambda_{1,G}\geq\num{\numLambdaOne}.
\end{align}
\end{corollary}

\section{Perturbation theory of the complex Monge--Ampere equation}
\label{section:perturbation-theory}

Let $(M,I)$ be a compact connected complex manifold of complex dimension $m\geq2$, and fix $0<\alpha<1$. Suppose that  the first Chern class of $M$ vanishes. Let $\Omega$ be a holomorphic volume form. Then $$\mu_{\Omega} = i^{m^2} \Omega\wedge\bar{\Omega}$$
is a positive real volume form on $M$. Let $\omega$ be a K\"ahler metric on $M$. Denote its volume form by $\mu_{\omega} = \frac{1}{m!} \omega^m$ and assume $\vol(M,\mu_\omega)=1$. To find a Ricci-flat K\"ahler metric, whose existence is provided by the proof of the Calabi conjecture, we need to solve the equation $\mu_{\omega} = \hat{\mu}$, where $\hat{\mu}$ is the normalized volume form
\begin{align}
    \label{eqn:volume-normalisation}
    \hat{\mu} = \frac{\vol(M, \mu_{\omega})}{\vol(M,\mu_{\Omega})} \mu_{\Omega}.
\end{align}
Here we write $\vol(M,\mu) = \int_M \mu$ for the volume of $M$ with respect to a volume form $\mu$. Following Donaldson, we write
\begin{align*}
    \eta = \frac{\mu_{\omega}}{\hat{\mu}}.
\end{align*}
With this equation, the Calabi--Yau equation then becomes $\eta \equiv 1$.

To find a Ricci flat K\"ahler metric $\tilde{\omega}= \omega + i\partial\bar{\partial}u$ in the de Rham cohomology class of $\omega$, we need to solve the complex Monge--Ampere equation
\begin{align*}
    \frac{(\omega+i\partial\bar{\partial}u)^m}{\omega^m}
    =
    \frac{\hat{\mu}}{\mu_{\omega}}.
\end{align*}
In this section we will discuss how to solve this equation perturbatively if $\omega$ is a good approximation to a Calabi--Yau metric.
To set up the problem, introduce the operator 
\begin{align}
    \label{equation:monge-ampere}
    P(u)=\frac{(\omega+i\partial\bar{\partial}u)^m}{\omega^m}- \frac{\hat{\mu}}{\mu_{\omega}},
\end{align}
and use 
the function $f$ 
\begin{align}
    \label{eqn:residual-def}
    f = 1-\frac{\hat{\mu}}{\mu_{\omega}}=1-\frac{1}{\eta}.
\end{align}
to measure the failure of $\omega$ to be Ricci-flat.
Thus  we want to solve $P(u)=0$ given that $P(0)=f$ is small in a certain sense. We want to apply a variant of the inverse function theorem in Banach spaces. 
The linearisation of $P$ at $u=0$ is $\frac{1}{2}\Delta_g$ and therefore we can expand $P$ as
\begin{align*}
    P(u)= f + \frac{1}{2}\Delta_g u + Q(u),
\end{align*}
where the operator 
\begin{align}
    Q(u) 
    =
    \sum_{k=2}^m \binom{m}{k} \frac{(i\partial\bar{\partial}u)^k\wedge \omega^{m-k}}{\omega^m}.
    \label{non-linear-part}
\end{align}
is the nonlinear part of $P$.

To solve the perturbation problem we will apply the following variant of the Banach fixed point theorem:

\begin{lemma}
\label{lemma:inverse-function-thm}
    Let $P:\mathscr{X}\xrightarrow{}\mathscr{Y}$ be a smooth operator between Banach spaces and write $P(u)=P(0)+L(u)+Q(u)$, where $L$ is linear and $Q$ is the nonlinear part of $P$. Assume there exist constants $C,t,q>0$ such that
    \begin{itemize}
        \item $L$ is invertible with $\|{L^{-1}}\leq C$,
        \item $\|{Q(u)-Q(v)}_{\mathscr{Y}}\leq q \|{u-v}_{\mathscr{X}} (\|{u}_{\mathscr{X}} + \|{v}_{\mathscr{X}})$ for all 
             $u$ and $v$ with $\|{u}_{\mathscr{X}}, \|{v}_{\mathscr{X}} < t$,
        \item $\|{P(0)}_{\mathscr{Y}} < \min\{\tfrac{t}{2C},\tfrac{1}{4q C^2} \}$.
    \end{itemize}
    Then there exists a unique $u\in \mathscr{X}$ with $\|{u}_{\mathscr{X}} \leq 2C \|{P(0)}_{\mathscr{Y}}$
    such that $P(u)=0$.
\end{lemma}

Lemma \ref{lemma:inverse-function-thm} is proved by applying the Banach fixed point theorem to the map $T(u)=-L^{-1}(P(0)+Q(u))$. The solution is the limit of the iteration $u_{k+1} = T(u_k), u_0 = 0$.

\subsection{Conventions}

Throughout, let $n=2m$ be the real dimension of $M$.
By $\omega_{n-1}=|S^{n-1}|$ we denote the volume of the $(n-1)$-dimensional unit sphere.

Next we define Hölder spaces on $B_R$, the ball of radius $R$ in $\mathbb{R}^n$. We will work with two conventions for scale invariant Hölder spaces for functions on $B_R$. We define
\begin{align*}
    \begin{split}
    \|{u}'_{C^{k,\alpha}(B_R)}
    &= 
    \sum_{l=0}^k R^l [u]_{l,0;B_R} + R^{k+\alpha} [u]_{k,\alpha;B_R},
    \\
    \|{u}'_{C^{k}(B_R)}
    &
    =
    \sum_{l=0}^k R^l [u]_{l,0;B_R},
    \end{split}
\end{align*}
where 
\begin{align*}
    [u]_{k,0;\Omega} &=  \sup_{x\in\Omega} |D^k u(x)|,
    \\
    [u]_{k,\alpha;\Omega} &=  \sup_{x\neq y\in\Omega} 
    \frac{|D^k u(x)-D^k u(y)|}{|x-y|^{\alpha}}.
\end{align*}
We will frequently work on a ball of radius $2R$. To avoid factors of $2$ we will use a different convention for the weight factors in the definition of the Hölder norm:
\begin{align*}
    \begin{split}
    \|{u}''_{C^{k,\alpha}(B_{2R})}
    &= 
    \sum_{l=0}^k R^l [u]_{l,0;B_{2R}} + R^{k+\alpha} [u]_{k,\alpha;B_{2R}},
    \\
    \|{u}''_{C^{k}(B_{2R})}
    &
    =
    \sum_{l=0}^k R^l [u]_{l,0;B_{2R}}.
    \end{split}
\end{align*}

Throughout, we use the Euclidean tensor norm $|T|^2=\sum_I |T_I|^2$ and write $|A|_{\mathrm{op}}$ for the operator norm.
If $\gamma$ is a differential form of degree $k$, by $|\gamma|$
we denote the norm of $\gamma$ as a differential form, which is its tensor norm divided by $\sqrt{k!}$.

Single bars denote pointwise norms and double bars denote function-space norms.

\subsection{Construction of ball charts}
\label{section:analytic-set-up}

In this section we will for each $p\in M$ construct a 
holomorphic chart $\phi_p: \mathcal{B}_p \rightarrow B_{2R}$ 
such that if we denote by $(g_{\mathcal B_p})_{ij}$ the metric tensor in these charts and by $g_{\mathcal B_p}^{ij}$ its inverse, we have $(g_{\mathcal B_p})_{ij}(0)= \delta_{ij}$ and $g_{\mathcal B_p}$ is close to the Euclidean metric. This will allow us to prove a Schauder estimate for $\Delta_g$ by coefficient freezing from the local Schauder interior estimate in Euclidean space.

The ball charts will be constructed from the atlas $(U_a,\psi_a)_{a\in\mathcal A}$ of holomorphic charts $\psi_a:U_a\to V_a\subset\mathbb C^m$ introduced in \cref{subsection:atlas-bounds}. We recall that 
we are also given an assignment $a:M\rightarrow\mathcal A$ that assigns to each point in $M$ a parent chart such that $p\in U_{a(p)}$ and a constant $\rho>0$ from \cref{prop:atlas-rho-bound} which measures how far inside a point is in its assigned chart. 
Furthermore, in \cref{subsection:atlas-bounds} we introduced a square root $Q_{a,x}$ of $G_a(\psi_a(x))$,  bounds $\lambda, \Lambda$ from \cref{equation:eigenvalues-metric} measuring how degenerate the metric tensor $G_a$ is in coordinates, and
the bound $\tilde K$ from  \cref{equation:3rd-der-potential} measuring the derivative of the metric tensor in coordinates.

The idea to construct the ball chart around $p\in M$ is to compose the chart $\psi_{a(p)}$ with a complex affine transformation such that $p$ gets mapped to $0$ and $(g_{\mathcal B_p})_{ij}(0)= \delta_{ij}$.

\begin{definition}
\label[definition]{definition:ball-charts}
Let $0<R$ such that
\begin{align}
    R\leq R_0 :=\frac12\rho\sqrt{\lambda}. \label{radius-cond-1}
\end{align}
For every $p\in M$, write $a=a(p)$ and define
\begin{align*}
 \mathcal B_p
 &:=\{x\in U_a:|Q_{a,p}(\psi_a(x)-\psi_a(p))|<2R\},\\
 \phi_p:\mathcal B_p
 &\to B_{2R},
 \qquad x\mapsto Q_{a,p}(\psi_a(x)-\psi_a(p)).
\end{align*}
\end{definition}

\begin{proposition}
\label[proposition]{proposition:ball-chart-properties}
The charts of \cref{definition:ball-charts} have the following properties:
\begin{itemize}
 \item For every $p\in M$, with $a=a(p)$, the ellipsoid $\psi_a(p)+Q_{a,p}^{-1}B_{2R}$ is contained in $V_a$.
 \item Each $\phi_p:\mathcal B_p\to B_{2R}$ is biholomorphic.
 \item At the centre, $\phi_p(p)=0$ and $g_{\mathcal B_p}^{ij}(0)=\delta^{ij}$.
\end{itemize}
\end{proposition}

\begin{proof}
Fix $p\in M$ and write $a=a(p)$.
By \eqref{equation:eigenvalues-metric}, $|Q_{a,p}^{-1}|_{\mathrm{op}}\leq\lambda^{-1/2}$.
For $y\in B_{2R}$ we therefore have $|Q_{a,p}^{-1}y|<2R/\sqrt{\lambda}\leq\rho$.
\Cref{equation:parent-chart-margin} implies that $\psi_a(p)+Q_{a,p}^{-1}y\in V_a$.
Since $Q_{a,p}$ is complex linear and invertible, $\phi_p$ is biholomorphic.
By construction, $\phi_p(p)=0$, and the inverse metric at the centre satisfies
\begin{align*}
 \bigl(g_{\mathcal B_p}^{ij}(0)\bigr)_{i,j}
 =Q_{a,p}G_a(\psi_a(p))^{-1}Q_{a,p}^{T}
 =I.
\end{align*}
\end{proof}

Our motivation for the construction of the ball charts $\mathcal{B}_p$ is to obtain coordinates in which the metric is close to the Euclidean metric. We will now give a quantitative bound for the difference of $g_{\mathcal{B}_p}$ and the Euclidean metric. This bound becomes arbitrarily small if $R$ is sufficiently small.

\begin{lemma}
\label[lemma]{lemma:metric-bounds-ball-charts}
If in addition to $R\leq R_0$ we have 
\begin{align}
    \widetilde K R<\sqrt2, \label{radius-cond-2}
\end{align}
then 
\begin{align}
 \sup_{p\in M}
 \|{g_{\mathcal B_p}^{-1}-I}_{C^0(B_{2R})}
 \leq\varepsilon_0,
 \quad
 R^\alpha\sup_{p\in M}
 [g_{\mathcal B_p}^{-1}-I]_{0,\alpha;B_{2R}}
 \leq\varepsilon_\alpha,
 \label{equation:metric-coordinate-bounds}
\end{align}
where
\begin{align}
 \varepsilon_0
 =\sqrt2\left(
 \left(1-\frac{\widetilde K R}{\sqrt2}\right)^{-2}-1
 \right),
 \quad
 \varepsilon_\alpha
 =4^{1-\alpha}\widetilde K R
 \left(1-\frac{\widetilde K R}{\sqrt2}\right)^{-3}.
 \label{equation:ball-metric-constants}
\end{align}
\end{lemma}

\begin{proof}
Fix $p\in M$, write $a=a(p)$, and let $g_{\mathcal B_p}$ denote the covariant metric matrix as a function on $B_{2R}$.
In the ball coordinates we have
\begin{align}
 g_{\mathcal B_p}(y)
 =Q_{a,p}^{-T}
 G_a\bigl(\psi_a(p)+Q_{a,p}^{-1}y\bigr)
 Q_{a,p}^{-1}.
 \label{formula-gBp}
\end{align}
Put $\xi=\phi_p^{-1}(y)$ and $A=Q_{a,\xi}Q_{a,p}^{-1}$, so that $A^TA=g_{\mathcal B_p}(y)$.
Applying the chain rule to \eqref{formula-gBp}, using the identities $Dg_{\mathcal B_p}^{-1} = - g_{\mathcal B_p}^{-1} (Dg_{\mathcal B_p}) g_{\mathcal B_p}^{-1} $ and $Ag_{\mathcal B_p}(y)^{-1} = A^{-T}$, and using the bound \eqref{equation:3rd-der-potential}, we get
\begin{align}
\label{eqn:Dg-Dginv-bounds}
\begin{split}
 |Dg_{\mathcal B_p}(y)|
 &=\left|
DG_a\bigl(\psi_a(\xi)\bigr)
\left[
Q_{a,\xi}^{-1}A\,\cdot\, ,
Q_{a,\xi}^{-1}A\,\cdot\, ,
Q_{a,\xi}^{-1}A\,\cdot\,
\right]
\right|
\leq\widetilde K |A|_{\mathrm{op}}^3
 =\widetilde K |g_{\mathcal B_p}(y)|_{\mathrm{op}}^{3/2},
 \\
 |D(g_{\mathcal B_p}^{-1})(y)|
 &=\left|
DG_a\bigl(\psi_a(\xi)\bigr)
\left[
Q_{a,\xi}^{-1}A^{-T}\,\cdot\, ,
Q_{a,\xi}^{-1}A^{-T}\,\cdot\, ,
Q_{a,\xi}^{-1}A\,\cdot\,
\right]
\right|
\\
&
\leq\widetilde K
 |A^{-T}|_{\mathrm{op}}^2
 |A|_{\mathrm{op}}
 =\widetilde K
 |g_{\mathcal B_p}(y)^{-1}|_{\mathrm{op}}
 |g_{\mathcal B_p}(y)|_{\mathrm{op}}^{1/2}.
\end{split}
\end{align}
The eigenvalues of the real Hermitian matrices $g_{\mathcal B_p}^{\pm1}-I$ occur in pairs, so their operator norms are at most their tensor norms divided by $\sqrt2$.
For a unit vector $v$, set
\begin{align*}
 h(t):=\max\left\{
 |g_{\mathcal B_p}(tv)-I|,
 |g_{\mathcal B_p}(tv)^{-1}-I|
 \right\}.
\end{align*}
This function satisfies $h(0)=0$, and using \cref{eqn:Dg-Dginv-bounds} we find that it satisfies
\begin{align*}
 h'(t)\leq\widetilde K
 \left(1+\frac{h(t)}{\sqrt2}\right)^{3/2}.
\end{align*}
(Strictly speaking, we only know $h$ to be locally Lipschitz and this estimate to hold for almost every $0<t<2R$ because the tensor norm $|\cdot|$ is not smooth.)
Integrating
$\frac{d}{dt}(1+h(t)/\sqrt2)^{-1/2}
\geq-\widetilde K/(2\sqrt2)$ gives
\begin{align*}
 h(t)\leq\sqrt2\left[
 \left(1-\frac{\widetilde K t}{2\sqrt2}\right)^{-2}-1
 \right].
\end{align*}
Taking $t=|y|$ proves the first bound.
By the definition of $h$ and the fact that the operator norms of
$g_{\mathcal B_p}^{\pm1}-I$ are at most their tensor norms divided by
$\sqrt2$, we have
\begin{align}
 |g_{\mathcal B_p}(tv)^{\pm1}|_{\mathrm{op}}
 &\leq 1+\frac{h(t)}{\sqrt2}
 \leq
 \left(1-\frac{\widetilde Kt}{2\sqrt2}\right)^{-2}.
\end{align}
Plugging this into \cref{eqn:Dg-Dginv-bounds} gives
\begin{align*}
 |Dg_{\mathcal B_p}(tv)|
 &\leq\widetilde K
 \left[
 \left(1-\frac{\widetilde Kt}{2\sqrt2}\right)^{-2}
 \right]^{3/2}
 =\widetilde K
 \left(1-\frac{\widetilde Kt}{2\sqrt2}\right)^{-3}.
\end{align*}
Similarly, for the inverse we get
\begin{align*}
 |D(g_{\mathcal B_p}^{-1})(tv)|
 &\leq\widetilde K
 \left(1-\frac{\widetilde Kt}{2\sqrt2}\right)^{-2}
 \left(1-\frac{\widetilde Kt}{2\sqrt2}\right)^{-1}
 =\widetilde K
 \left(1-\frac{\widetilde Kt}{2\sqrt2}\right)^{-3}.
\end{align*}
Altogether, writing $y=tv$ with $t=|y|$, we get
\begin{align*}
 \max\left\{
 |Dg_{\mathcal B_p}(y)|,
 |D(g_{\mathcal B_p}^{-1})(y)|
 \right\}
 \leq\widetilde K
 \left(1-\frac{\widetilde K|y|}{2\sqrt2}\right)^{-3}.
\end{align*}

For $x,y\in B_{2R}$, integration along the segment gives
\begin{align*}
 |g_{\mathcal B_p}(x)^{-1}-g_{\mathcal B_p}(y)^{-1}|
 &\leq\widetilde K
 \left(1-\frac{\widetilde K R}{\sqrt2}\right)^{-3}|x-y|
 \leq
 \widetilde K
 \left(1-\frac{\widetilde K R}{\sqrt2}\right)^{-3}
 (4R)^{1-\alpha}|x-y|^\alpha.
\end{align*}
Taking the Hölder seminorm and multiplying by $R^\alpha$ proves the second bound.
\end{proof}
\begin{remark}
    For $A=Q_{a,\xi}Q_{a,p}^{-1}$ we have
    $|A|_{\operatorname{op}} \leq
    |Q_{a,\xi}|_{\operatorname{op}}
    |Q_{a,p}^{-1}|_{\operatorname{op}}
    \leq \bigl(\frac{\Lambda}{\lambda}\bigr)^{1/2}$
    and thus a simpler estimate than in the proof of
    Lemma \ref{lemma:metric-bounds-ball-charts} gives
    $|D g_{\mathcal{B}_p}| \leq
    \tilde{K} \bigl(\frac{\Lambda}{\lambda}\bigr)^{3/2}$ and
    \begin{align*}
        |g_{\mathcal{B}_p}(y) - I|
        = |g_{\mathcal{B}_p}(y) - g_{\mathcal{B}_p}(0)|
        \leq \tilde{K}
        \left(\frac{\Lambda}{\lambda}\right)^{3/2} |y|.
    \end{align*}
    Thus by making the extra assumption \eqref{radius-cond-2},
    with Lemma \ref{lemma:metric-bounds-ball-charts}
    we make the bounds $\varepsilon_0$ and $\varepsilon_{\alpha}$
    smaller as $R\to0$ by a factor of
    $\bigl(\frac{\Lambda}{\lambda}\bigr)^{3/2}
    \approx \num{\fpeval{round(
        (\numLambdaUpper/\numLambdaLower)^(3/2),0)}}$.
\end{remark}

With the bounds from Lemma \ref{lemma:metric-bounds-ball-charts} we get
\begin{lemma}\label{lemma:eigenvalues-metric}
We have
\begin{align*}
    K_0^{-1} \delta^{ij} \leq g^{ij}_{\mathcal{B}_p} \leq K_0 \delta^{ij}
\end{align*}
with
\begin{align*}
    K_0 = 1 + \frac{\varepsilon_0}{\sqrt2}.
\end{align*}
\end{lemma}

\begin{proof}
The preceding proof bounds the operator norms of both
$g_{\mathcal B_p}$ and $g_{\mathcal B_p}^{-1}$ by
$1+\varepsilon_0/\sqrt2$.
Inverting the former bound gives the lower bound.
\end{proof}

\begin{lemma}
\label{lemma:Holder-bound-omega}
    We have
    \begin{align*}
        R^{\alpha} [\omega_{\mathcal{B}_p}]_{0,\alpha;B_R} \leq S,
    \end{align*}
    where
    \begin{align*}
        S = 2^{\frac{1}{2}-\alpha}\,\widetilde K R
        \left(1-\frac{\widetilde K R}{2\sqrt2}\right)^{-3}.
    \end{align*}
    \end{lemma}

\begin{proof}
For the induced Euclidean norm on two-forms,
$|D\omega_{\mathcal B_p}|=|Dg_{\mathcal B_p}|/\sqrt2$.
The derivative estimate in the proof of
Lemma \ref{lemma:metric-bounds-ball-charts}, restricted to $B_R$, gives
\begin{align*}
 R^\alpha[\omega_{\mathcal B_p}]_{0,\alpha;B_R}
 &\leq 2^{1-\alpha} R\sup_{B_R}|D\omega_{\mathcal B_p}|
 \leq2^{\frac{1}{2}-\alpha}\,\widetilde K R
 \left(1-\frac{\widetilde K R}{2\sqrt2}\right)^{-3}
 =S.
 \qedhere
\end{align*}
\end{proof}

Recall the bounds $\tilde c_1,\tilde c_2$ for the transitions of parent charts from \cref{equation:parent-transition-second-derivative-bound,equation:parent-transition-third-derivative-bound}.

\begin{lemma}
\label[lemma]{lemma:ball-transition-bounds}
\begin{enumerate}[label= (\roman*), nosep, wide]
    \item Let $p,q\in M$ with $\mathcal{B}_p \cap \mathcal{B}_q \neq \emptyset$. Write $\phi_{pq}=\phi_p \circ \phi^{-1}_q$ for the transition map.
    We have
\begin{align*}
    |D\phi_{pq}|_{\mathrm{op}} \leq c_0,
    \quad
 R|D^2\phi_{pq}| \leq c_1,
\end{align*}
where
\begin{align*}
 c_0
 =K_0,
 \quad
 c_1
 =R K_0^{3/2}\tilde c_1.
\end{align*}
\item Let $x,y\in B_{2R}$ with $|x-y| < \frac{1}{c_0} 2R$. Write $p=\phi_q^{-1}((x+y)/2)$ and $\gamma(t) = \phi_q^{-1}((1-t)x+ty)$ for $0\leq t\leq 1$. Then the image of $\gamma$ is contained in $\mathcal{B}_p$ and the image of $\phi_p\circ\gamma$ is contained in $B_R$ with
\begin{align*}
    |\phi_{pq}(x)-\phi_{pq}(y)| \leq c_0 |x-y|.
\end{align*}
Furthermore we have
\begin{gather*}
    |D\phi_{pq}(x)-D\phi_{pq}(y)|_{\operatorname{op}} \leq \frac{1}{R} c_1 |x-y|,
    \\
    R^{1+\alpha} |D^2\phi_{pq}(x)-D^2\phi_{pq}(y)| \leq c_{1,\alpha} |x-y|^{\alpha},
\end{gather*}
where
\begin{align*}
    c_{1,\alpha} = \left(\frac{2}{c_0}\right)^{1-\alpha} R^2 K_0^2 \tilde c_2.
\end{align*}
\end{enumerate}
\end{lemma}

\begin{proof}
Proof of (i):
Write $g_{\mathcal B_p}$ for the matrix of the metric tensor with respect to the coordinate vector fields of $\phi_p$.
Let $y\in\phi_q(\mathcal B_p\cap\mathcal B_q)$, and set $\xi:=\phi_q^{-1}(y)$ and $T:=D\phi_{pq}(y)$.
Lemma \ref{lemma:eigenvalues-metric} implies
\begin{align}
 g_{\mathcal B_q}(y)
 &\leq K_0 I,
 \qquad
 g_{\mathcal B_p}\bigl(\phi_p(\xi)\bigr)
 \geq K_0^{-1}I,
 \qquad
 g_{\mathcal B_p}\bigl(\phi_p(\xi)\bigr)^{-1}
 \leq K_0 I,
 \label{equation:one-sided-metric-bounds}
\end{align}
where the third inequality follows from the second by taking inverses.
Hence
\begin{align*}
 K_0^{-1}T^TT
 \leq T^Tg_{\mathcal B_p}\bigl(\phi_p(\xi)\bigr)T
 =g_{\mathcal B_q}(y)
 \leq K_0 I,
\end{align*}
where in the first step we used \eqref{equation:one-sided-metric-bounds} after multiplication on the left by $T^T$ and on the right by $T$, in the second we applied the coordinate change to the metric, and in the third we used \eqref{equation:one-sided-metric-bounds} again.
Taking the largest eigenvalue gives $|T|_{\mathrm{op}}^2\leq K_0^2$, which proves the bound for $c_0$.

For $c_1$, write $a=a(p)$ and $b=a(q)$.
By the chain rule,
\begin{align}
 D^2\phi_{pq}(y)
 =Q_{a,p}D^2F_{ab}\bigl(\psi_b(\xi)\bigr)
 \left[Q_{b,q}^{-1}\,\cdot,Q_{b,q}^{-1}\,\cdot\right].
 \label{equation:ball-transition-second-jet}
\end{align}
We have
\begin{align}
 g_{\mathcal B_p}\bigl(\phi_p(\xi)\bigr)
 &=Q_{a,p}^{-T}G_a\bigl(\psi_a(\xi)\bigr)Q_{a,p}^{-1}
 =Q_{a,p}^{-T}Q_{a,\xi}^TQ_{a,\xi}Q_{a,p}^{-1}
 =(Q_{a,\xi}Q_{a,p}^{-1})^T(Q_{a,\xi}Q_{a,p}^{-1}).
 \label{equation:ball-metric-normaliser-identity}
\end{align}
Here the first equality is the coordinate-change formula for the metric, and the second uses $Q_{a,\xi}=G_a(\psi_a(\xi))^{1/2}$.
The analogous identity holds for the $\phi_q$-coordinates.
\Cref{equation:one-sided-metric-bounds,equation:ball-metric-normaliser-identity} therefore give\begin{align*}
 |Q_{a,p}Q_{a,\xi}^{-1}|_{\mathrm{op}}
 &\leq K_0^{1/2},
 &
 |Q_{b,\xi}Q_{b,q}^{-1}|_{\mathrm{op}}
 &\leq K_0^{1/2}.
\end{align*}
Combining \eqref{equation:ball-transition-second-jet} with the normalised parent bound \eqref{equation:parent-transition-second-derivative-bound}, we obtain
\begin{align*}
 |D^2\phi_{pq}(y)|
 &\leq
 |Q_{a,p}Q_{a,\xi}^{-1}|_{\mathrm{op}}
 \tilde c_1
 |Q_{b,\xi}Q_{b,q}^{-1}|_{\mathrm{op}}^2
 \leq K_0^{3/2}\tilde c_1.
\end{align*}
Multiplying by $R$ proves the second bound.

Proof of (ii): Suppose the curve $\tilde{\gamma}(t) = \phi_{pq}((x+y)/2+t(x-y))$ first exits $B_R$ at a time $\tau\in (0,1/2]$. We have
\begin{align*}
    |\tilde{\gamma}(\tau)|
    =
    |\tilde{\gamma}(\tau)-\tilde{\gamma}(0)|
    \leq
    \int_0^{\tau} |D\phi_{pq}((x+y)/2+t(x-y))(x-y)| dt
    \leq c_0 \tau |x-y| < \tau 2R \leq R,
\end{align*}
where we have used the assumption $|x-y| < \frac{1}{c_0}2R$. This is a contradiction. Similarly, we argue for the segment emanating from $p$ to $\phi_q^{-1}(y)$. This proves  $(\phi_p\circ\gamma)([0,1]) \subset B_R$.

If $|v|=1$ we have
    \begin{align*}
        |(D\phi_{pq}(x)-D\phi_{pq}(y))v| \leq \int_0^1 |D^2\phi_{pq}(y+t(x-y))[x-y,v]| dt \leq \frac{1}{R} c_1 |x-y|,
    \end{align*}
    and thus $|D\phi_{pq}(x)-D\phi_{pq}(y)|_{\operatorname{op}} \leq \frac{1}{R} c_1 |x-y|$.

For $\xi=\gamma(t)$, suppressing the coordinate arguments of the metric matrices, we have
\begin{align}
\begin{split}
 |Q_{a,p}Q_{a,\xi}^{-1}|_{\mathrm{op}}^2
 &=|g_{\mathcal B_p}^{-1}|_{\mathrm{op}}
 \leq K_0,\\
 |Q_{b,\xi}Q_{b,q}^{-1}|_{\mathrm{op}}^2
 &=|g_{\mathcal B_q}|_{\mathrm{op}}
 \leq K_0,
\end{split}
 \label{equation:normalising-matrix-bounds}
\end{align}
where the equalities follow from \eqref{equation:ball-metric-normaliser-identity} and the inequalities follow from Lemma \ref{lemma:eigenvalues-metric}.
Taking square roots we can estimate:
\begin{align*}
 |D^3\phi_{pq}(\phi_q(\xi))|
 &\leq
 |Q_{a,p}Q_{a,\xi}^{-1}|_{\mathrm{op}}
 |Q_{b,\xi}Q_{b,q}^{-1}|_{\mathrm{op}}^3
 \tilde c_2
 \leq
 K_0^{2}\tilde c_2.
\end{align*}
Here we differentiated \eqref{equation:ball-transition-second-jet} once more and used the bound \eqref{equation:parent-transition-third-derivative-bound} in the first inequality, and \eqref{equation:normalising-matrix-bounds} in the second.
Integrating along the segment gives
\begin{align*}
 |D^2\phi_{pq}(x)-D^2\phi_{pq}(y)|
 \leq
 K_0^{2}
 \tilde c_2|x-y|.
\end{align*}
Dividing by $|x-y|^\alpha$, using $|x-y|<\frac{1}{c_0}2R$, and multiplying by $R^{1+\alpha}$ gives
\begin{equation*}
 R^{1+\alpha}
 \frac{|D^2\phi_{pq}(x)-D^2\phi_{pq}(y)|}{|x-y|^\alpha}
 \leq
 \left(\frac{2}{c_0} \right)^{1-\alpha}R^2
K_0^{2}\tilde c_2
 =c_{1,\alpha}.\qedhere
\end{equation*}
\end{proof}

\subsection{Function spaces and residual bound}
\label{section:functionspaces}

We will now define the function spaces $\mathscr{X}$ and $\mathscr{Y}$ to which we will apply the contraction mapping theorem \ref{lemma:inverse-function-thm}.

We first define global Hölder spaces.
With the local Hölder norms as defined in \cref{section:functionspaces}, set
\begin{align*}
    \|{u}'_{C^{2,\alpha}(M)} = \sup_{p\in M} \|{u \circ \phi_p^{-1}}'_{C^{2,\alpha}(B_{R})},
    \quad \|{u}''_{C^{2,\alpha}(M)} = \sup_{p\in M} \|{u \circ \phi_p^{-1}}''_{C^{2,\alpha}(B_{2R})},
    \\
    \|{u}'_{C^{0,\alpha}(M)} = \sup_{p\in M} \|{u\circ \phi_p^{-1}}'_{C^{0,\alpha}(B_R)}, \quad  \|{u}''_{C^{0,\alpha}(M)} = \sup_{p\in M} \|{u\circ \phi_p^{-1}}''_{C^{0,\alpha}(B_{2R})}.
\end{align*}
The Sobolev norm on $L^2_2(M)$ is defined by
\[
\|{u}_{L^2_2(M)}
=
\|{u}_{L^2(M)}+\|{\nabla u}_{L^2(M)}+\|{\nabla^2 u}_{L^2(M)}.
\]
Let $\mathscr{X}$ be the space of $G$-invariant functions on $M$ which are in the H\"older space $C^{2,\alpha}$ and are $L^2$-orthogonal to the constant functions, or equivalently have mean zero:
\begin{align*}
    \mathscr{X} = \{ u\in C^{2,\alpha}(M): \int_M u \mu_{\omega} = 0, u\,\text{is}\, G-\text{invariant} \}.
\end{align*}
We will equip $\mathscr{X}$ with the norm
\begin{align*}
    \|{u}_{\mathscr{X}} = 
    R^{-2}  \|{u}'_{C^{2,\alpha}(M)}
    +
    a R^{-2-\frac{n}{q}} \|{u}_{L^2_2(M)}
\end{align*}
for $R\leq R_0$, a constant $a>0$ and $1\leq q < \infty$ to be determined later.
Define $\mathscr{Y}$ as the space of $G$-invariant $C^{0,\alpha}$-functions on $M$ with mean zero:
\begin{align*}
    \mathscr{Y} = \{ u\in C^{0,\alpha}(M): \int_M u \mu_{\omega} = 0, u\,\text{is}\, G-\text{invariant} \}.
\end{align*}
We will equip $\mathscr{Y}$ with the norm
\begin{align*}
    \|{u}_{\mathscr{Y}} = 
   \|{u}'_{C^{0,\alpha}(M)}
    +
    a R^{-2-\frac{n}{q}} \|{u}_{L^2(M)}.
\end{align*}

Using rigorous computing, in \cref{proposition:certified-residual} we bounded the Riemannian gradient of the residual.
In the following proposition we translate this into a bound for the Hölder seminorm with respect to the small ball charts $\{\mathcal{B}_p\}_{p \in M}$:

\begin{proposition}
\label{prop-residual-bound}
The residual $f$ from \cref{eqn:residual-def} satisfies $f \in \mathscr{Y}$ and 
\begin{align*}
    \sup_{p\in M} R^{\alpha} [f\circ \phi_p^{-1}]_{0,\alpha;B_R} \leq 2^{1-\alpha} R K_0^{\frac{1}{2}} \sup_{M} |\nabla f|_g.
\end{align*}
    \end{proposition}
\begin{proof}
The function $f$ is smooth and $G$-invariant.
The normalisation \cref{eqn:volume-normalisation} gives $\int_M f\,\mu_\omega=\int_M\mu_\omega-\int_M\hat\mu=0$, hence $f \in \mathscr{Y}$.

For $p\in M$ and $x,y\in B_R$ there exists $\xi$ on the line segment between $x$ and $y$ such that by Lemma \ref{lemma:eigenvalues-metric}
\begin{align*}
    |(f\circ \phi_p^{-1})(x)-(f\circ \phi_p^{-1})(y)| \leq |D(f\circ \phi_p^{-1})(\xi)|_{g_{\mathcal{B}_p}} |x-y|_{g_{\mathcal{B}_p}}
    \leq |(\nabla^g f)(\phi_p^{-1}(\xi))|_g K_0^{\frac{1}{2}} |x-y|.
\end{align*}
The statement follows.
\end{proof}

\subsection{Estimate of the inverse of the Laplacian}
\label{subsection:estimate-inverse-laplacian}

In this section we will prove an explicit bound on the inverse of the Laplace operator $\Delta_g: \mathscr{X}\rightarrow \mathscr{Y}$. Because the linearisation of our perturbation problem is a multiple of the Laplace operator, we then get an explicit bound for $C$ in our fixed point method from Lemma \ref{lemma:inverse-function-thm}. The main ingredient is a Schauder estimate for $\Delta_g$. Because the Schauder estimate is an interior estimate, we need to compare the Hölder norms for functions on $M$ obtained from ball charts of radius $2R$ with the Hölder norms from ball charts with radius $R$. For this we need the control on the coordinate transitions of our ball charts from Lemma \ref{lemma:ball-transition-bounds}.

\begin{lemma}
\label{lemma:doubling-holo-atlas}
\begin{enumerate}[label= (\roman*), nosep, wide]
\item
     We have 
        \begin{align*}
          \|{u }''_{C^{2,\alpha}(M)}
         \leq \Cl{doubling-2-alpha}
        \|{u }'_{C^{2,\alpha}(M)},
        \end{align*}
        where 
        \begin{align*}
        \Cr{doubling-2-alpha}
        =
      1+c_0+\max\{c_0^2,c_1\} +\max\left\{  2^{1-\alpha}c_0^{2+\alpha},  3\cdot2^{1-\alpha}c_0^\alpha c_1,  c_{1,\alpha} \right\}.
        \end{align*}
        \item
    We have 
    \begin{align*}
          \|{u }''_{C^{0,\alpha}(M)}
         \leq \Cl{doubling-0-alpha}
         \|{u }'_{C^{0,\alpha}(M)},
        \end{align*}
        where 
        \begin{align*}
        \Cr{doubling-0-alpha}
        =
        1+2^{1-\alpha} c_0^{\alpha}.
        \end{align*}
    \end{enumerate}
    \end{lemma}

\begin{proof}
    Proof of (i): 
    \\
    \textbf{Estimate of first derivative:}
    First suppose that $\phi_p^{-1}(x) \in \phi_p^{-1}(B_{2R}) \cap \phi_q^{-1}(B_R)$.

    We have 
\begin{align*}
   (\nabla(u\circ \phi_p^{-1})(x))^T = (\nabla(u\circ \phi_q^{-1})(\phi_{qp}(x)))^T D\phi_{qp}(x).
\end{align*}
This gives
\begin{align*}
    R\left|\nabla (u \circ \phi_p^{-1})\right|(x)
    \leq  c_0 R [u\circ \phi_q^{-1}]_{1,0;B_R}
     .
\end{align*}
    
    \noindent    
\textbf{Estimate of second derivative:}
We have
    \begin{align}
D^2(u \circ \phi_p^{-1})(x)[\, \cdot\, , \, \cdot\,] =& D^2(u\circ \phi_q^{-1})(\phi_{qp}(x))[D\phi_{qp}(x)\, \cdot\,, D\phi_{qp}(x)\, \cdot\, ] \nonumber \\
&+ D(u\circ \phi_q^{-1})(\phi_{qp}(x))\circ D^2\phi_{qp}(x)[\, \cdot\, , \, \cdot\,].
\label{two-der}
\end{align}
Now for $n\times n$- matrices $A$ and $M$ we have 
\begin{align}
|A^T M A| \leq |A|_{\operatorname{op}}^2 |M|. 
\label{matrix-inequality}
\end{align}
Indeed, by the singular value decomposition of $A$ we write $A=U \Sigma V^T$ with orthogonal matrices $U$ and $V$. Using $|A|_{\operatorname{op}} = |\Sigma|_{\operatorname{op}}$ and $|VBV^T| = |B|$ for an orthogonal matrix $V$, we have
\begin{align*}
    |A^TMA| = |V\Sigma U^T M U \Sigma V^T|= |\Sigma U^T M U \Sigma| \leq |\Sigma|_{\operatorname{op}}^2 |U^T M U| = |A|_{\operatorname{op}}^2 |M|.
\end{align*}
Applying \eqref{matrix-inequality} to \eqref{two-der} we get
\begin{align*}
R^2\left|D^2(u\circ \phi_p^{-1})\right|(x)
\leq 
 c_0^2 R^2 [u\circ \phi_q^{-1}]_{2,0; B_R} +  c_1 R[u\circ \phi_q^{-1}]_{1,0;B_R}
 .
\end{align*}

   \noindent   
    \textbf{Hölder bound of second derivative:}
    For $p\in M$ let $\hat{x},\hat{y}\in \mathcal{B}_p$
    and write $x=\phi_p(\hat{x}), y = \phi_p(\hat{y}) \in B_{2R}$. We will distinguish two cases.
    \\
    Case 1: $|x-y| < \frac{1}{c_0}2R$. 
    Set $q = \phi_p^{-1}((x+y)/2)$ and write $x' = \phi_{q}(\hat{x})$ and $y' = \phi_{q}(\hat{y})$. Then by Lemma \ref{lemma:ball-transition-bounds}(ii)
    the points $x',y'$ are in $B_R$ and we have
    \begin{gather*}
         |x'-y'| \leq c_0 |x-y|,
         \\
         |D\phi_{qp}(x)-D\phi_{qp}(y)|_{\operatorname{op}} \leq \frac{1}{R} c_1 |x-y| \leq \frac{c_1}{R} \left(\frac{2R}{c_0} \right)^{1-\alpha} |x-y|^{\alpha} = 2^{1-\alpha} c_0^{\alpha-1} c_1R^{-\alpha} |x-y|^{\alpha} ,
         \\
         R^{1+\alpha} |D^2\phi_{qp}(x)-D^2\phi_{qp}(y)| \leq c_{1,\alpha} |x-y|^{\alpha}.
    \end{gather*}
    From \eqref{two-der} we get
    \begin{align*}
        D^2(u\circ \phi_p^{-1})(x) - D^2(u\circ \phi_p^{-1})(y)
        =&
        \{ D^2(u\circ \phi_q^{-1})(x') - D^2 (u\circ \phi_q^{-1})(y')\}[D\phi_{qp}(x)\, \cdot\, , D\phi_{qp}(x)\, \cdot\,]
        \\
        &+ D^2(u\circ \phi_q^{-1})(y') [\{D\phi_{qp}(x)-D\phi_{qp}(y)\}\, \cdot\, , D\phi_{qp}(x) \, \cdot\,]
        \\
        &+ D^2(u\circ \phi_q^{-1})(y') [D\phi_{qp}(y) \, \cdot\,, \{D\phi_{qp}(x)-D\phi_{qp}(y)\}\, \cdot\,]
        \\
        &+ \{ D(u\circ \phi_q^{-1})(x')-D(u\circ \phi_q^{-1})(y')\}\circ D^2\phi_{qp}(x) [\, \cdot\, , \, \cdot\,]
        \\
        &+ D(u\circ \phi_q^{-1})(y') \circ \{D^2\phi_{qp}(x) - D^2\phi_{qp}(y) \} [\, \cdot\, , \, \cdot\,].
    \end{align*}
    Thus we get
    \begin{align*}
        R^{2+\alpha} \frac{|D^2(u\circ \phi_p^{-1})(x) - D^2(u\circ \phi_p^{-1})(y)|}{|x-y|^{\alpha}}
        \leq&
        c_0^{2+\alpha} R^{2+\alpha} [u\circ \phi_q^{-1}]_{2,\alpha;B_R}
        +
       2 \cdot2^{1-\alpha} c_0^{\alpha} c_1 R^2 [u\circ \phi_q^{-1}]_{2,0;B_R}
        \\
        &+
        c_0^{\alpha} c_1 R^{1+\alpha} [u\circ \phi_q^{-1}]_{1,\alpha;B_R}
        \\
        &+  c_{1,\alpha} R [u\circ \phi_q^{-1}]_{1,0;B_R}
        \\
        \leq &
        c_0^{2+\alpha} R^{2+\alpha} [u\circ \phi_q^{-1}]_{2,\alpha;B_R}
        +
        3 \cdot2^{1-\alpha} c_0^{\alpha} c_1 R^2 [u\circ \phi_q^{-1}]_{2,0;B_R}
        \\
        &+ c_{1,\alpha} R [u\circ \phi_q^{-1}]_{1,0;B_R}
       .
        \end{align*}
    
Case 2: $|x-y| \geq \frac{1}{c_0}2R$. 
Let $q_1 = \phi_p^{-1}(x)$ and $q_2=\phi_p^{-1}(y)$.
Then with the $C^2$-estimate above we get
\begin{align*}
     R^{2+\alpha} \frac{|D^2(u\circ \phi_p^{-1})(x) - D^2(u\circ \phi_p^{-1})(y)|}{|x-y|^{\alpha}}
    \leq 
     \left(\frac{c_0}{2}\right)^{\alpha}  \left\{ R^2\left|D^2(u\circ \phi_p^{-1})\right|(x) + R^2\left|D^2(u\circ \phi_p^{-1})\right|(y) \right\}
    \\
    \leq 
     2^{1-\alpha} \sup_{q\in M}\{ c_0^{2+\alpha} R^2 [u\circ \phi_q^{-1}]_{2,0;B_R}
    + c_0^{\alpha}  c_1 R [u\circ \phi_q^{-1}]_{1,0;B_R}\}.
\end{align*}
Combining the two cases gives
\begin{align*}
 R^{2+\alpha}[u\circ\phi_p^{-1}]_{2,\alpha;B_{2R}}
 \leq
 \max\left\{
  2^{1-\alpha}c_0^{2+\alpha},
  3\cdot2^{1-\alpha}c_0^\alpha c_1,
  c_{1,\alpha}
 \right\}
 \|{u}'_{C^{2,\alpha}(M)}.
\end{align*}
The remaining three terms in
$\|{u\circ\phi_p^{-1}}''_{C^{2,\alpha}(B_{2R})}$
are bounded by $1$, $c_0$ and $\max\{c_0^2,c_1\}$
times $\|{u}'_{C^{2,\alpha}(M)}$, respectively.
Adding these bounds and taking the supremum over $p\in M$
proves (i).
\\
Proof of (ii):  For $p\in M$ let $\hat{x},\hat{y}\in \mathcal{B}_p$
    and write $x=\phi_p(\hat{x}), y = \phi_p(\hat{y}) \in B_{2R}$.
Like in the Hölder bound of (i) we distinguish two cases.
\\
Case 1: $|x-y| < \frac{1}{c_0}2R$. 
Analogous to the proof of (i) set $q = \phi_p^{-1}((x+y)/2)$ and write $x' = \phi_{q}(\hat{x})$ and $\phi_{q}(\hat{y})$. By Lemma \ref{lemma:ball-transition-bounds}(ii)
    the points $x',y'$ are in $B_R$ and we have
   $|x'-y'| \leq c_0 |x-y|$.
Then 
\begin{align*}
    \frac{|u(\hat{x})-u(\hat{y})|}{|x-y|^{\alpha}} \leq c_0^{\alpha}  \frac{|u(\hat{x})-u(\hat{y})|}{|x'-y'|^{\alpha}}.
\end{align*}
Case 2: $|x-y| \geq \frac{1}{c_0}2R$. 
Then
\begin{align*}
    R^{\alpha}  \frac{|u(\hat{x})-u(\hat{y})|}{|x-y|^{\alpha}} \leq R^{\alpha} \left(\frac{c_0}{2R} \right)^{\alpha} 2 \|{u}_{C^0(M)}
    \leq c_0^{\alpha} 2^{1-\alpha} \|{u}_{C^0(M)}. 
\end{align*}
Thus
\begin{align*}
    &[u\circ \phi_p^{-1}]_{0,0;B_{2R}} + R^{\alpha} [u\circ \phi_p^{-1}]_{0,\alpha;B_{2R}}
    \\
    \leq&\, \sup_{q\in M} [u\circ \phi_q^{-1}]_{0,0;B_{R}} + \max\{ \sup_{q\in M} c_0^{\alpha} R^{\alpha} [u\circ \phi_q^{-1}]_{0,\alpha;B_R}, \sup_{q\in M} c_0^{\alpha}2^{1-\alpha}  [u\circ \phi_q^{-1}]_{0,0;B_{R}} \}   
    \\
    \leq&\, (1+2^{1-\alpha} c_0^{\alpha} ) \|{u}'_{C^{0,\alpha}(M)}.
\end{align*}
\end{proof}

\begin{lemma}
\label{lemma:Laplace-difference}
   We have
    \begin{align*}
        R^2 \|{(\Delta_g u)\circ \phi_p^{-1} - \Delta_0(u\circ \phi_p^{-1})  }''_{C^{0,\alpha}(B_{2R})} \leq
        (\varepsilon_0 +  \varepsilon_{\alpha})
         \|{u\circ \phi_p^{-1}}''_{C^{2,\alpha}(B_{2R})},
    \end{align*}
    where
    $\Delta_0$ denotes the Laplace operator on flat $\mathbb{R}^n$.
\end{lemma}
\begin{proof}
    Because the coordinate map $\phi_p$ is holomorphic, the coordinate functions are harmonic with respect to $g$ and thus $(\Delta_g u)\circ \phi_p^{-1} - \Delta_0(u\circ \phi_p^{-1})  = (g^{ij}_{\mathcal{B}_p}-\delta^{ij}) \partial_{ij}(u\circ \phi_p^{-1})$. The statement then follows from the inequality $|\operatorname{Tr}(AB)| \leq |A| \cdot |B|$ for $A,B\in \mathbb{R}^{n\times n}$ and the bounds \eqref{equation:metric-coordinate-bounds}.
\end{proof}

\begin{proposition}
\label{proposition-global-Schauder}
    If we have
    \begin{align}
          \varepsilon_0 +  \varepsilon_{\alpha} \leq \frac{1}{4 \Cr{loc-Schauder-1} \Cr{doubling-2-alpha}},
 \quad
 \varepsilon_0 \leq \frac{n}{2 \Cr{loc-Schauder-2}},
         \label{smallness-condition}
    \end{align}
    then
    \begin{align*}
  R^{-2} \|{u}'_{C^{2,\alpha}(M)}
         \leq 
         \Cl{Lq-Schauder-1}
         \|{\Delta_g u}'_{C^{0,\alpha}(M)}
         + \Cl{Lq-Schauder-2} R^{-2-\frac{n}{q}} \|{u}_{L^q(M)}
         .    
\end{align*}
    where
    \begin{align*}
        \Cr{Lq-Schauder-1} = 2(\Cr{loc-Schauder-1} \Cr{doubling-0-alpha} + \tfrac{1}{2n} \Cr{loc-Schauder-2}),
        \quad
        \Cr{Lq-Schauder-2} = 2 \Cr{loc-Schauder-2} n^{\frac{1}{q}} \omega_{n-1}^{-\frac{1}{q}} K_0^{\frac{n}{2q}}
        .
    \end{align*}
    Here $\omega_{n-1} = \vol_{n-1}(S^{n-1})$.
    \end{proposition}

\begin{remark}
    As $R\to0$, $\Cr{doubling-2-alpha}\to3+2^{1-\alpha}$ and $\varepsilon_0+\varepsilon_\alpha\to0$.
    Hence \eqref{smallness-condition} holds for sufficiently small $R$.
\end{remark}
    
\begin{proof}[Proof of Proposition \ref{proposition-global-Schauder}]
    Applying the local Schauder estimate from Proposition \ref{prop:local-Schauder} on a coordinate ball, Lemma \ref{lemma:Laplace-difference} and Lemma \ref{lemma:doubling-holo-atlas} gives
    \begin{gather*}
        R^{-2} \|{u\circ \phi_p^{-1}}'_{C^{2,\alpha}(B_R)} 
        \leq 
        \Cr{loc-Schauder-1} \|{\Delta(u\circ \phi_p^{-1})}''_{C^{0,\alpha}(B_{2R})} + \Cr{loc-Schauder-2} R^{-2} \|{u\circ \phi_p^{-1}}_{C^0(B_{2R})}
        \\
        \leq
        \Cr{loc-Schauder-1} \|{(\Delta_g u)\circ \phi_p^{-1}}''_{C^{0,\alpha}(B_{2R})}
        +
        \Cr{loc-Schauder-2} R^{-2} \|{u}_{C^0(M)}
        +
        \Cr{loc-Schauder-1} (\varepsilon_0 +  \varepsilon_{\alpha})  R^{-2} \|{u\circ \phi_p^{-1}}''_{C^{2,\alpha}(B_{2R})}.
        \\
        \leq
       \Cr{loc-Schauder-1} \Cr{doubling-0-alpha} \|{\Delta_g u}'_{C^{0,\alpha}(M)}
        +
        \Cr{loc-Schauder-2} R^{-2} \|{u}_{C^0(M)}
        +
        \Cr{loc-Schauder-1} (\varepsilon_0 +  \varepsilon_{\alpha})\Cr{doubling-2-alpha}  R^{-2} \|{u}'_{C^{2,\alpha}(M)}.
    \end{gather*}
    Taking the supremum over all coordinate balls gives
    \begin{align*}
        R^{-2} \|{u}'_{C^{2,\alpha}(M)} \leq  \Cr{loc-Schauder-1} \Cr{doubling-0-alpha} \|{\Delta_g u}'_{C^{0,\alpha}(M)}
        +
        \Cr{loc-Schauder-2} R^{-2} \|{u}_{C^0(M)}
        +
        \Cr{loc-Schauder-1} (\varepsilon_0 +  \varepsilon_{\alpha})\Cr{doubling-2-alpha} R^{-2} \|{u}'_{C^{2,\alpha}(M)}.
    \end{align*}
    Applying Lemma \ref{lemma-meanvalue} at a coordinate ball centred at $p$ shows
    \begin{align*}
        |u(p)| &\leq  \left( \fint_{B_R} |u\circ \phi_p^{-1}|^q dx \right)^{\frac{1}{q}} + \frac{R^2}{2n} \|{\Delta (u\circ \phi_p^{-1})}_{C^0(B_R)}
        \\
        &\leq n^{\frac{1}{q}}
        \omega_{n-1}^{-\frac{1}{q}} K_0^{\frac{n}{2q}}R^{-\frac{n}{q}} \|{u}_{L^q(M)} + \frac{1}{2n} R^2 \|{\Delta_g u}_{C^0(M)} + \frac{1}{2n} \varepsilon_0  \|{u}'_{C^{2,\alpha}(M)}.
    \end{align*}
    Here we used $\vol(B_R) = R^n \frac{1}{n}\omega_{n-1}$.
     Thus we get
     \begin{align*}
         R^{-2} \|{u}'_{C^{2,\alpha}(M)} \leq&\,  (\Cr{loc-Schauder-1} \Cr{doubling-0-alpha} + \tfrac{1}{2n} \Cr{loc-Schauder-2}) \|{\Delta_g u}'_{C^{0,\alpha}(M)}
        +
        \Cr{loc-Schauder-2} n^{\frac{1}{q}} \omega_{n-1}^{-\frac{1}{q}} K_0^{\frac{n}{2q}}R^{-2-\frac{n}{q}} \|{u}_{L^q(M)}
        \\
        &+
        \Cr{loc-Schauder-1} (\varepsilon_0 +  \varepsilon_{\alpha}) \Cr{doubling-2-alpha} R^{-2} \|{u}'_{C^{2,\alpha}(M)}
        +
        \Cr{loc-Schauder-2} \frac{\varepsilon_0}{2n} R^{-2} \|{u}'_{C^{2,\alpha}(M)}
        .
     \end{align*}

    The statement follows.
\end{proof}

\begin{proposition}
\label{proposition:Laplace-inverse}
    Assume \eqref{radius-cond-1}, \eqref{radius-cond-2} and \eqref{smallness-condition}, and let
    $\frac{1}{2}-\frac{2}{n}<\frac{1}{q}<\frac{1}{2}-\frac{1}{n}$.
    Suppose that $\kappa>0$, $\diam(M,g)<D$, $0<\mu\leq\lambda_{1,G}$ for $-\Delta_g$, $a\geq\Cr{L22-Sobolev-update}$, and
    \begin{align}
        \operatorname{Ric}_g\geq-(n-1)\kappa g.
    \end{align}
     Then we have
     \begin{align*}
    \|{u}_{\mathscr{X}} \leq C_{\Delta^{-1}} \|{\Delta_g u}_{\mathscr{Y}},
\end{align*}
where
\begin{align*}
    C_{\Delta^{-1}} = \max\{ \Cr{Lq-Schauder-1},  (\Cr{Lq-Schauder-2}  + 1)  \Cr{L2-ell-estimate} (1+\mu^{-1}) \}.
\end{align*}
\end{proposition}
\begin{proof}
By \cref{proposition-global-Schauder} and \cref{cor:L22-Sobolev-embedding-update} we have
\begin{align*}
  R^{-2} \|{u}'_{C^{2,\alpha}(M)}
         \leq 
         \Cr{Lq-Schauder-1}
         \|{\Delta_g u}'_{C^{0,\alpha}(M)}
         + \Cr{Lq-Schauder-2} \Cr{L22-Sobolev-update} R^{-2-\frac{n}{q}} \|{u}_{L^2_2(M)}
         .    
\end{align*}
By \cref{prop:L2-ell-estimate} we have
    \begin{align*}
         \|{u}_{L^2_2(M)} \leq  \Cr{L2-ell-estimate} ( \|{\Delta_g u}_{L^2(M)} + \|{u}_{L^2(M)})
        \leq  \Cr{L2-ell-estimate} (1+\mu^{-1})  \|{\Delta_g u}_{L^2(M)}.
    \end{align*}   
Combining these estimates gives 
\begin{align*}
    \|{u}_{\mathscr{X}} \leq \Cr{Lq-Schauder-1}
         \|{\Delta_g u}'_{C^{0,\alpha}(M)} + (\Cr{Lq-Schauder-2} \Cr{L22-Sobolev-update} + a) \Cr{L2-ell-estimate} (1+\mu^{-1}) R^{-2-\frac{n}{q}}\|{\Delta_g u}_{L^2(M)}.
\end{align*}
Since $a\geq\Cr{L22-Sobolev-update}$, we get the statement.
\end{proof}

\subsection{Estimate of the non-linear term}
\label{subsection:estimate-non-linear}

In this section we will derive the bound on the nonlinear term
\begin{align}
    Q(u) 
    =
    \sum_{k=2}^m \binom{m}{k} \frac{(i\partial\bar{\partial}u)^k\wedge \omega^{m-k}}{\omega^m},
    \end{align}
which is required to apply Lemma \ref{lemma:inverse-function-thm}.
Throughout we take $m=2$, i.e. real dimension $n=4$.

Set $B = B_R$. For a differential form $\rho$  
on $B$
we define
\begin{align*}
    [\rho]_{0;B} = \sup_{x\in B} |\rho(x)|,
    \quad
    [\rho]_{0,\alpha;B} = \sup_{x\neq y\in B} \frac{|\rho(x)-\rho(y)|}{|x-y|^{\alpha}}.
\end{align*}
Here $|\rho(x)|$ denotes the norm as a differential form determined by the flat metric. We will also need an adaptation by using the Riemannian metric $g$ in the coordinate system defined by $\phi_p$:
\begin{align*}
    [\rho]_{0;B,g} = \sup_{x\in B} |\rho(x)|_{g_{\mathcal{B}_p}(x)},
    \quad
    [\rho]_{0,\alpha;B,g} = \sup_{x\neq y\in B} \frac{|\rho(x)-\rho(y)|_{g_{\mathcal{B}_p}(x)}}{|x-y|^{\alpha}}.
\end{align*}
We emphasize that in this definition the difference $\rho(x)-\rho(y)$ is still defined by the parallel transport of the flat metric.

In dimension $4$ for $\sigma,\tau\in \Lambda^2(\mathbb{R}^4)^*$ we have
\begin{align*}
    |\sigma\wedge\tau| \leq |\sigma| \cdot |\tau|,
\end{align*}
which gives the product estimates
\begin{align*}
[\sigma\wedge\tau]_{0;B}
&\leq [\sigma]_{0;B}[\tau]_{0;B},\\
[\sigma\wedge\tau]_{0,\alpha;B}
&\leq
[\sigma]_{0,\alpha;B}[\tau]_{0;B}
+[\sigma]_{0;B}[\tau]_{0,\alpha;B}.
\end{align*}
The $C^0$ estimate also holds with respect to $g$.
The H\"older estimate becomes
\begin{align}
\label{eqn:holder-wedge-estimate}
[\sigma\wedge\tau]_{0,\alpha;B,g}
&\leq
[\sigma]_{0,\alpha;B,g}[\tau]_{0;B,g}
+\left(\sup_{x,y\in B}|\sigma(y)|_{g_{\mathcal{B}_p}(x)}\right)
[\tau]_{0,\alpha;B,g}\\
&\leq
[\sigma]_{0,\alpha;B,g}[\tau]_{0;B,g}
+K_0[\sigma]_{0;B}[\tau]_{0,\alpha;B,g},
\end{align}
where we used \cref{lemma:eigenvalues-metric} in the last step.

In the following we will often need
\begin{align}
    |\omega|_g=\sqrt{2},
    \qquad
    |\omega^2|_g=2.
    \label{norm-omega}
\end{align}

\begin{lemma}
    Let $m=2$ and $\gamma,\eta\in\Omega^2(B)$. Then
    \begin{subequations}
        \begin{align}
            [\gamma^2]_{0,\alpha;B} &\leq 2 [\gamma]_{0,\alpha;B} [\gamma]_{0;B},
            \\
            [\gamma^2-\eta^2]_{0;B} &\leq [\gamma-\eta]_{0;B} ([\gamma]_{0;B} + [\eta]_{0;B}),
            \\
            R^{\alpha}[\gamma^2-\eta^2]_{0,\alpha;B} &\leq \|{\gamma-\eta}'_{0,\alpha;B} \|{\gamma+\eta}'_{0,\alpha;B}.
        \end{align}
    \end{subequations}
\end{lemma}
\begin{proof}

Apply the Euclidean product estimates to $\gamma\wedge\gamma$
and $(\gamma-\eta)\wedge(\gamma+\eta)=\gamma^2-\eta^2$.
\end{proof}

\begin{lemma}
\label{lemma:estimate-nonlinear-term}
   Write $\gamma = (\phi_p^{-1})^*(i\partial\bar{\partial}u)$ and $\eta = (\phi_p^{-1})^*(i\partial\bar{\partial}v)$. By $K_0$ denote  the bound on the eigenvalues of $g_{\mathcal{B}_p}$ from Lemma \ref{lemma:eigenvalues-metric}
   and by $S$  the Hölder bound on $\omega$ from Lemma 
   \ref{lemma:Holder-bound-omega}
  .
Let $m=2$. Then we have
\begin{align*}
    \|{(Q(u)-Q(v))\circ \phi_p^{-1}}'_{0,\alpha;B}
    \leq \frac{K_0^2}{2} (2+\sqrt{2}K_0^3 S) \|{\gamma-\eta}'_{0,\alpha;B} \|{\gamma+\eta}'_{0,\alpha;B}.
\end{align*}
\end{lemma}
\begin{proof}

We have
\[
(Q(u)-Q(v))\circ\phi_p^{-1}
=
\frac{\gamma^2-\eta^2}{\omega^2}.
\]
Since $|\omega^2|_g=2$, the Euclidean product estimate gives
\begin{align*}
\left[\frac{\gamma^2-\eta^2}{\omega^2} \right]_{0;B}    
\leq \frac{K_0^2}{2} [\gamma-\eta]_{0;B} [\gamma+\eta]_{0;B}.
\end{align*}
For scalar functions $f,h$ with $h$ nowhere zero and $x,y\in B$,
\[
\frac{f(x)}{h(x)}-\frac{f(y)}{h(y)}
=
\frac{f(x)-f(y)}{h(x)}
+\frac{f(y)}{h(y)}\frac{h(y)-h(x)}{h(x)}.
\]
Applying this quotient identity, using $|\omega^2|_g=2$
and Lemma \ref{lemma:eigenvalues-metric}, we obtain
\begin{align}
    \left[\frac{\gamma^2-\eta^2}{\omega^2}\right]_{0,\alpha;B}
    \leq
    \frac{K_0^2}{2}[\gamma^2-\eta^2]_{0,\alpha;B}
    +\frac{K_0^2}{4}[\omega^2]_{0,\alpha;B,g}
    [\gamma^2-\eta^2]_{0;B}.
    \label{nonlinear-Hoelder}
\end{align}
Since $[\omega]_{0;B}\leq\sqrt{2}K_0$, we have
\begin{align*}
R^\alpha[\omega^2]_{0,\alpha;B,g}
&\leq K_0^2R^\alpha[\omega^2]_{0,\alpha;B}
\leq 2K_0^2R^\alpha
[\omega]_{0,\alpha;B}[\omega]_{0;B}
\leq 2^{3/2}K_0^3S.
\end{align*}
Thus 
\begin{align*}
   \|{(Q(u)-Q(v))\circ\phi_p^{-1}}'_{0,\alpha;B}
    \leq \frac{K_0^2}{2} (2+\sqrt{2}K_0^3 S) \|{\gamma-\eta}'_{0,\alpha;B} \|{\gamma+\eta}'_{0,\alpha;B}.
\end{align*}
\end{proof}

\begin{lemma}
\label{lemma:complexHessian-realHessian}
    For $f\in C^{2,\alpha}(B)$ we have
    \begin{align}
        |i\partial\bar{\partial}f|_g \leq |\nabla^g d f|_g
        \label{equ:complexHessian-realHessian}
    \end{align}
    and
    \begin{align*}
        [i\partial\bar{\partial}f]_{0;B} \leq  [f]_{2,0;B},
        \quad
        [i\partial\bar{\partial}f]_{0,\alpha;B} \leq  [f]_{2,\alpha;B}.
    \end{align*}
\end{lemma}

\begin{proof}
    Let $z_1=x_1+iy_1$, $z_2=x_2+iy_2$ be holomorphic coordinates on $U\subset M$.
    We have
    \begin{align*}
        \frac{\partial^2f}{\partial z_k \partial\bar{z}_l} 
        &= 
        \frac{1}{4}(\partial_{x_k}-i\partial_{y_k})(\partial_{x_l}+i\partial_{y_l})f=
        \frac{1}{4} \{\partial_{x_k}\partial_{x_l}f+\partial_{y_k}\partial_{y_l}f+i(\partial_{x_k}\partial_{y_l}f-\partial_{y_k}\partial_{x_l}f)\},
        \\[10pt]
        dz_k\wedge d\bar{z}_l 
        &= 
        (dx_k+idy_k)\wedge(dx_l-i dy_l)
        =
        dx_k\wedge dx_l + dy_k\wedge dy_l + i (-dx_k\wedge dy_l+dy_k\wedge dx_l)
    \end{align*}
    and thus
    \begin{align*}
        i\partial\bar{\partial}f
        =&
        -\mathrm{Im}(\partial\bar{\partial}f)
        \\
        =&
        -\frac{1}{4}\sum_{k,l}\Bigl[        (\partial_{x_k}\partial_{x_l}f+\partial_{y_k}\partial_{y_l}f)
        (-dx_k\wedge dy_l+dy_k\wedge dx_l)
        +
        (\partial_{x_k}\partial_{y_l}f-\partial_{y_k}\partial_{x_l}f)
        (dx_k\wedge dx_l + dy_k\wedge dy_l)\Bigr]
        \\
        =&
        \frac{1}{2}
        \sum_{k,l}(\partial_{x_k}\partial_{x_l}f+\partial_{y_k}\partial_{y_l}f)
        dx_k\wedge dy_l 
        -
        \frac{1}{2}
        \sum_{k<l}(\partial_{x_k}\partial_{y_l}f-\partial_{y_k}\partial_{x_l}f)
        (dx_k\wedge dx_l + dy_k\wedge dy_l).
    \end{align*}
    Thus  we have
    \begin{align*}
        |i\partial\bar{\partial}f|^2
        &=
        \frac{1}{4} 
        \sum_{k,l}\Bigl[(\partial_{x_k}\partial_{x_l}f+\partial_{y_k}\partial_{y_l}f)^2
        +
        (\partial_{x_k}\partial_{y_l}f-\partial_{y_k}\partial_{x_l}f)^2\Bigr]
        \\
        &\leq
        \frac{1}{2}  
        \sum_{k,l}\Bigl[(\partial_{x_k}\partial_{x_l}f)^2
        +        (\partial_{y_k}\partial_{y_l}f)^2
        +       2 (\partial_{x_k}\partial_{y_l}f)^2\Bigr]
        \\
        &\leq  |\nabla df|^2.
    \end{align*}
    Now if we apply this in holomorphic normal coordinates centred at a certain point, we get the pointwise inequality \eqref{equ:complexHessian-realHessian}. If we apply it in $U=B$ we get $[i\partial\bar{\partial}f]_{0;B} \leq  [f]_{2,0;B}$.
    
    Similarly, if $q$ is another point in the coordinate chart we have
    \begin{align*}
        |i\partial\bar{\partial}f(p)-i\partial\bar{\partial}f(q)|^2
        =
        &\frac{1}{4} 
        \sum_{k,l}\Bigl[\{(\partial_{x_k}\partial_{x_l}f(p)-\partial_{x_k}\partial_{x_l}f(q))+(\partial_{y_k}\partial_{y_l}f(p)-\partial_{y_k}\partial_{y_l}f(q))\}^2
        \\
        &+        \{(\partial_{x_k}\partial_{y_l}f(p)-\partial_{x_k}\partial_{y_l}f(q))-(\partial_{y_k}\partial_{x_l}f(p)-\partial_{y_k}\partial_{x_l}f(q))\}^2\Bigr] 
        \\
        \leq&
        \frac{1}{2}  
        \sum_{k,l}\Bigl[(\partial_{x_k}\partial_{x_l}f(p)-\partial_{x_k}\partial_{x_l}f(q))^2
        +        (\partial_{y_k}\partial_{y_l}f(p)-\partial_{y_k}\partial_{y_l}f(q))^2
        \\
        &+        2(\partial_{x_k}\partial_{y_l}f(p)-\partial_{x_k}\partial_{y_l}f(q))^2\Bigr]
        \\
        \leq&
        |\nabla df(p)-\nabla df(q)|^2.
        \end{align*}
        The statement is a direct consequence.
    \end{proof}

\begin{proposition}
\label{lemma:nonlinear-estimate-holomprphic-atlas}
    For $u,v\in\mathscr{X}$, we have
     \begin{align*}
        \|{Q(u)-Q(v)}_{\mathscr{Y}} \leq q_{\mathrm{nl}}\|{u-v}_{\mathscr{X}} (\|{u}_{\mathscr{X}}+\|{v}_{\mathscr{X}}),
    \end{align*}
    where
    \begin{align*}
        q_{\mathrm{nl}} = \frac{K_0^2}{2}(2+\sqrt{2}K_0^3S)+\frac{K_0}{2}.
    \end{align*}
\end{proposition}
\begin{proof}
From \cref{lemma:estimate-nonlinear-term,lemma:complexHessian-realHessian},
we get
\begin{align*}
     \|{Q(u)-Q(v)}'_{C^{0,\alpha}(M)} 
     &\leq   \frac{K_0^2}{2} (2+\sqrt{2}K_0^3 S) R^{-2} \|{u-v}'_{C^{2,\alpha}(M)} R^{-2} \|{u+v}'_{C^{2,\alpha}(M)}
     \\
     &\leq   \frac{K_0^2}{2} (2+\sqrt{2}K_0^3 S)  \|{u-v}_{\mathscr{X}} \|{u+v}_{\mathscr{X}}.
\end{align*}
Writing $\gamma=i\partial\bar{\partial}u$ and $\eta=i\partial\bar{\partial}v$, we have by \eqref{equ:complexHessian-realHessian}
\begin{align*}
    |Q(u)-Q(v)| \leq \frac{1}{2} |\gamma-\eta|_g |\gamma+\eta|_g
    \leq
    \frac{1}{2} |i\partial\bar{\partial}(u-v)|_g |\nabla^g d (u+v)|_g,
\end{align*}
which gives with the proof of Lemma \ref{lemma:complexHessian-realHessian}
\begin{align*}
    aR^{-2-\frac{n}{q}}\|{Q(u)-Q(v)}_{L^2(M,g)}
    &\leq 
    \frac{1}{2}  K_0 \sup_p [(\phi_p^{-1})^*(\gamma-\eta)]_{0;B_R} aR^{-2-\frac{n}{q}} \|{\nabla^g d (u+v)}_{L^2(M,g)}
    \\
    &\leq
    \frac{1}{2}K_0 R^{-2} \|{u-v}'_{C^{2,\alpha}(M)}  \|{u+v}_{\mathscr{X}}
    \\
    &\leq 
    \frac{1}{2} K_0 \|{u-v}_{\mathscr{X}} \|{u+v}_{\mathscr{X}}.
\end{align*}

Adding the two estimates gives the claim.
\end{proof}

\subsection{Summary}

In this section we put together the pieces of this chapter to 
apply the Banach fixed point method \ref{lemma:inverse-function-thm} to perturb our approximate metric to a Calabi--Yau metric.

\begin{proof}[Proof of \cref{theorem:main}]
To fix our function spaces $\mathscr{X}$ and $\mathscr{Y}$,
which we introduced in section \ref{section:functionspaces}, we choose $\alpha=\frac{1}{2}$ and need to fix $R\leq R_0$, the radius of our coordinate balls, and the constant $a > 0$. 
We define $R$ close to the maximal $R$ which satisfies all three conditions \eqref{radius-cond-1}, \ref{radius-cond-2} and \ref{smallness-condition} for $R$.  For our approximate metric the most stringent is \ref{smallness-condition}, which is required to make the rearrangement in the global Schauder estimate work. 
We therefore choose
\begin{align*}
    R=7.85\cdot10^{-6}.
\end{align*}
In the proof of Proposition \ref{proposition:Laplace-inverse}
we used $a \geq \Cr{L22-Sobolev-update}$, the constant from the Sobolev embedding $ L^2_2(M) \hookrightarrow L^q(M)$ from Corollary \ref{cor:L22-Sobolev-embedding-update}, but we have not specified $q$ yet. By the definition of the $\mathscr{Y}$-norm,
the radius $R$ enters with exponent $2+\frac{n}{q}$. We want this exponent to be as small as possible while keeping the constant from the Sobolev embedding theorem reasonable. We will use $q=24$ which allows to choose $a=44$.

For $\omega$ the bounds for the inverse of the Laplacian from Proposition \cref{proposition:Laplace-inverse} and the non-linear term of $P$ from \cref{lemma:nonlinear-estimate-holomprphic-atlas} are
\begin{align*}
    C_{\Delta^{-1}}\leq153,
    \qquad q_{\mathrm{nl}}\leq1.503.
\end{align*}
Because the linearisation of $P$ is $\frac{1}{2}\Delta_g$, we can set $C = 2 C_{\Delta^{-1}}$ in Lemma \ref{lemma:inverse-function-thm}. Setting $t=\frac{1}{2Cq_{\mathrm{nl}}}$ in Lemma \ref{lemma:inverse-function-thm}, we need to check the following inequality:
\begin{align*}
   \|{f}_{\mathscr{Y}}
       < 
       \frac{1}{4 q_{\mathrm{nl}} C^2} 
       \geq 1.7763\cdot10^{-6}.
\end{align*}
Therefore, by \cref{prop-residual-bound} it is enough to check
\begin{align}
        \|{f}_{C^0(M)}
        <
        \chi_1 \frac{1}{4q_{\mathrm{nl}} C^2},
        \quad
        \sup_{M} |\nabla f|_g \leq 2^{\alpha-1} K_0^{-\frac{1}{2}} R^{-1} \chi_2\frac{1}{4q_{\mathrm{nl}} C^2},
        \quad
        \|{f}_{L^2(M)} \leq \frac{1}{a} R^{2+\frac{n}{q}} \chi_3\frac{1}{4q_{\mathrm{nl}} C^2}
        \label{required-residual}
\end{align}
for some $\chi_1,\chi_2,\chi_3 > 0$ with $\chi_1 + \chi_2 + \chi_3 =1$. Because the $L^2$-bound has the highest power of $R$ we expect it to be the most stringent bound and test this criterion with $\chi_1 = \chi_2 = \frac{1}{200}$ and $\chi_3 = \frac{99}{100}$. Substituting the values the required residual bounds are
\begin{align*}
     \|{f}_{C^0(M)}
        \leq
        8.88 \cdot 10^{-9},
        \quad
        \sup_{M} |\nabla f|_g \leq 7.99\cdot 10^{-4},
        \quad
        \|{f}_{L^2(M)} \leq  3.47\cdot 10^{-19}.
\end{align*}
The certified residual bounds from Proposition \ref{proposition:certified-residual} clearly meet this.
Here we estimate the $L^2$-norm by the $C^0$-norm using $\operatorname{vol}(M,\mu_{\omega})=1$.
By the definition of $\mathscr Y$ and Proposition \ref{prop-residual-bound}, we obtain
\begin{align*}
    \|{f}_{\mathscr Y}
    &\leq
    \left(1+aR^{-2-n/q}\right)\|{f}_{C^0(M)}
    +2^{1-\alpha}R K_0^{\frac12}\sup_M|\nabla f|_g
    \\
    &\leq 2.533\cdot10^{-7}.
\end{align*}
Thus Lemma \ref{lemma:inverse-function-thm} gives a solution $u$ such that
\begin{align*}
    \frac{1}{m!} (\omega + i\partial\bar{\partial}u)^m = \hat{\mu}
\end{align*}
with
\begin{align}
    \|{u}_{\mathscr{X}} \leq 2C \|{f}_{\mathscr{Y}}
    \leq 1.551\cdot10^{-4}.
    \label{error-bound}
\end{align}

We next show that $\omega+ i\partial\bar{\partial}u$ is positive definite.

    On $B_R$ in a coordinate chart $\mathcal{B}_p$ write $g_{\mathcal{B}_p} + g_{u,p}$ for the bilinear form defined by $\omega+ i\partial\bar{\partial}u$. Analogously to Lemma \ref{lemma:complexHessian-realHessian} we have $|g_{u,p}| \leq |D^2 (u\circ \phi_p^{-1})|$. Because the eigenvalues of $g_{u,p}$ appear in pairs we have $|g_{u,p}|_{\operatorname{op}} \leq \frac{1}{\sqrt{2}} |g_{u,p}|$. 
    By Lemma \ref{lemma:eigenvalues-metric}, $g_{\mathcal{B}_p}\geq K_0^{-1}I$.
    Therefore by \eqref{error-bound}
    \begin{align*}
        g_{\mathcal{B}_p} + g_{u,p} &\geq K_0^{-1} - |g_{u,p}|_{\operatorname{op}} 
        \geq  K_0^{-1} - \frac{1}{\sqrt{2}} |D^2(u\circ \phi_p^{-1})|
        \geq K_0^{-1} - \frac{1}{\sqrt{2}} \|{u}_{\mathscr{X}} \\
        &\geq
        K_0^{-1} - \sqrt{2}C \|{f}_{\mathscr{Y}} 
        \geq 
        K_0^{-1} - \frac{1.551\cdot10^{-4}}{\sqrt2}>0.9992.
\end{align*}
By \cref{lemma:eigenvalues-metric,lemma:complexHessian-realHessian},
we have
$\|{i\partial\bar\partial u}_{C^0(M,g)}
\leq K_0\|{u}_{\mathscr X}<1.553\cdot10^{-4}$.
\end{proof}

We conclude this section by explaining our motivation for the definition of the norms on $\mathscr{X}$ and $\mathscr{Y}$. 
We could also obtain bounds $C$ and $q_{\mathrm{nl}}$ if instead we take $\|{u}_{\mathscr{X}} = R^{-2}\|{u}_{C^{2,\alpha}(M)} $ and $\|{u}_{\mathscr{Y}} = \|{u}_{C^{0,\alpha}(M)}$, but the residual requirement gets significantly worse. Schematically, the estimate for $\Delta_g^{-1}$ then is
\begin{align*}
    R^{-2}\|{u}_{C^{2,\alpha}(M)} &\lesssim \|{\Delta_g u}_{C^{0,\alpha}(M)} + R^{-2}\|{u}_{C^0(M)}
    \\
    &\lesssim \|{\Delta_g u}_{C^{0,\alpha}(M)} + R^{-2} (\|{\Delta_g u}_{C^0(M)} + \|{u}_{L^2(M)})
    \\
    &\lesssim
     \|{\Delta_g u}_{C^{0,\alpha}(M)} + R^{-2} \|{\Delta_g u}_{C^0(M)}.
\end{align*}
Note that on the right hand side  $\|{\Delta_g u}_{C^0(M)}$
appears once with a factor $R^{-2}$ and once without the factor. Thus we need to include the factor $R^{-2}$ into $C$. Because $C$ appears quadratically in the residual requirement,
we then need $\|{f}_{C^0(M)} \lesssim R^4$, which is significantly worse. The best improvement for the exponent of $R$ in the residual requirement is $R^2$ if in Proposition \ref{proposition:Laplace-inverse} we embed $L^2_3 \hookrightarrow L^{\infty} $ instead. Then we would need to take the $L^2_3$-norm in the definition of $\|{\, \cdot\,}_{\mathscr{X}}$ and the $L^2_1$-norm in the definition of $\|{\, \cdot\,}_{\mathscr{Y}}$. However, for the elliptic estimate of $\Delta_g: L^2_3 \rightarrow L^2_1$ we need a bound on the sectional curvature which cannot be small by the Gauss--Bonnet theorem. Thus the marginal gain in the exponent of $R$ is not justified.

Next we explain our reason for including the factor $a$ in the definition of the norms. In Proposition \ref{proposition:Laplace-inverse} setting $a = \Cr{L22-Sobolev-update}$ equal to the Sobolev embedding constant can improve $C$ by up to a factor of $a$. We have to pay for this by introducing a factor $\frac{1}{a}$ in the required $L^2$-bound in \eqref{required-residual}. Because $q_{\mathrm{nl}}$ does not depend on $a$ and $C$ appears quadratically in the residual bound, overall this can improve the required $L^2$-bound by up to a factor of $a$.

\appendix

\section{Analysis on $\mathbb{R}^n$}

\subsection{Mean value inequality}

\begin{lemma}
\label{lemma-meanvalue}

    Let $1\leq q<\infty$ and $u\in C^2(B_r(0))$. Then we have
    \begin{align*}
        u(0) \leq \left( \fint_{B_r} |u|^q dx \right)^{\frac{1}{q}} + \frac{r^2}{2n} \sup |\Delta u|.
    \end{align*}
\end{lemma}
\begin{proof}
    Set $w(x)=r^2-|x|^2$. Then $\tilde{u}:= u-\frac{\sup|\Delta u|}{2n} w$ is a sub-harmonic function. By the mean value theorem we have
    \begin{align*}
        u(0)-\frac{r^2}{2n} \sup |\Delta u| = \tilde{u}(0) \leq \fint_{B_r} \tilde{u} dx \leq \fint_{B_r} u dx, 
    \end{align*}
    as $w\geq 0$ on $B_r$. The statement follows from Hölder's inequality.
\end{proof}

\subsection{Schauder estimate}

In this section we prove the Schauder interior estimate on a ball in $\mathbb{R}^n, n\geq 3$, with explicit constants.
For $f\in C^{0,\alpha}(\overline{B_2})$, we define its Newtonian potential as
\begin{align*}
    w(x) = \int_{B_2} \Gamma(x-y) f(y) dy,
\end{align*}
where $\Gamma$ is the fundamental solution of the Laplace operator given by
\begin{align*}
\Gamma(x-y) = \frac{1}{(2-n) \omega_{n-1}} |x-y|^{2-n}.   
\end{align*}
Throughout this section $\omega_{n-1}=|S^{n-1}|$ denotes the volume of the $(n-1)$-dimensional unit sphere.

Writing
$\operatorname{sym}(\theta\otimes\operatorname{Id})_{ijk}
= \theta_i\delta_{jk}+\theta_j\delta_{ik}+\theta_k\delta_{ij}$
we note here the first three derivatives of $\Gamma(x-y)$:
\begin{align*}
    D\Gamma(x-y) &= \frac{1}{\omega_{n-1}} (x-y) |x-y|^{-n},
    \\
    D^2\Gamma(x-y) &= \frac{1}{\omega_{n-1}} \left(\operatorname{Id} -n |x-y|^{-2} (x-y)^{\otimes 2}  \right) |x-y|^{-n},
    \\
    D^3\Gamma(x-y) &= \frac{n}{\omega_{n-1}} \{(n+2) |x-y|^{-3} (x-y)^{\otimes 3}  - |x-y|^{-1}  \operatorname{sym}((x-y) \otimes \operatorname{Id})\} |x-y|^{-n-1}.
\end{align*}

For $A,B\in\mathbb{R}^{n\times n}$ denote by $\langle A, B\rangle = \operatorname{Tr}(A B^T)$ their inner product. In the proof of the Schauder estimate we will bound $|D^2 w(x)|$ by pairing $D^2w(x)$ with the symmetric matrix $A=D^2 w(x)/|D^2 w(x)$. In polar coordinates $y = x + \rho \theta$, where $\rho > 0$ and $\theta\in S^{n-1}$, we have
\begin{align}
    \langle D^2 \Gamma (x-y), A\rangle = - \frac{n}{\omega_{n-1}} (\theta^T A_0 \theta) \rho^{-n},
\end{align}
where $A_0$ is the trace-free part of $A$.
$D^3 \Gamma(x-y)$ we will pair with $a\in \mathbb{R}^n$ and a symmetric matrix $A$ and in polar coordinates get
\begin{align}
    D^3 \Gamma (x-y)[a,A] = - \frac{n}{\omega_{n-1}} \{ (n+2) \theta^{\otimes 3} - \operatorname{sym}(\theta\otimes \operatorname{Id})\} [a,A] \rho^{-n-1}
 \end{align}

Denote by $d\sigma_{n-1}$ the normalised surface measure on $S^{n-1}$. Denote by $\theta_1, \dots, \theta_n$ the restriction of the standard coordinates of $\mathbb{R}^n$ to $S^{n-1}$. We have the following identities:
\begin{subequations}
\begin{align}
    \int_{S^{n-1}} \theta_i \theta_j\, d\sigma_{n-1}(\theta) &= \frac{1}{n} \delta_{ij},
    \label{moment2}
    \\
    \int_{S^{n-1}}  \theta_i \theta_j \theta_k \theta_l \, d\sigma_{n-1}(\theta) &= \frac{1}{n(n+2)} (\delta_{ij}\delta_{kl} + \delta_{ik}\delta_{jl}+ \delta_{il}\delta_{jk}),
    \label{moment4}
    \\
     \int_{S^{n-1}} \theta_{i_1}\theta_{i_2}\theta_{i_3}\theta_{i_4}\theta_{i_5}\theta_{i_6}\,d\sigma_{n-1}(\theta)
&= \frac{1}{n(n+2)(n+4)}
\sum_{\substack{
    \{\{a_1,a_2\},\{a_3,a_4\},\{a_5,a_6\}\}\\
    \text{a partition of }\{1,\dots,6\}
}}
\delta_{i_{a_1}i_{a_2}}\delta_{i_{a_3}i_{a_4}}\delta_{i_{a_5}i_{a_6}}.
     \label{moment6}
\end{align}
\end{subequations}

\begin{lemma}
\label{lemma:integral-id}
\begin{enumerate}[label= (\roman*), nosep, wide]
    \item Suppose that $A$ is symmetric and trace-free. Then
    \begin{align*}
        \int_{S^{n-1}} (\theta^TA\theta)^2 \, d\sigma_{n-1}(\theta)  = \frac{2}{n(n+2)} |A|^2.
    \end{align*}
    \item
    If $A$ is a symmetric $n\times n$ matrix and $a$ a unit length vector, then 
    \begin{align*}
        \int_{S^{n-1}} \left(\{(n+2) \theta^{\otimes 3} - \operatorname{sym}(\theta \otimes\operatorname{Id} )\}[a,A] \right)^2 d\sigma_{n-1}(\theta) 
        \leq \frac{6}{n+4} |A|^2.
    \end{align*}
    \end{enumerate}
\end{lemma}
\begin{proof}
Proof of (i):
Now we have
\begin{align*}
\int_{S^{n-1}} (\theta^TA\theta)^2 \, d\sigma_{n-1}(\theta) 
&= \sum_{i,j,k,l} A_{ij} A_{kl}   \int_{S^{n-1}} \theta_{i}\theta_{j}\theta_{k}\theta_{l} \, d\sigma_{n-1}(\theta)
\\
&=\frac{1}{n(n+2)} \sum_{i,j,k,l} A_{ij} A_{kl} (\delta_{ij}\delta_{kl} + \delta_{ik}\delta_{jl}+ \delta_{il}\delta_{jk}) 
\\
&= \frac{1}{n(n+2)} \left(\sum_{ik} A_{ii} A_{kk} + \sum_{ij} A_{ij} A_{ij} + \sum_{i,j} A_{ij} A_{ji} \right) 
\\
&= \frac{2}{n(n+2)} |A|^2.
\end{align*}
\\
Proof of (ii):
The tensor $\{(n+2)\theta^{\otimes3}-\operatorname{sym}(\theta\otimes\operatorname{Id})\}$ annihilates the scalar part of $A$. Since removing this part does not increase $|A|$, we may assume without loss of generality that $A$ is trace-free.
Because $A$ is symmetric we have
\begin{align*}
\{(n+2) \theta^{\otimes 3} - \operatorname{sym}(\theta \otimes\operatorname{Id} )\}[a,A] = (n+2) (\theta^T a) (\theta^T A \theta) - 2 (\theta^T A a).     
\end{align*}
We have
\begin{align*}
    ((n+2) (\theta^T a) (\theta^T A \theta) - 2 (\theta^T A a))^2 =&\, (n+2)^2 (\theta^Ta)^2 (\theta^TA\theta)^2 
    \\ & -4 (n+2) (\theta^Ta) (\theta^TA\theta) (\theta^T A a) + 4 (\theta^T A a)^2.
\end{align*}
By \eqref{moment2} we have
\begin{align*}
\int_{S^{n-1}}  (\theta^T A a)^2 \, d\sigma_{n-1}(\theta) = \sum_{i,j,k,l} A_{ij} a_j A_{kl} a_l \int_{S^{n-1}}  \theta_i \theta_k \, d\sigma_{n-1}(\theta)    
=
\frac{1}{n} \sum_{i,j,k,l} A_{ij} a_j A_{kl} a_l \delta_{ik} = \frac{1}{n} |Aa|^2.
\end{align*}
By \eqref{moment4} we have
\begin{gather*}
\int_{S^{n-1}} (\theta^Ta) (\theta^TA\theta) (\theta^T A a) \, d\sigma_{n-1}(\theta) 
= \sum_{i,j,k,l,p} A_{jk} A_{lp} a_i a_p \int_{S^{n-1}}  \theta_i \theta_j \theta_k \theta_l \, d\sigma_{n-1}(\theta)
\\
=
\frac{1}{n(n+2)} \sum_{i,j,k,l,p} A_{jk} A_{lp} a_i a_p (\delta_{ij}\delta_{kl} + \delta_{ik}\delta_{jl}+ \delta_{il}\delta_{jk}) = \frac{2}{n(n+2)} |Aa|^2.
\end{gather*}
By \eqref{moment6} we have
\begin{align*}
    \int_{S^{n-1}} (\theta^Ta)^2 (\theta^TA\theta)^2 \, d\sigma_{n-1}(\theta) &= \sum_{i,j,k,l,p,q} a_i a_j A_{kl} A_{pq} \int_{S^{n-1}} \theta_i \theta_j \theta_k \theta_l \theta_p \theta_q \, d\sigma_{n-1}(\theta)
    \\
    &= \frac{1}{n(n+2)(n+4)} (2 |A|^2 |a|^2 + 8 |Aa|^2).
\end{align*}
Thus we have
\begin{align*}
 \int_{S^{n-1}} ((n+2) (\theta^T a) (\theta^T A \theta) - 2 (\theta^T A a))^2 \, d\sigma_{n-1}(\theta) = \frac{2(n+2)}{n(n+4)} |A|^2 + \frac{4n}{n(n+4)} |Aa|^2.
\end{align*}
Because $A$ is trace-free and symmetric we have $|A|_{\mathrm{op}} \leq \sqrt{\frac{n-1}{n}} |A|$.  The statement follows.  
\end{proof}

\begin{lemma}
\label{lemma:integral-id-2}
    If $x\in B_1$, then
    \begin{align*}
        \frac{1}{\omega_{n-1}} \int_{\partial B_2} |x-y|^{1-n} ds_y \leq M_n,
    \end{align*}
    where 
    \begin{align*}
        M_n = 2^{n-1} \frac{\omega_{n-2}}{\omega_{n-1}} \int_{-1}^1 \frac{(1-t^2)^{\frac{n-3}{2}}}{(5-4t)^{\frac{n-1}{2}}} dt.
    \end{align*}
\end{lemma}
\begin{remark}
    Using rigorous quadrature from \cite{python-flint}, we compute $M_4 < 1.112$. 
\end{remark}
\begin{proof}
    For $x\in B_1$ write $F(x) = \int_{\partial B_2} |x-y|^{1-n} ds_y$. The integrand is smooth and satisfies $\Delta_x |x-y|^{1-n} = (n-1) |x-y|^{-n-1}$. Therefore, $F$ is subharmonic and attains its maximum on the boundary of $B_1$. By rotational symmetry we can assume $x = e_n$. 
    First write $y=2\theta$ for $\theta\in S^{n-1}$. 
    We can write $\theta = ((1-t^2)^{\frac{1}{2}} \xi, t) $ with $t\in (-1,1)$ and $\xi\in S^{n-2}$. We have $|x-y|^2 = 5-4t $ and thus
    \begin{align*}
        \int_{\partial B_2} |e_n-y|^{1-n} ds_y = 2^{n-1} \int_{-1}^1 \int_{S^{n-2}} (5-4t)^{\frac{1-n}{2}} ds_{S^{n-1}}(t,\xi).
    \end{align*}
    We are integrating over the surface levels of the function $g(\theta) = \theta_n$ and each surface level is a $(n-2)-$dimensional sphere of radius $(1-t^2)^{\frac{1}{2}}$. By the co-area formula the volume form is given by
    \begin{align*}
        ds_{S^{n-1}}(t,\xi) = (1-t^2)^{\frac{n-2}{2}} \frac{1}{|\nabla^{S^{n-1}} g(\theta)|} dt ds_{S^{n-2}}(\xi).
    \end{align*}
    $\nabla^{S^{n-1}} g(\theta)$ is the projection of $\nabla g$ onto the tangent planes of $S^{n-1}$ so we have
    \begin{align*}
        \nabla^{S^{n-1}} g(\theta) = \nabla g(\theta) - \langle \nabla g(\theta),\theta \rangle \theta  = e_n - \langle e_n,\theta\rangle \theta 
    \end{align*}
    and thus
    \begin{align*}
        | \nabla^{S^{n-1}} g(\theta)|^2 = 1 - \theta_n^2 = 1-t^2.
    \end{align*}
    The statement follows.
\end{proof}

\begin{lemma}
\label{lemma:integral-id-3}
    Let $e\in S^{n-1}$. Then
    \begin{align*}
        \int_{S^{n-1}} |\langle e,\theta\rangle| d\sigma_{n-1} = \frac{2}{n-1} \frac{\omega_{n-2}}{\omega_{n-1}}.
    \end{align*}
\end{lemma}
\begin{proof}
    By rotational symmetry we can assume $e = e_n$. 
    As in the proof of Lemma \ref{lemma:integral-id-2} we integrate over the level sets of the function $g(\theta) = \theta_n$ and get
    \begin{align*}
     \int_{S^{n-1}} |\theta_n| d\sigma_{n-1}(\theta) = \frac{\omega_{n-2}}{\omega_{n-1}} \int_{-1}^1 |t| (1-t^2)^{\frac{n-3}{2}} dt = 2   \frac{\omega_{n-2}}{\omega_{n-1}} \int_{0}^1 t (1-t^2)^{\frac{n-3}{2}} dt.
     \end{align*}
     We use the substitution $t = \sin \phi$ to compute
     \begin{align*}
     \int_{0}^1 t (1-t^2)^{\frac{n-3}{2}} dt = \int_0^{\frac{\pi}{2}} \sin \phi (\cos \phi)^{n-2} d\phi = \left.- \frac{1}{n-1} (\cos\phi)^{n-1}\right|^{\phi=\frac{\pi}{2}}_{\phi=0} = \frac{1}{n-1}.   
     \end{align*}
     The statement follows.
\end{proof}

\begin{lemma}
\label{lemma:integral-Newton-kernel}
    For $n\geq 3$ and $x\in \partial B_2$ we have
    \begin{align*}
        \int_{B_2} |\Gamma(x-y)| dy = \frac{1}{(n-2)\omega_{n-1}} \int_{B_2} |x-y|^{2-n} dy = \frac{4}{n(n-2)}.
    \end{align*}
\end{lemma}
\begin{proof}
    Write $F(x)=\int_{B_2} |\Gamma(x-y)| dy $. This function is radially symmetric. Write $r=|x|$. Differentiating the integral gives
    \begin{align*}
        F'(r) = - \frac{1}{\omega_{n-1}} \int_{B_2} |r e-y|^{-n} (re-y)\cdot e dy,
    \end{align*}
    where $e$ is any unit vector. For $r=0$ the integral vanishes because the integrand is odd with respect to the symmetry $y\mapsto -y$. Thus $F'(0)=0$.

    Because $F$ is radially symmetric we have $\Delta F = F''(r) + \frac{n-1}{r} F'(r)$. Because $\Gamma$ is the fundamental solution of the Laplace operator we have $\Delta F = -1$. Therefore we have
    \begin{align*}
        F''(r)+ \frac{n-1}{r} F'(r) = -1,
    \end{align*}
    and together with $F'(0)=0$ this gives $F'(r) = -\frac{1}{n} r$ and $F(r) = F(0) - \frac{1}{2n} r^2 = \frac{2}{n-2}-\frac{1}{2n}r^2$ and $F(2)=\frac{4}{n(n-2)}$.
\end{proof}

\begin{lemma}
Let $w$ be the Newton potential of $f$ in $B_2$. 
\begin{enumerate}[label= (\roman*), nosep, wide]
\item We have 
    \begin{align*}
        [w]_{1,0;B_1} \leq \Cl{Schauder-1-0} [f]_{0;B_2},
    \end{align*}
    where
    \begin{align*}
        \Cr{Schauder-1-0} = \frac{6}{n-1} \frac{\omega_{n-2}}{\omega_{n-1}}.
    \end{align*}
    \item
    We have
    \begin{align*}
        [w]_{2,0;B_1} \leq \Cl{Schauder-2-0-1} [f]_{0;B_2} + \Cl{Schauder-2-0-2} [f]_{0,\alpha;B_2},
    \end{align*}
    where
    \begin{align*}
        \Cr{Schauder-2-0-1} = M_n, \quad \Cr{Schauder-2-0-2} = \frac{n}{\alpha} 3^{\alpha} \left( \frac{2}{n(n+2)}\right)^{\frac{1}{2}}
        \end{align*}
    \item 
    We have
    \begin{align*}
        [w]_{2,\alpha;B_1} \leq \Cl{Schauder-2-alpha-1} [f]_{0;B_2} + \Cl{Schauder-2-alpha-2} [f]_{0,\alpha;B_2},
    \end{align*}
    where
    \begin{align*}
        \Cr{Schauder-2-alpha-1} &= 2^{1-\alpha} (n-1) M_n, \\ \Cr{Schauder-2-alpha-2} &= 3 M_n + 2 \frac{n}{\alpha} \left( \frac{3}{2}\right)^{\alpha} \left( \frac{2}{n(n+2)}\right)^{\frac{1}{2}} + n \frac{1}{1-\alpha} 3^{\alpha} 2^{1-\alpha} \left( \frac{6}{n+4}\right)^{\frac{1}{2}}.
        \end{align*}
\end{enumerate}
\end{lemma}
\begin{proof}
Proof of (i):
Let $x\in B_1$ with $Dw(x)\neq 0$.
Set $a=Dw(x)/|Dw(x)|$.
By \cite[Lemma 4.1]{Gilbarg1998} we have
\begin{align*}
    \langle D w(x),a \rangle = \int_{B_2} \langle D\Gamma(x-y),a\rangle f(y) dy.
\end{align*}
Thus by Lemma \ref{lemma:integral-id-3}
\begin{align*}
    |Dw(x)| \leq [f]_{0;B_2} \int_{B_3(x)} |\langle D\Gamma(x-y),a\rangle| dy
    = [f]_{0;B_2} \int_0^3 d\rho  \int_{S^{n-1}} |\theta_1| d\sigma_{n-1}(\theta)  
    \leq \frac{6}{n-1} \frac{\omega_{n-2}}{\omega_{n-1}} [f]_{0;B_2}.
\end{align*}
Proof of (ii): 
Let $x\in B_1$ with $D^2w(x)\neq 0$. Set $A=D^2w(x)/|D^2w(x)|$.
In particular $A$ is symmetric with $|A|=1$.
By \cite[Lemma 4.2]{Gilbarg1998} we have
\begin{align*}
    \langle D^2w(x),A \rangle = \int_{B_2} \langle D^2\Gamma(x-y),A \rangle (f(y)-f(x)) dy - f(x) \int_{\partial B_2} D\Gamma(x-y)^T A \nu(y) ds_y.
\end{align*}
By Lemma \ref{lemma:integral-id} we have
\begin{align*}
\int_{B_2} |\langle D^2\Gamma(x-y),A \rangle (f(y)-f(x))| dy &\leq  [f]_{0,\alpha;B_2}  \int_{B_3(x)} |\langle D^2\Gamma(x-y),A \rangle| |x-y|^{\alpha} dy
\\
&=
n[f]_{0,\alpha;B_2} \int_0^3 \rho^{\alpha-1} d\rho  \int_{S^{n-1}} |\theta^T A_0 \theta| d\sigma_{n-1}(\theta) 
\\
&\leq 
n[f]_{0,\alpha;B_2} \frac{1}{\alpha} 3^{\alpha} \left( \frac{2}{n(n+2)}\right)^{\frac{1}{2}}
\end{align*}
and by Lemma \ref{lemma:integral-id-2}
\begin{align*}
            \int_{\partial B_2} |D\Gamma(x-y)^T A \nu(y)| ds_y \leq \frac{1}{\omega_{n-1}} |A| \int_{\partial B_2}|x-y|^{1-n} \leq M_n |A|.
\end{align*}
\\
Proof of (iii):
Let $x,\bar{x}\in B_1$ be distinct and write $\delta= |x-\bar{x}| $ and $\xi = \frac{1}{2}(x+\bar{x}) $.
    As in the proof of \cite[Lemma 4.4]{Gilbarg1998} we write
    \begin{align*}
        D^2w(\bar{x})-D^2w(x) = f(x) I_1 + (f(x)-f(\bar{x})) I_2 + I_3-I_4+ (f(x)-f(\bar{x})) I_5 + I_6,
    \end{align*}
    where
    \begin{align*}
        I_1 &= \int_{\partial B_2} (D\Gamma(x-y) - D\Gamma(\bar{x}-y)) \otimes \nu(y) ds_y,
        \\
        I_2 &= \int_{\partial B_2} D\Gamma(\bar{x}-y) \otimes \nu(y) ds_y,
        \\
        I_3 &= \int_{B_2\cap B(\xi,\delta)} D^2\Gamma(x-y) (f(x)-f(y)) dy,
        \\
        I_4 &= \int_{B_2\cap B(\xi,\delta)} D^2\Gamma(\bar{x}-y) (f(\bar{x})-f(y)) dy,
        \\
        I_5 &= \int_{B_2 - B(\xi,\delta)} D^2\Gamma(x-y) dy,
        \\
        I_6 &= \int_{B_2 - B(\xi,\delta)} (D^2\Gamma(x-y) - D^2\Gamma(\bar{x}-y))(f(\bar{x})-f(y))dy.
    \end{align*}
    Assume without loss of generality that $D^2w(x) \neq D^2w(\bar{x})$ and set $A=(D^2w(x)-D^2w(\bar{x}))/|D^2w(x)-D^2w(\bar{x})|$.
    \\
    \textbf{Estimate of $I_1$:}
    There exists $\hat{x}$ on the line segment between $x$ and $\bar{x}$ such that
    \begin{align*}
        |\langle I_1,A\rangle| \leq |x-\bar{x}| \cdot |A| \int_{\partial B_2} |D^2\Gamma(\hat{x}-y)[a,\, \cdot\, ]| ds_y
    \end{align*}
    for the unit vector $a=(x-\bar{x})/|x-\bar{x}|$.
    We have
    \begin{align*}
        D^2 \Gamma(\hat{x}-y)[a,\,\cdot\,] = \frac{1}{\omega_{n-1}} |\hat{x}-y|^{-n} (a-n\langle a, \hat{x}-y\rangle (\hat{x}-y) |\hat{x}-y|^{-2}).
    \end{align*}
    For unit vectors $a,b$ we have
    \begin{align*}
        |n\langle a,b\rangle b - a| \leq n-1.
    \end{align*}
    Thus
    \begin{align*}
        |\langle I_1,A\rangle| &\leq |x-\bar{x}| \cdot |A| (n-1) \frac{1}{\omega_{n-1}} \int_{\partial B_2} |\hat{x}-y|^{-n} ds_y
        \\
        &\leq |x-\bar{x}| \cdot |A| (n-1) \frac{1}{\omega_{n-1}} \int_{\partial B_2} |\hat{x}-y|^{1-n} ds_y
        \\
        &\leq |x-\bar{x}| \cdot |A| (n-1) M_n \leq  2^{1-\alpha} (n-1) M_n \delta^{\alpha}.
    \end{align*}
    Here we have used $|\hat{x}-y| > 1$. The integral can probably be improved.
    \\
    \textbf{Estimate of $I_2$:} We have
    \begin{align*}
        |\langle I_2,A\rangle| \leq \int_{\partial B_2} |D\Gamma(\bar{x}-y)^T A \nu(y)| ds_y \leq \frac{1}{\omega_{n-1}} |A| \int_{\partial B_2}|\bar{x}-y|^{1-n} \leq M_n |A|.
    \end{align*}
    This can be improved.
    \\
    \textbf{Estimate of $I_3$:}
    Because we have $B(\xi, \delta) \subset B(x,\tfrac{3}{2}\delta)$ we have
        \begin{align*}
            |\langle I_3, A\rangle| &\leq [f]_{0,\alpha} \int_{B(x,\tfrac{3}{2}\delta)} |\langle D^2 \Gamma(x-y), A \rangle|\cdot |x-y|^\alpha dy 
            \\
            &\leq [f]_{0,\alpha} \int_0^{\frac{3}{2}\delta} \rho^{\alpha-1} d\rho  \cdot n \int_{S^{n-1}} |\theta^T A_0 \theta| d\sigma_{n-1}(\theta)
            \\
            &\leq [f]_{0,\alpha} \frac{n}{\alpha} \left( \frac{3}{2}\delta\right)^{\alpha} \left( \frac{2}{n(n+2)}\right)^{\frac{1}{2}} |A_0|,
        \end{align*}
        where in the last step we have used Hölder's inequality and Lemma \ref{lemma:integral-id}.
        \\
        \textbf{Estimate of $I_4$:} Analogous to $I_3$.
        \\
        \textbf{Estimate of $I_5$:} Write $V_A(y) = - A D\Gamma(x-y)$ so that $\operatorname{div} V_A(y) = \langle D^2\Gamma(x-y),A\rangle$. Then by the divergence theorem we have
        \begin{align*}
           |\langle I_5,A\rangle| &\leq \int_{\partial B_2} |\langle V_A(y),\nu(y)\rangle| ds_y +    \int_{\partial B_{\delta}(\xi)} |\langle V_A(y),\nu(y)\rangle| ds_y 
           \\
           &\leq |A| \frac{1}{\omega_{n-1}} \int_{\partial B_2} |x-y|^{1-n} ds_y +  |A| \frac{1}{\omega_{n-1}} \int_{\partial B_{\delta}(\xi)} |x-y|^{1-n} ds_y
           \\
           &\leq 2 M_n.
        \end{align*}
        \\
        \textbf{Estimate of $I_6$:} There exists $\hat{x}$ on the line between $x$ and $\bar{x}$ such that 
        \begin{align*}
            |\langle I_6,A\rangle| &\leq \delta \int_{B_2-B(\xi,\delta)} | D^3\Gamma(\hat{x}-y)[a,A_0]| |f(\bar{x})-f(y)| dy 
            \\
            &\leq \delta [f]_{0,\alpha;B_2} \int_{B_2-B(\hat{x},\delta/2)} | D^3\Gamma(\hat{x}-y)[a,A_0]| |\bar{x}-y|^{\alpha} dy
            \\
            &\leq \delta [f]_{0,\alpha;B_2} 3^{\alpha} \int_{B(\hat{x},3)-B(\hat{x},\delta/2)} | D^3\Gamma(\hat{x}-y)[a,A_0]| |\hat{x}-y|^{\alpha} dy
            \\
            &\leq \delta [f]_{0,\alpha;B_2} 3^{\alpha}n \int_{\frac{\delta}{2}}^3 \rho^{\alpha-2} d\rho \int_{S^{n-1}} |\{ (n+2)\theta^{\otimes 3}-\operatorname{sym}(\theta\otimes\operatorname{Id})\}[a,A_0]| d\sigma_{n-1}(\theta)
            \\
            &\leq \delta [f]_{0,\alpha;B_2} 3^{\alpha}n \frac{1}{1-\alpha} \left( \frac{\delta}{2}\right)^{\alpha-1} \left( \frac{6}{n+4}\right)^{\frac{1}{2}} |A_0|
            \\
            &\leq n \frac{1}{1-\alpha} 3^{\alpha} 2^{1-\alpha} \left( \frac{6}{n+4}\right)^{\frac{1}{2}} |A_0| \delta^{\alpha} [f]_{0,\alpha;B_2}.
        \end{align*}
\end{proof}

\begin{lemma}
    Let $v$ be a harmonic function on $B_{2R}$. Then
    \begin{align*}
        R [v]_{1,0;B_R} \leq \Cl{harmonic-1}  \|{v}_{C^0(B_{2R})},
        \\
        R^2 [v]_{2,0;B_R} \leq \Cl{harmonic-2}   \|{v}_{C^0(B_{2R})},
        \\
        R^{2+\alpha} [v]_{2,\alpha;B_R} \leq \Cl{harmonic-3} \|{v}_{C^0(B_{2R})}.
    \end{align*}
    where
    \begin{align*}
        \Cr{harmonic-1} &= \frac{2n}{n-1} \frac{\omega_{n-2}}{\omega_{n-1}},
        \\
        \Cr{harmonic-2} &= \sqrt{2n (n+2)},
        \\
        \Cr{harmonic-3} &= 2^{\frac{3}{2}-\alpha} 3^{\frac{\alpha}{2}} n^{\frac{1+\alpha}{2}} (n+2)^{\frac{1-\alpha}{2}} (n+4)^{\frac{\alpha}{2}}.
    \end{align*}
\end{lemma}    
\begin{proof}
Let $x\in B_R$. Applying Poisson's formula and its derivatives to a ball of radius $R$ centred at $x$ gives 
\begin{subequations}
\begin{align}
    R Dv(x) &=  n \int_{S^{n-1}} v(x+R\theta) \theta\,  d\sigma_{n-1}(\theta),
    \label{Poisson1}
    \\
    R^2 D^2v(x) &=  (n+2) \int_{S^{n-1}} v(x+R\theta) (n \theta \otimes \theta - \operatorname{Id})\, d\sigma_{n-1}(\theta),
    \label{Poisson2}
    \\
    R^3 D^3 v(x) &= n (n+4) \int_{S^{n-1}} v(x+R\theta) \{(n+2) \theta^{\otimes 3} - \operatorname{sym}(\theta \otimes\operatorname{Id} )\} \, d\sigma_{n-1}(\theta).
    \label{Poisson3}
\end{align}
\end{subequations}
\textbf{Estimate of the first derivative:} 
Let $a\in \mathbb{R}^n$ be a unit length vector. We test \eqref{Poisson1} with $a$ and get
\begin{align*}
    R |\langle D v(x), a \rangle| \leq n \sup_{B_{2R}} |v|  \int_{S^{n-1}} |\langle \theta, a \rangle| \, d\sigma_{n-1}(\theta)
    =
    n \sup_{B_{2R}} |v|  \int_{S^{n-1}} |\theta_1| \, d\sigma_{n-1}(\theta)
    =
     \frac{2n}{n-1} \frac{\omega_{n-2}}{\omega_{n-1}} \sup_{B_{2R}} |v|.
\end{align*}
Here we have used Lemma \ref{lemma:integral-id-3}. Applying this formula with $a = Dv(x)/|Dv(x)|$ gives the bound for the first derivative.
\\
\textbf{Estimate of the second derivative:} 
Let $A$ be a symmetric, trace-free $n\times n$-matrix with $|A|=1$. Testing \eqref{Poisson2} with $A$ gives
\begin{align*}
    R^2 |\operatorname{Tr}(D^2 v(x) A)| \leq (n+2) \sup_{B_{2R}} |v|  \int_{S^{n-1}}  |\operatorname{Tr}((n \theta \otimes \theta - \operatorname{Id})A)|\, d\sigma_{n-1}(\theta).
\end{align*}
Because of $\operatorname{Tr}(A)=0$ we have 
\begin{align*}
\operatorname{Tr}((n \theta \otimes \theta - \operatorname{Id})A) = n \theta^T A \theta - \operatorname{Tr}(A) =   n \theta^T A \theta.  
\end{align*}
Using Hölder's inequality gives
\begin{align*}
R^2 |\operatorname{Tr}(D^2 v(x) A)| \leq  n  (n+2) \sup_{B_{2R}} |v|   \left( \int_{S^{n-1}} (\theta^TA\theta)^2 \, d\sigma_{n-1}(\theta)\right)^{\frac{1}{2}}.
\end{align*}
Because $v$ is harmonic we have $\operatorname{Tr}(D^2v) = \Delta v = 0$ and thus we can apply Lemma \ref{lemma:integral-id} with $A = D^2v(x)/|D^2 v(x)|$ which gives the estimate for the second derivative.
\\
\textbf{Lipschitz bound for second derivative:} Let $a\in\mathbb{R}^n$ be a unit length vector and $A$ a symmetric, trace-free $n\times n$-matrix with $|A|=1$. Testing \eqref{Poisson3} with $(a,A)$ gives 
\begin{align*}
R^3 |D^3 v(x)[a,A]| \leq n (n+4) \sup_{B_{2R}} |v| \int_{S^{n-1}} |\{(n+2) \theta^{\otimes 3} - \operatorname{sym}(\theta \otimes\operatorname{Id} )\}[a,A]| \, d\sigma_{n-1}(\theta).   
\end{align*}
By Lemma \ref{lemma:integral-id} and Hölder's inequality
\begin{align*}
    R^3 |D^3v(x)[a,A]| \leq n \sqrt{6(n+4)} \sup_{B_{2R}} |v|.
\end{align*}
For $x,y\in B_R$ we have
\begin{align*}
   D^2v(x) - D^2v(y) = \int_0^{|x-y|} D^3 v(y + t \tfrac{x-y}{|x-y|})[\tfrac{x-y}{|x-y|},\,\cdot\,,\,\cdot\,]\,dt.
\end{align*}
Using the above estimate with $a = \tfrac{x-y}{|x-y|}$ and $A= (D^2v(x) - D^2v(y))/|D^2v(x) - D^2v(y)| $ gives
the Lipschitz bound
\begin{align*}
    |D^2v(x)-D^2v(y)| \leq \underbrace{R^{-3} n \sqrt{6(n+4)} \sup_{B_{2R}} |v|}_{=L}\, |x-y|.
\end{align*}
\textbf{Hölder bound for second derivative:}
We have
\begin{align*}
    |D^2v(x)-D^2v(y)| \leq \min\{2 [v]_{2,0;B_R}, L |x-y| \} \leq (2 [v]_{2,0;B_R})^{1-\alpha} L^{\alpha} |x-y|^{\alpha}. 
    \end{align*}
This gives the desired bound.
\end{proof}

\begin{proposition}\label{prop:local-Schauder}
    We have
    \begin{align*}
        \begin{pmatrix}
            [u]_{0;B_R}
            \\
            R[u]_{1,0;B_R}
            \\
            R^2 [u]_{2,0;B_R}
            \\
            R^{2+\alpha} [u]_{2,\alpha;B_{R}}
        \end{pmatrix}
        \leq 
        \begin{pmatrix}
            1 & 0 & 0
            \\
            \Cr{harmonic-1} & \Cr{Schauder-1-0}+ \tfrac{4}{n(n-2)} \Cr{harmonic-1} & 0
            \\
             \Cr{harmonic-2} & \Cr{Schauder-2-0-1}+ \tfrac{4}{n(n-2)} \Cr{harmonic-2} & \Cr{Schauder-2-0-2}
             \\
             \Cr{harmonic-3} & \Cr{Schauder-2-alpha-1}+ \tfrac{4}{n(n-2)} \Cr{harmonic-3} & \Cr{Schauder-2-alpha-2}
        \end{pmatrix}
        \begin{pmatrix}
            \|{u}_{C^0(B_{2R})}
            \\
            R^2 [\Delta u]_{0;B_{2R}}
            \\
            R^{2+\alpha} [\Delta u]_{0,\alpha;B_{2R}}
        \end{pmatrix},
    \end{align*}
    and
    \begin{align*}
        \|{u}'_{C^{2,\alpha}(B_R)} \leq \Cl{loc-Schauder-1} R^2\|{\Delta u}''_{C^{0,\alpha}(B_{2R})} + \Cl{loc-Schauder-2} \|{u}_{C^0(B_{2R})},
    \end{align*}
    where
    \begin{align*}
        \Cr{loc-Schauder-1} = \max\{\Cr{Schauder-2-0-2} + \Cr{Schauder-2-alpha-2} , \Cr{Schauder-1-0} + \Cr{Schauder-2-0-1} + \Cr{Schauder-2-alpha-1} +\frac{4}{n(n-2)} (\Cr{harmonic-1} + \Cr{harmonic-2} + \Cr{harmonic-3}) \} , \quad \Cr{loc-Schauder-2} = 1 + \Cr{harmonic-1} + \Cr{harmonic-2} + \Cr{harmonic-3}.
    \end{align*}
\end{proposition}
\begin{proof}
    Set $f = \Delta u$ and denote by $w$ the Newton potential of $f$ in $B_2$. Set $v = u-w$, so that $v$ is harmonic. 
    We first note that $[v]_{0;B_2}  \leq \|{u}_{C^0(B_2)} +  \frac{4}{n(n-2)} [f]_{0;B_2}$, indeed by the maximum principle we have 
    $$\sup_{B_2} |v| \leq \sup_{\partial B_2}|v| \leq \sup_{\partial B_2}|u| + \sup_{\partial B_2}|w| \leq \|{u}_{C^0(B_2)} +  \sup_{\partial B_2}|w| $$ 
   and for $x\in \partial B_2$ by Lemma \ref{lemma:integral-Newton-kernel}
    \begin{align*}
       |w(x)| \leq [f]_{0;B_2} \int_{B_2} |\Gamma(x-y)| dy \leq \frac{4}{n(n-2)} [f]_{0;B_2} .
    \end{align*}
    Thus we have 
    \begin{align*}
        [u]_{1,0;B_1} &\leq [w]_{1,0;B_1} + [v]_{1,0;B_1} \leq \Cr{Schauder-1-0} [f]_{0;B_2} + \Cr{harmonic-1} [v]_{0;B_2} 
        \\
        &\leq (\Cr{Schauder-1-0} + \tfrac{4}{n(n-2)} \Cr{harmonic-1}) [f]_{0;B_2} + \Cr{harmonic-1} [u]_{0;B_2},
        \\
        [u]_{2,0;B_1} &\leq [w]_{2,0;B_1} + [v]_{2,0;B_1} \leq \Cr{Schauder-2-0-1} [f]_{0;B_2} + \Cr{Schauder-2-0-2} [f]_{0, \alpha;B_2} + \Cr{harmonic-2} [v]_{0;B_2}
        \\
        & \leq (\Cr{Schauder-2-0-1} + \tfrac{4}{n(n-2)} \Cr{harmonic-2}) [f]_{0;B_2} + \Cr{Schauder-2-0-2} [f]_{0, \alpha;B_2} + \Cr{harmonic-2} [u]_{0;B_2},
        \\
        [u]_{2,\alpha;B_1} &\leq [w]_{2,\alpha;B_1} + [v]_{2,\alpha;B_1} \leq \Cr{Schauder-2-alpha-1} [f]_{0;B_2} + \Cr{Schauder-2-alpha-2} [f]_{0, \alpha;B_2} + \Cr{harmonic-3} [v]_{0;B_2}
        \\
        & \leq (\Cr{Schauder-2-alpha-1} +  \tfrac{4}{n(n-2)} \Cr{harmonic-3}) [f]_{0;B_2} + \Cr{Schauder-2-alpha-2} [f]_{0, \alpha;B_2} + \Cr{harmonic-3} [u]_{0;B_2},
    \end{align*}
\end{proof}
\begin{remark}
    Numerically for $n=4$ and $\alpha=\frac{1}{2}$ we have
    \begin{align*}
        \begin{pmatrix}
            [u]_{0;B_R}
            \\
            R[u]_{1,0;B_R}
            \\
            R^2 [u]_{2,0;B_R}
            \\
            R^{2+\alpha} [u]_{2,\alpha;B_{R}}
        \end{pmatrix}
        \leq 
        \begin{pmatrix}
            1 & 0 & 0
            \\
            1.70 & 2.13 & 0
            \\
             6.93 & 4.58 & 4
             \\
             19.60 & 14.52 & 25.97
        \end{pmatrix}
        \begin{pmatrix}
            \|{u}_{C^0(B_{2R})}
            \\
            R^2 [\Delta u]_{0;B_{2R}}
            \\
            R^{2+\alpha} [\Delta u]_{0,\alpha;B_{2R}}
        \end{pmatrix},
    \end{align*}
    and 
    \begin{align*}
        \Cr{loc-Schauder-1}\leq29.97,\qquad
        \Cr{loc-Schauder-2}\leq29.23.    \end{align*}
\end{remark}

\section{Analysis on manifolds}

\subsection{Sobolev embedding theorem}

The following is \cite[Theorem~5]{Ilias1983}, after rescaling and using the geometric bounds in \cite[Theorems~1(i) and~2(i), p.~157]{Ilias1983}, except that the coefficient of $\|{f}_{L^q}$ is improved from $2$ to $2^{1/q}$.

For $n\geq2$ and $1<q<n$, define
\[
    \omega_{n-1}=\vol_{n-1}(S^{n-1})=\frac{2\pi^{n/2}}{\Gamma(n/2)},
    \qquad C^*=n^{1-1/n}\omega_{n-1}^{1/n},
\]
\[
    K(n,q)=n^{-1/q}
    \left(\frac{q-1}{n-q}\right)^{1-1/q}
    \left[
        \frac{\Gamma(n+1)}
        {\omega_{n-1}\Gamma(n/q)\Gamma(n+1-n/q)}
    \right]^{1/n}.
\]
For $r>0$, define
\[
    h_r=\left[\int_0^{r/2}\cosh^{n-1}t\,dt\right]^{-1},
    \qquad
    C_{1,r}=\left[
        \int_0^r\left(\frac{\cosh t}{h_r}+\frac{\sinh t}{n}\right)^{n-1}dt
    \right]^{-1/n}.
\]

\begin{theorem}
\label{theorem:Ilias-Thm-5-update}
Let $(M,g)$ be a closed connected Riemannian manifold of dimension $n\geq2$ and volume $1$.
Suppose that $\kappa,D>0$ and
\[
    \Ric_g\geq-(n-1)\kappa g,
    \qquad \diam(M,g)<D.
\]
For $1<q<n$, put $p=nq/(n-q)$.
Then every $f\in L^q_1(M)$ satisfies
\[
    \|{f}_{L^p(M)}
    \leq
    A_q
    \|{\nabla f}_{L^q(M)}
    +2^{1/q}\|{f}_{L^q(M)}
    \quad
    \text{for}
    \quad
    A_q = \frac{C^*K(n,q)}{\sqrt\kappa\,C_{1,\sqrt\kappa D}}.
\]
\end{theorem}

\begin{proof}
As in \cite[p.~157]{Ilias1983}, define the constants
\[
    h(M,g)=\inf_{0<\vol_g(\Omega)\leq\vol_g(M)/2}
    \frac{\vol_{n-1,g}(\partial\Omega)}{\vol_g(\Omega)},
    \qquad
    C_1(M,g)=\inf_{0<\vol_g(\Omega)\leq\vol_g(M)/2}
    \frac{\vol_{n-1,g}(\partial\Omega)}{\vol_g(\Omega)^{(n-1)/n}},
\]
where $\Omega$ ranges over smooth domains in $M$.
Set $\widetilde g=\kappa g$.
Then
\[
    \Ric_{\widetilde g}\geq-(n-1)\widetilde g,
    \qquad \vol_{\widetilde g}(M)=\kappa^{n/2},
    \qquad \diam(M,\widetilde g)<\sqrt\kappa D.
\]
By \cite[Theorem~1(i)]{Ilias1983},
\[
    h(M,\widetilde g)
    \geq\left[\int_0^{\diam(M,\widetilde g)/2}\cosh^{n-1}t\,dt\right]^{-1}
    \geq\left[\int_0^{\sqrt\kappa D/2}\cosh^{n-1}t\,dt\right]^{-1}
    =h_{\sqrt\kappa D}.
\]
We have:
\begin{align*}
    C_1(M,\widetilde g)
    &\geq\vol_{\widetilde g}(M)^{1/n}
    \left[\int_0^{\diam(M,\widetilde g)}
    \left(\frac{\cosh t}{h(M,\widetilde g)}+\frac{\sinh t}{n}\right)^{n-1}dt\right]^{-1/n}
    \\
    &=\sqrt\kappa
    \left[\int_0^{\diam(M,\widetilde g)}
    \left(\frac{\cosh t}{h(M,\widetilde g)}+\frac{\sinh t}{n}\right)^{n-1}dt\right]^{-1/n}
    \\
    &\geq\sqrt\kappa
    \left[\int_0^{\sqrt\kappa D}
    \left(\frac{\cosh t}{h_{\sqrt\kappa D}}+\frac{\sinh t}{n}\right)^{n-1}dt\right]^{-1/n}
    =\sqrt\kappa\,C_{1,\sqrt\kappa D},
\end{align*}
where in the first step we applied \cite[Theorem~2(i)]{Ilias1983}, in the second step we used $\vol_{\widetilde g}(M)^{1/n}=\sqrt\kappa$, and in the third step we used the preceding bound on $h(M,\widetilde g)$ and $\diam(M,\widetilde g)<\sqrt\kappa D$.
The constants $h$ and $C_1$ scale as $h(M,\widetilde g)=\kappa^{-1/2}h(M,g)$ and $C_1(M,\widetilde g)=C_1(M,g)$.
We therefore have $h(M,g)\geq\sqrt\kappa\,h_{\sqrt\kappa D}$ and $C_1(M,g)\geq\sqrt\kappa\,C_{1,\sqrt\kappa D}$.

Suppose first that $f\geq0$ is a smooth Morse function and let $a\geq0$ so that
\[
\mathrm{vol}_g(\{f>a\})
=
\mathrm{vol}_g(\{f<a\})
=
\frac{1}{2}.
\]
As in \cite[proof of Theorem~5, p.~161, third displayed formula]{Ilias1983}, one derives
\[
    \|{f}_{L^p(M)}
    \leq\frac{C^*K(n,q)}{C_1(M,g)}\|{\nabla f}_{L^q(M)}+a.
\]
Since $f^q\geq a^q$ on $\{f\geq a\}$, a set of volume $1/2$,
\[
    \frac{a^q}{2}
    = a^q\vol_g\{f\geq a\}
    \leq\int_{\{f\geq a\}}f^q\,d\vol_g
    \leq\|{f}_{L^q(M)}^q.
\]
Taking $q$-th roots gives $a\leq2^{1/q}\|{f}_{L^q(M)}$.
Substituting this estimate and the lower bound for $C_1(M,g)$ gives the asserted inequality.
As in \cite[\S~2.1, p.~158]{Ilias1983}, approximation and passage to $|f|$ extend the estimate to every $f\in L^q_1(M)$.
\end{proof}

\begin{corollary}
\label[corollary]{cor:L22-Sobolev-embedding-update}
Let $(M,g)$ be a closed connected Riemannian manifold of dimension $n\geq3$ and volume $1$.
Suppose that $\kappa,D>0$, $\Ric_g\geq-(n-1)\kappa g$, and $\diam(M,g)<D$.
Let $1\leq p<\infty$ satisfy $1/2-2/n<1/p<1/2-1/n$.
Then every $f\in L^2_2(M)$ satisfies
\[
    \|{f}_{L^p(M)}
    \leq\Cr{L22-Sobolev-update}\|{f}_{L^2_2(M)},
\]
where $\Cr{L22-Sobolev-update}$ is defined in \cref{eqn:L22-Sobolev-update}.
\end{corollary}
\begin{proof}
Set $t=np/(n+p)$ and $\theta=(n-2)/2-n/p$.
Then $2<t<n$, $0<\theta<1$, and $1/t=(1-\theta)/2+\theta(n-2)/(2n)$.
Applying Theorem~\ref{theorem:Ilias-Thm-5-update} to $|T|$, and using Kato's inequality, gives
\begin{equation}
\label{eqn:L2-first-order-Sobolev-update}
    \|{T}_{L^{2n/(n-2)}(M)}
    \leq A_2\|{\nabla T}_{L^2(M)}+\sqrt2\,\|{T}_{L^2(M)}
    \qquad
    \text{ for }
    \qquad
    T=f,\nabla f.
\end{equation}
For every $\gamma>0$, we have
\begin{align*}
    \|{T}_{L^t(M)}
    &\leq\|{T}_{L^2(M)}^{1-\theta}
    \|{T}_{L^{2n/(n-2)}(M)}^\theta
    \\
    &\leq\|{T}_{L^2(M)}^{1-\theta}
    \left(A_2\|{\nabla T}_{L^2(M)}+\sqrt2\,\|{T}_{L^2(M)}\right)^\theta
    \\
    &\leq\gamma^{\theta-1}\left[
    \bigl((1-\theta)\gamma+\sqrt2\,\theta\bigr)\|{T}_{L^2(M)}
    +A_2\theta\|{\nabla T}_{L^2(M)}\right],
\end{align*}
where in the first step we used H\"older's inequality, in the second step we used \cref{eqn:L2-first-order-Sobolev-update}, and in the third step we used Young's inequality $a^{1-\theta}b^\theta\leq(1-\theta)a+\theta b$, inserting the factor $(\gamma^\theta)^{1-\theta}(\gamma^{\theta-1})^\theta=1$.
Substituting this bound for $T=f$ and $T=\nabla f$ into the Sobolev embedding $\|{f}_{L^p(M)}\leq A_t\|{\nabla f}_{L^t(M)}+2^{1/t}\|{f}_{L^t(M)}$ from Theorem~\ref{theorem:Ilias-Thm-5-update} gives
\begin{align*}
    \|{f}_{L^p(M)}
    &\leq\gamma^{\theta-1}\Bigl[
    2^{1/t}\bigl((1-\theta)\gamma+\sqrt2\,\theta\bigr)
    \|{f}_{L^2(M)}
    \\
    &\qquad+\bigl[
    A_t\bigl((1-\theta)\gamma+\sqrt2\,\theta\bigr)+2^{1/t}A_2\theta
    \bigr]\|{\nabla f}_{L^2(M)}
    \\
    &\qquad+A_tA_2\theta\|{\nabla^2f}_{L^2(M)}\Bigr].
\end{align*}
We now choose
\begin{equation}
\label{eqn:L22-Sobolev-update}
\begin{gathered}
    \gamma=\max\left\{\sqrt2,\frac{\theta(A_2-\sqrt2)}{1-\theta}\right\},
    \\[1mm]
    \Cl{L22-Sobolev-update}
    =\gamma^{\theta-1}\max\left\{
    \begin{gathered}
        2^{1/t}\bigl((1-\theta)\gamma+\sqrt2\,\theta\bigr),\\
        A_t\bigl((1-\theta)\gamma+\sqrt2\,\theta\bigr)+2^{1/t}A_2\theta
    \end{gathered}
    \right\}.
\end{gathered}
\end{equation}
Since $(1-\theta)\gamma+\sqrt2\,\theta\geq A_2\theta$, the coefficient of $\|{\nabla^2f}_{L^2(M)}$ is no larger than that of $\|{\nabla f}_{L^2(M)}$, and the asserted estimate follows.
\end{proof}

\subsection{$L^2$-elliptic estimate}

\begin{proposition}
\label{prop:L2-ell-estimate}
    Suppose that $\kappa\geq0$ and 
        \begin{align}
    \operatorname{Ric}_g\geq-(n-1)\kappa g.
    \label{Ricci-bound}
    \end{align}
    Then
    \begin{align*}
        \|{u}_{L^2_2(M)} \leq \Cl{L2-ell-estimate} (\|{\Delta u}_{L^2(M)} + \|{u}_{L^2(M)}),
    \end{align*}
    where
    \begin{align*}
        \Cr{L2-ell-estimate} = \left(\frac{3}{2} (3+(n-1)\kappa) \right)^{\frac{1}{2}}.
    \end{align*}
\end{proposition}
\begin{proof}

In the following write $\nabla$ for $\nabla^g$ and $\Delta$ for $\Delta_g$.

By integration by parts and Cauchy--Schwarz we have
\begin{align}
    \|{\nabla u}^2_{L^2(M)} \leq \|{\Delta u}_{L^2(M)} \|{u}_{L^2(M)} \leq \frac{1}{2} \|{\Delta u}^2_{L^2(M)} +  \frac{1}{2}\|{u}^2_{L^2(M)}.  
    \label{one-der-L2}
\end{align}

The Bochner formula is
\begin{align*}
    \frac{1}{2} \Delta |\nabla u|^2 = g(\nabla \Delta u, \nabla u) + |\nabla^2 u|_g^2 + \operatorname{Ric}(\nabla u, \nabla u).
\end{align*}
Integrating this identity over $M$, integrating by parts, using \eqref{Ricci-bound} and \eqref{one-der-L2} gives
\begin{align*}
 \|{\nabla^2 u}^2_{L^2(M)} \leq  \|{\Delta u}^2_{L^2(M)} + (n-1) \kappa \|{\nabla u}^2_{L^2(M)}  
 \leq \left(1+\frac{n-1}{2} \kappa \right) \|{\Delta u}^2_{L^2(M)} + \frac{n-1}{2} \kappa \|{u}^2_{L^2(M)}.
\end{align*}
This gives
\begin{align*}
    \|{u}^2_{L^2_2(M)} \leq 3 (\|{\nabla^2 u}^2_{L^2(M)} + \|{\nabla u}^2_{L^2(M)} + \|{u}^2_{L^2(M)})
    \leq 3 \frac{(3+(n-1)\kappa)}{2} (\|{\Delta u}^2_{L^2(M)} + \|{u}^2_{L^2(M)}).
\end{align*}
\end{proof}

\section{Details of the numerical construction}
\label{section:numerical-details}

This appendix gives the details of the method of \cref{subsection:computing-metrics}, and it shows that the method reaches smaller residuals.
As in \cref{subsection:computing-metrics}, the residual numbers are non-rigorous Monte-Carlo estimates.

\subsection{Invariant Hermitian forms}
\label{subsection:invariant-forms}

As in \cref{subsection:computing-metrics}, the coefficient matrix $(N_{AB})$ of an invariant form $N(X)=\sum_{A,B}N_{AB}X^A\bar X^B$ of bidegree $(k,k)$ has $64$ blocks, one for each residue $r=A\bmod 4$, and here $|r|\equiv k\pmod 4$.
In the block of $r$ we write $X^A=X^ru^C$ with $u_a=X_a^4$.
On $M$ the variables $u$ satisfy $u_0+u_1+u_2+u_3=0$, so the block acts on the polynomials of degree $m_r=(k-|r|)/4$ on this plane, a space of dimension $\binom{m_r+2}{2}$.
The permutations in $S_4$ permute the blocks, so one block of each orbit carries the free parameters, and the stabiliser of $r$ constrains them.
The number of free parameters is therefore
\[
    n_k=\sum_{[r]}\dim\bigl(\Sym^2\Sym^{m_r}W\bigr)^{\Stab(r)},
\]
where $[r]$ runs over the $S_4$-orbits of residues and $W$ is the three-dimensional standard representation of $S_4$ on the plane.
\Cref{tab:parameter-counts} gives $n_k$ for the degrees that we use.

\begin{table}[ht]
\centering
\begin{tabular}{@{}lcccccccc@{}}
\toprule
$k$ & $12$ & $16$ & $20$ & $24$ & $30$ & $32$ & $34$ & $36$\\
\midrule
$n_k$ & $56$ & $135$ & $285$ & $542$ & $1{,}225$ & $1{,}559$ & $1{,}958$ & $2{,}432$\\
\bottomrule
\end{tabular}
\caption{The number $n_k$ of free parameters of an invariant form of bidegree $(k,k)$ on $M$.}
\label{tab:parameter-counts}
\end{table}

\subsection{The weight of the levels}
\label{subsection:level-weight}

This subsection derives the statement on the weight of the levels in \cref{subsection:computing-metrics}.
Let $h$ be a Hermitian metric on $\mathcal O(1)$ over $M$ with positive curvature form $\omega_h$, and let $s_1,\dots,s_m$ be a basis of $H^0(M,\mathcal O(L))$ that is orthonormal for the $L^2$ inner product of $h^L$ and $\omega_h^2/2$.
By the Tian--Yau--Zelditch expansion \cite{Tian1990,Zelditch1998}, the Bergman density $\rho_L=\sum_i|s_i|^2_{h^L}$ satisfies $\rho_L=c\,(1+O(L^{-1}))$ with all derivatives, where $c$ is a constant of the order of $L^2$.

Let $|s|^2_h=|s|^2e^{-\Phi_{\mathrm{CY}}/k_0}$ for a potential $\Phi_{\mathrm{CY}}$ of the Calabi--Yau metric, and $N=\sum_i|s_i|^2$.
Since $\Phi_{\mathrm{CY}}$ is invariant under $G$ and complex conjugation, $N$ is a form of the same kind as $N_L$, and $\rho_L=Ne^{-(L/k_0)\Phi_{\mathrm{CY}}}$.
So, for every base,
\[
    e^{(L/k_0)(\Phi_{\mathrm{CY}}-\log p_0)}=\frac{N\,p_0^{-L/k_0}}{\rho_L},
\]
where $Np_0^{-L/k_0}$ lies in $V_L$ and $\rho_L$ is almost constant.
The base is close to the Calabi--Yau metric, so for a suitable constant in $\Phi_{\mathrm{CY}}$ the correction is small, and the linearisation of the left side gives
\[
    \Phi_{\mathrm{CY}}-\log p_0=\frac{k_0}{Lc}\,N\,p_0^{-L/k_0}-\frac{k_0}{L}+O(L^{-2})+O\bigl((\Phi_{\mathrm{CY}}-\log p_0)^2\bigr).
\]
The first term lies in $V_L$ and is almost the constant $k_0/L$, and the constant $-k_0/L$ does not change the metric.
With the weight $|X|^{-2L}$, the identity reads $e^{(L/k_0)(\Phi_{\mathrm{CY}}-\log p_0)}=N|X|^{-2L}(|X|^2p_0^{-1/k_0})^L/\rho_L$, and the factor $(|X|^2p_0^{-1/k_0})^L$ varies over $M$ by an amount that grows exponentially with $L$, so forms of degree $L$ approximate it poorly.

\subsection{The orthonormal bases of the levels}
\label{subsection:level-bases}

We write each function $N_Lp_0^{-L/k_0}$ of the level $V_L$ as $W_Lq$ with the weight $W_L=(|X|^2p_0^{-1/k_0})^L$ and $q=N_L/|X|^{2L}$, and we call $Q_L$ the space of all these functions $q$ on $M$.
The recurrence below builds the functions $q$, and in it the weight $W_L$ enters only through the inner product.
Because $N/|X|^{2d}=N|X|^2/|X|^{2(d+1)}$, these spaces are nested, $Q_0\subset Q_1\subset\dots\subset Q_L$, and $Q_d$ has the dimension $n_d$ of \cref{subsection:invariant-forms}.
The basis of the level is a sequence of functions $q_t\in Q_L$ that are orthonormal for the inner product
\[
    \langle a,b\rangle_L=\sum_iw_i\,W_L(x_i)^2\,a(x_i)\,b(x_i)
\]
on a sample of build points $x_i$, with the weights $w_i$ of \cref{subsection:computing-metrics}.
So the functions $W_Lq_t$ are orthonormal on the sample.

The multipliers of the recurrence are the functions
\[
    y_{AB}=\Re\sum_{\pi\in S_4}\frac{X^{\pi(A)}\bar X^{\pi(B)}}{|X|^{2d}},
    \qquad
    |A|=|B|=d\leq12,
    \quad
    A\equiv B\pmod 4,
\]
which are $G$-invariant and lie in $Q_d$.
The basis starts with the constant function
\[
    q_0=\Bigl(\sum_iw_i\,W_L(x_i)^2\Bigr)^{-1/2},
\]
which spans $Q_0$ and has norm $1$, and we write $d_t$ for the degree of the step that builds $q_t$, so $d_0=0$.

The basis of a stage comes from two pivoted Gram--Schmidt processes, which differ in their inputs, their inner products and their rules to stop:
\begin{enumerate}
\item
\emph{The basis of one level $L$.}
At each degree $d=1,\dots,L$ in turn, the inputs are the candidates $y_{AB}\,q_j$ with $|A|+d_j=d$; for $|A|=d$ the earlier function is $q_0$, so the multipliers themselves are candidates.
Each step takes the candidate whose component orthogonal to all functions of this level so far is largest relative to its norm, in $\langle\cdot,\cdot\rangle_L$ on the build sample of the level.
If this candidate is $y_{AB}\,q_j$ and the new function has the index $t$, we call $y_t=y_{AB}$ its multiplier and $q_j$ its parent, and we write $p(t)=j$, so $p(t)<t$.
Then
\[
    q_t=\frac{1}{h_{tt}}\Bigl(y_t\,q_{p(t)}-\sum_{s<t}h_{st}\,q_s\Bigr),
    \qquad
    h_{st}=\langle y_t\,q_{p(t)},q_s\rangle_L,
\]
where $h_{tt}>0$ gives $q_t$ the norm $1$, and the step removes the component along $q_t$ from the other candidates.
The degree $d$ is complete when the level has $n_d=\dim Q_d$ functions in total; with multipliers up to degree $12$ the candidates reach this number at every degree that we use.
The stored indices $p(t)$, multipliers $y_t$ and coefficients $h_{st}$ then give the functions $q_t$ at any point of $M$.
\item
\emph{The basis of a stage.}
The inputs are the functions $W_Lq_t$ of all levels $L\in\mathcal L_1$; in the second stage of \cref{subsection:two-stage-results} they are the functions $|X|^{2L}e^{-(L/k_0)\Phi_1}q_t$ with $L\in\mathcal L_2$ and the same $q_t$.
The process takes the levels in turn from the highest degree down, and inside a level it uses the same pivot rule, now for the inner product $\sum_iw_i\,a(x_i)\,b(x_i)$ on $60{,}000$ build points, against all functions that it has accepted so far.
A function enters only if its new part $\nu$, the norm of its component orthogonal to these functions relative to its own norm, is larger than the cutoff $\tau$.
Each accepted function gets one column of the least-squares matrix, and the orthogonal projection of the constant function $1$ onto their span gets one more.
This function carries the main part of the correction, and without its own column the least-squares solve in double precision cannot resolve it.
\end{enumerate}

In double precision the error of the recurrence grows by a factor of about $3.2$ for each degree, and in double-double arithmetic by about $2.6$.
So the bases must be built and evaluated in double-double arithmetic: one level of degree $36$ gives $\operatorname{mean}|f|=\num{3.5e-16}$ with the basis functions in double precision, and $\num{1.1e-17}$ in double-double arithmetic.
To limit the memory, process 1 runs first in double precision, and its run in double-double arithmetic then takes only the multipliers that the first run picks.

Three details of the double-double arithmetic matter:
\begin{itemize}
\item
All products use the exact product of Dekker \cite{Dekker1971}.
A faster product with a relative error of about $2^{-79}$ puts a floor near $10^{-24}$ on $f$.
\item
The sampled points are double-precision points about $10^{-16}$ away from $M$.
At such points a basis function of degree $36$ changes by up to $7\%$ of its largest value.
So one Newton step in double-double arithmetic moves each point onto $M$, to a relative error of about $10^{-31}$.
\item
The Gram--Schmidt processes split the matrix products into four slices with an error-free transformation \cite{OzakiOgitaOishiRump2012}.
\end{itemize}

\subsection{Smaller residuals}
\label{subsection:two-stage-results}

The same construction gives a metric for each required accuracy.
\Cref{tab:other-degrees} uses the eight consecutive levels $D,D-1,\dots,D-7$ below the top degree $D$ in the last stage, and a cutoff near the error of the level functions.
With one stage the best base degree is about $2D/3$, and four more degrees gain a factor of $\num{1.8e4}$ from $D=16$ to $D=20$, but only $5$ to $7$ from $D=32$ to $D=36$, because the cutoff must grow with the error of the level functions.
The metric \eqref{eqn:certified-potential} of this paper has twelve levels, $24$ to $13$, in place of eight, and $\operatorname{mean}|f|=\num{1.73e-22}$.

\begin{table}[ht]
\centering
\begin{tabular}{@{}lllll@{}}
\toprule
$D$ & one stage & base degree & two stages & base degree; first stage\\
\midrule
$16$ & \num{1.4e-14} & $10$ & -- & --\\
$20$ & \num{7.3e-19} & $14$ & -- & --\\
$24$ & \num{2.9e-22} & $16$ & \num{1.1e-22} & $8$; levels $16$, $12$\\
$28$ & \num{1.3e-24} & $18$ & \num{1.3e-25} & $10$; levels $20$, $16$, $12$\\
$32$ & \num{2.4e-26} & $28$ & \num{1.5e-26} & $12$; levels $24$, $20$, $16$\\
$36$ & \num{3.5e-27} & $24$ & \num{1.3e-27} & $12$; levels $24$ to $17$\\
$36$, $\tau=10^{-20}$ & -- & -- & \num{2.7e-28} & $12$; levels $24$ to $17$\\
\bottomrule
\end{tabular}
\caption{$\operatorname{mean}|f|$ on $400{,}000$ hold-out points for the top degree $D$. The last row is the two-stage metric.}
\label{tab:other-degrees}
\end{table}

A second stage corrects $\Phi_1$ of \eqref{eqn:one-stage-potential} with further terms,
\begin{align}
\label{eqn:second-stage}
    \Phi_2=\Phi_1+\sum_{L\in\mathcal L_2}N'_L\,e^{-(L/k_0)\Phi_1},
\end{align}
with a second set $\mathcal L_2$ of degrees and forms $N'_L$ of the same kind as the $N_L$, fitted in one linear least-squares problem with $\Phi_1$ fixed.
With the weights $p_0^{-L/k_0}$ in place of $e^{-(L/k_0)\Phi_1}$, a second stage adds no function to the span of one stage with the levels of $\mathcal L_2$, because the levels are nested and the problem is linear.
At the same cutoff, two stages are better than one by a factor of $1.6$ to $10.6$, and the best base degree is then about $D/3$.
Our most accurate metric of this kind has $k_0=12$, $\mathcal L_1=\{17,18,\dots,24\}$, $\mathcal L_2=\{29,30,\dots,36\}$ and $\tau=10^{-20}$ in both stages.
It gives $\operatorname{mean}|f|=\num{2.74e-28}$ on $8\cdot10^6$ hold-out points, and $(2.71\pm0.16)\cdot10^{-28}$ on $256$ random points in mpmath \cite{mpmath} at $100$ significant digits.

At this level the rounding noise of $f$ in double-double arithmetic is no longer small against $f$: it is about $\num{1.5e-29}$ on average, and up to $\num{8e-29}$ near the curves $M\cap\{X_a=0\}$.
The fit points and the hold-out points write each point in the same chart, so noise that depends on the chart does not show in the hold-out values.
So we also compute $f$ in the chart of the second largest coordinate and at $12$ test points in mpmath, and we accept a fit only if both change $f$ by at most about $10\%$ of $\operatorname{mean}|f|$ and the mpmath sample agrees with the double-double value.
With the second-stage cutoff $10^{-21}$ in place of $10^{-20}$ the double-double value falls to $\num{1.1e-28}$, but the second chart changes $f$ by $17\%$ of $\operatorname{mean}|f|$ and the test points by $35\%$, so this gain is rounding noise.

\section{Rigorous error bounds using Taylor models}
\label{section:rigorous-error-bounds}

This appendix states the steps of the proof of Proposition~\ref{proposition:certified-residual} and of the computer-assisted bounds of \cref{section:bounding-quantities} as definitions and lemmas.
The arithmetic of Taylor models in double-double and triple-double numbers on GPUs, with the bounds of its rounding errors, is the subject of the technical report \cite{technical-report-reference}.
The code and further details of the computation are in our repository \cite{code-repository}.
A lemma marked \emph{a computation} is a statement about stored numbers that a program checks in exact or ball arithmetic.
Throughout, $u=2^{-53}$ is the unit roundoff of double precision, and $\mathbb D^k_\delta=\{x\in\CC^k:|x_1|,\dots,|x_k|\leq\delta\}$ is the closed polydisc.
We use the affine coordinates $y=(X_1,X_2,X_3)$ of the chart $X_0=1$, in which $M$ is the zero set of $W(y)=1+X_1^4+X_2^4+X_3^4$, the coordinates $t=(X_2,X_3)$, and the region $\mathcal F$ of \cref{subsection:certified-errors-bounds}.

\subsection{Taylor models and rounding}
\label{subsection:taylor-models}

A \emph{Taylor model} $T$ of order $d$ on $\mathbb D^k_\delta$ is a polynomial $P_T(x)=\sum_{|\alpha|\leq d}c_\alpha x^\alpha$, whose coefficients are double-double or triple-double numbers taken as exact, together with a number $\epsilon_T\geq0$.
It \emph{encloses} a function $\varphi$ if $\varphi$ is continuous on $\mathbb D^k_\delta$, holomorphic in its interior, and $|\varphi-P_T|\leq\epsilon_T$ on $\mathbb D^k_\delta$.
Then
\begin{align}
\label{eqn:taylor-model-bound}
    \sup_{\mathbb D^k_\delta}|\varphi|
    \leq
    |c_0|+\sum_{\alpha\neq0}|c_\alpha|\,\delta^{|\alpha|}+\epsilon_T.
\end{align}
We use these operations \cite{technical-report-reference}; in each of them, the remainder bound of the result also contains the rounding errors of its coefficients, each multiplied by $\delta^{|\alpha|}$.
\begin{itemize}
\item
\emph{Sums and products.}
The polynomial of a product is the product of the polynomials, truncated at the order $d$.
Its remainder bound is $\epsilon_{\mathrm{tail}}+b_T\epsilon_{T'}+\epsilon_Tb_{T'}+\epsilon_T\epsilon_{T'}$, where $b_T=\sum_\alpha|c_\alpha|\delta^{|\alpha|}$ and $\epsilon_{\mathrm{tail}}$ bounds the dropped terms of degree above $d$.
\item
\emph{Series.}
If $T$ encloses $\varphi$ and $b_T+\epsilon_T<1$, then truncated binomial and logarithm series in $T$, with the Lagrange bounds for their rest, enclose $(1+\varphi)^a$ for rational $a$ and $\log(1+\varphi)$.
\item
\emph{Cauchy estimates.}
If $T$ encloses $\varphi$ on $\mathbb D^k_{\delta'}$ and $\delta<\delta'$, then $(\partial_iP_T,\epsilon_T/(\delta'-\delta))$ encloses $\partial_i\varphi$ on $\mathbb D^k_\delta$, and $(\partial_i\partial_jP_T,\epsilon_T/(\delta'-\delta)^2)$ encloses $\partial_i\partial_j\varphi$ for $i\neq j$.
\item
\emph{Polarisation.}
A real-analytic function $\varphi(w)$ of $w\in\CC^2$ near $0$ is the restriction to the slice $w'=\bar w$ of a holomorphic function $\tilde\varphi(w,w')$ of four complex variables.
For example, $|X_1|^2$ is the restriction of $X_1(w)\overline{X_1(\overline{w'})}$.
We model $\tilde\varphi$ on $\mathbb D^4_\delta$, and then $\partial_{w_i}\partial_{\bar w_j}\varphi(w)=\partial_{w_i}\partial_{w'_j}\tilde\varphi(w,\bar w)$.
\end{itemize}
Each operation bounds its own rounding errors from its computed results: a rounded sum $r$ has an error of at most $u|r|$, and a rounded product, quotient or square root an error of at most $u|r|+2^{-1074}$ \cite{IEEE754}, and the bounds themselves are formed with the directed rounding of the GPU.
So no bound depends on a count of operations, on the order of the parts of a number, or on the absence of underflow.
The steps with the largest cancellations run in triple-double arithmetic: the solve for the points of $M$, the models of the monomials, the change of centre of the forms and the matrix products with the forms.
The numbers that the GPUs take as input come from ball arithmetic \cite{Johansson2017} or from exact rationals, with bounds of their rounding.
A step whose condition fails, for example a series with $b_T+\epsilon_T\geq1$, gives an infinite or undefined bound, and the later steps carry it to the bound of the box.
So a box with a finite bound satisfies the conditions of all lemmas of \cref{subsection:box-model}, except the condition $s_p>0$ of Lemmas~\ref{lemma:level-sum} and~\ref{lemma:metric-eta}, which follows from Lemma~\ref{lemma:sign}.

\subsection{The potential as exact data}
\label{subsection:exact-potential}

The potential \eqref{eqn:certified-potential} is given by stored numbers: the blocks of the matrix $H$ of $p_0$, each entry the exact sum of two double-precision numbers, and the matrices of the forms $N_L$, each entry a decimal number and so an exact rational.

\begin{lemma}[the stored numbers, a computation]
\label{lemma:stored-numbers}
Each stored matrix is exactly symmetric.
The lists of residues and of exponents of $p_0$ and of the $N_L$ are closed under $S_4$, and the blocks of $H$ in one $S_4$-orbit of residues are exact copies of each other under the permutation of the exponents, also for the permutations that fix a block.
\end{lemma}

\begin{lemma}[the symmetry of the potential]
\label{lemma:potential-symmetry}
For each $\gamma\in G$, $p_0\circ\gamma=p_0$ and $N_L\circ\gamma=N_L$, and $p_0(\bar X)=p_0(X)$ and $N_L(\bar X)=N_L(X)$.
So $\eta$ and $|\nabla\eta|_g$ are invariant under $G$ and under complex conjugation.
\end{lemma}

\begin{lemma}[the chart does not matter]
\label{lemma:chart-independence}
Where $p_0>0$ and $i\partial\bar\partial\Phi$ is positive, $\eta$ and $|\nabla\eta|_g$ have the same values in every chart that contains the point.
\end{lemma}

\subsection{The model of the residual on a box}
\label{subsection:box-model}

\begin{definition}[a box]
\label{definition:box}
A box has the exact centre $c\in\CC^2$ of a cell or a sub-cube (\cref{subsection:cover-and-sign}), a radius $\delta$, a bound $h_{\max}$, an order $d$, a gap fraction $\vartheta>0$ and a reference centre $c_r$.
Let $\delta<\delta_f<\delta_\Phi$ be the double-precision values of $(1+\vartheta)\delta$ and $(1+\vartheta)^2\delta$.
The models of the potential live on $\mathbb D^4_{\delta_\Phi}$, the models of the metric and of $\eta$ on $\mathbb D^4_{\delta_f}$, and the bounds hold on $\mathbb D^4_\delta$.
\end{definition}

\begin{definition}[the frame]
\label{definition:frame}
From $c$, ball arithmetic at $512$ bits computes the principal fourth root $X_1(c)$ of $-1-c_1^4-c_2^4$ and the double-double numbers $p_1\approx X_1(c)$, $d_i\approx-c_i^3/p_1^3$, $e_3\approx\overline{\nabla W(p)}/|\nabla W(p)|$ and $c_*^{-1}\approx1/(\nabla W(p)\cdot e_3)$, where $p=(p_1,c_1,c_2)$, $e_1=(d_1,1,0)$, $e_2=(d_2,0,1)$, $E=(e_1,e_2,e_3)$ and $\nabla W\cdot e=\sum_k\partial_kW\,e_k$.
These numbers are exact from here on; the proof uses only $\det E\neq0$ and $c_*^{-1}\neq0$, which the computation proves.
Ball arithmetic also gives $1/(16|\det E|^2\bar v)$, $c_*^{-1}W(p)$, $c_*^{-1}\nabla W(p)\cdot e_i$ and $1-c_*^{-1}\nabla W(p)\cdot e_3$, rounded to double-double or triple-double numbers with bounds of their errors.
For $w\in\CC^2$ and $h\in\CC$ let $y(w,h)=p+w_1e_1+w_2e_2+he_3$.
\end{definition}

\begin{lemma}[the volume formula]
\label{lemma:volume-formula}
Let $h$ be holomorphic on an open set $V\subset\CC^2$ with $W(Y(w))=0$ for $Y(w)=y(w,h(w))$, and let $\nabla W(Y)\cdot e_3$, $\partial W/\partial X_1(Y)$ and $\det E$ be non-zero.
Then $Y$ is a local holomorphic parametrisation of $M$, and with $\Phi_{w\bar w}=(\partial_{w_i}\partial_{\bar w_j}(\Phi\circ Y))_{ij}$,
\begin{align}
\label{eqn:eta-box}
    \eta\circ Y
    =
    \frac{\det(\Phi_{w\bar w})\,|\nabla W(Y)\cdot e_3|^2}{16\,|\det E|^2\,\bar v},
    \qquad
    |\nabla\varphi|_g^2
    =
    \frac{4}{c_0}\,(\partial_w\varphi)^*\,\Phi_{w\bar w}^{-1}\,\partial_w\varphi,
    \qquad
    c_0=\frac{1}{\sqrt2\,\pi k_0},
\end{align}
for real functions $\varphi$ on $M$, where $c_0$ is the factor that gives $g$ the volume $1$.
\end{lemma}

\begin{lemma}[the points of $M$ on a box]
\label{lemma:contraction}
Let $T_w(h)=h-c_*^{-1}W(y(w,h))$, let $P_h$ be a polynomial in $w$, and let $b_k$ bound the coordinates of $y(w,P_h(w))$ on $\mathbb D^2_{\delta_\Phi}$.
Let $\beta\geq\sup|T_w(P_h(w))-P_h(w)|$ and $\lambda_0\geq\sup|T'_w(P_h(w))|$ on $\mathbb D^2_{\delta_\Phi}$, and
\[
    \lambda(\varrho)=\lambda_0+\varrho\,|c_*^{-1}|\sum_{k=1}^312\,|e_{3,k}|^2\,(b_k+\varrho|e_{3,k}|)^2 .
\]
Suppose that $\lambda(\varrho)<1$ and $\lambda(\varrho)\varrho+\beta\leq\varrho$ for some $\varrho>0$.
\begin{enumerate}
\item
For each $w\in\mathbb D^2_{\delta_\Phi}$ there is exactly one $h^*(w)$ with $|h^*(w)-P_h(w)|\leq\varrho$ and $W(y(w,h^*(w)))=0$, and $(P_h,\varrho)$ encloses the holomorphic function $h^*$.
\item
If also $\lambda(R)<1$ for $R=h_{\max}+\sup|P_h|+\varrho$, then each $h$ with $|h|\leq h_{\max}$ and $W(y(w,h))=0$ is $h^*(w)$.
\end{enumerate}
\end{lemma}

The code iterates $T_w$ on Taylor models in $w$ to find $P_h$, with $T_w$ written around $p$, so that only its small nonlinear terms are rounded.
The contraction factors $\lambda(R)$ of the proof are at most $0.203$.

\begin{definition}[the local forms]
\label{definition:local-forms}
On $M$, the relation $X_1^4=-1-X_2^4-X_3^4$ writes each of the forms $p_0$ and $N_L$ in the chart $X_0=1$ as $\sum_{e=0}^{3}|X_1|^{2e}\,m(t)^TQ_e\,\overline{m(t)}$, with the vector $m(t)$ of the monomials in $t$ and real symmetric matrices $Q_e$ of exact rationals.
For a reference centre $c_r$, let $x_r$ and $\sigma_r$ be real double-double numbers near $p_0(c_r)^{-1/k_0}$ and $|X_1(c_r)|^2$, let $S$ be the exact matrix that writes the monomials of $t$ in the monomials of $t-c_r$, and let $Q'_e=SQ_eS^*$.
The \emph{local forms} of $c_r$ are the matrices
\[
    R_j=\sum_{e=j}^{3}\binom ej\sigma_r^{e-j}Q'_e,
    \qquad
    j=0,\dots,3,
\]
of the base and of the combinations $\widetilde N_m=\sum_{L\geq m}\binom Lm x_r^LN_L$, $m=0,\dots,24$.
Ball arithmetic at $256$ bits computes them from the stored numbers, and their triple-double values are stored with bounds of the error and of the modulus of each entry.
\end{definition}

\begin{lemma}[the identities of the local forms]
\label{lemma:local-forms}
Let $\Delta=|X_1|^2-\sigma_r$ and $m'(t)=((t-c_r)^\alpha)_\alpha$.
\begin{enumerate}
\item
On $M$, each form is $\sum_{j=0}^{3}\Delta^j\,m'(t)^TR_j\,\overline{m'(t)}$ with its local forms $R_j$.
\item
Let $x=p_0^{-1/k_0}$, let $\zeta\in\CC$ with $x_c=x_r(1+\zeta)\neq0$, and let $\xi=x/x_c-1$.
Then $\sum_LN_Lx^L=\sum_{k=0}^{24}\mathcal F_k\,\xi^k$, where $\mathcal F_k$ is the form with the matrices $\sum_{m\geq k}(1+\zeta)^k\binom mk\zeta^{m-k}R_j(\widetilde N_m)$.
\item
All these matrices are Hermitian.
\end{enumerate}
\end{lemma}

On a box, ball arithmetic gives the weights $(1+\zeta)^k\binom mk\zeta^{m-k}$ for a double-double number $\zeta$ near $x(c)/x_r-1$, and a matrix product with error bounds gives each matrix of $\mathcal F_k$ with a bound of the error of each entry.
So the cancellations between the levels occur only in ball arithmetic, and the GPU sums a series in the small function $\xi$.

\begin{lemma}[the models in $w$]
\label{lemma:models-w}
Under the conditions of Lemma~\ref{lemma:contraction}, Taylor models on $\mathbb D^2_{\delta_\Phi}$ enclose the coordinates $Y_k$ of $Y(w)=y(w,h^*(w))$, the monomials $(t(Y)-c_r)^\alpha$ for $|\alpha|$ up to the degree of the forms, and $\nabla W(Y)\cdot e_3$.
Taylor models on $\mathbb D^4_{\delta_\Phi}$ enclose holomorphic functions $\tilde\Delta$ and $\tilde r$ of $(w,w')$ that are $|X_1|^2-\sigma_r$ and $1+|y|^2$ at $Y(w)$ on the slice $w'=\bar w$.
\end{lemma}

\begin{lemma}[the base]
\label{lemma:base}
Let $\tilde p=\tilde p_0\,\tilde r^{-k_0}$, where
\[
    \tilde p_0=\sum_{j=0}^{3}\tilde\Delta^j\tilde\rho_j,
    \qquad
    \tilde\rho_j(w,w')=m'(t(Y(w)))^TR_j\,\overline{m'(t(Y(\overline{w'})))},
\]
with the local forms $R_j$ of the base.
Let $s_r$ and $s_p$ be double-double numbers near $1/\tilde r(0)$ and $1/\Re\tilde p(0)$, with $s_p$ real.
If the models of $s_r\tilde r-1$ and $s_p\tilde p-1$ have bounds below $1$ on $\mathbb D^4_{\delta_\Phi}$, then binomial and logarithm series enclose $\tilde r^{-1}$, $\log(s_r\tilde r)$, $(s_p\tilde p)^{-1/k_0}$ and $\log(s_p\tilde p)$.
On the slice, $\tilde p=p_0(1+|y|^2)^{-k_0}$ and $s_p\tilde p\in(0,2)$, so $p_0\neq0$ and $p_0$ has the sign of $s_p$.
\end{lemma}

\begin{lemma}[the level sum]
\label{lemma:level-sum}
Let $s_p>0$, and let $x=\tilde r^{-1}s_p^{1/k_0}(s_p\tilde p)^{-1/k_0}$ and $\xi=x/x_c-1$.
The models of the polarised forms $\mathcal F_k$, composed at the order $d$ for $k\leq3$ and $\max(d-k+3,4)$ for $k>3$ with Horner's rule in $\tilde\Delta$, and Horner's rule in the model of $\xi$ give a model of a function $\tilde\psi$ on $\mathbb D^4_{\delta_\Phi}$ that is $\Phi-\log p_0$ at $Y(w)$ on the slice.
\end{lemma}

\begin{lemma}[the metric and $\eta$]
\label{lemma:metric-eta}
Let $s_p>0$.
The model of $\tilde\Phi=k_0\log(s_r\tilde r)+\log(s_p\tilde p)+\tilde\psi$ encloses a function that is $\Phi\circ Y$ up to a constant on the slice.
Cauchy estimates from $\delta_\Phi$ to $\delta_f$ give models of $\partial_{w_i}\partial_{w'_j}\tilde\Phi$, which are $\Phi_{w\bar w}$ on the slice, and \cref{eqn:eta-box}, with the model of $\nabla W(Y)\cdot e_3$, gives a model of a function on $\mathbb D^4_{\delta_f}$ that is $\eta\circ Y-1$ on the slice.
\end{lemma}

\begin{lemma}[the gradient]
\label{lemma:gradient}
Let $b_1,b_2$ bound the models of $\partial_{w_i}\eta$ on $\mathbb D^4_\delta$, from a Cauchy estimate from $\delta_f$ to $\delta$.
Let $\ell_1$ and $\ell_2$ be lower bounds of the real parts of the models of $\det\Phi_{w\bar w}$ and $(\Phi_{w\bar w})_{11}$ on $\mathbb D^4_\delta$, and $\theta$ an upper bound of $\operatorname{tr}\Phi_{w\bar w}$ there.
If $\ell_1>0$ and $\ell_2>0$, then on the slice $\Phi_{w\bar w}$ is positive definite, $\lambda_{\min}(\Phi_{w\bar w})\geq\lambda_{\mathrm{low}}=\ell_1/\theta$, and
\[
    |\nabla\eta|_g\leq\frac{2\,(b_1^2+b_2^2)^{1/2}}{(c_0\lambda_{\mathrm{low}})^{1/2}}.
\]
\end{lemma}

\begin{theorem}[one box]
\label{theorem:one-box}
Let a box give a finite bound $B$ of the model of $\eta-1$ on $\mathbb D^4_\delta$, and $\lambda(R)<1$ in Lemma~\ref{lemma:contraction}.
Then for each $w\in\mathbb D^2_\delta$ and each $h$ with $|h|\leq h_{\max}$, $W(y(w,h))=0$ and $X_1\neq0$ at $y(w,h)$:
\begin{enumerate}
\item
$p_0(y(w,h))\neq0$, and it has the sign of $s_p$;
\item
if $s_p>0$, then $|\eta-1|\leq B$ at $y(w,h)$;
\item
if $s_p>0$ and $\lambda_{\mathrm{low}}>0$, then $g$ is positive definite at $y(w,h)$, and $|\nabla\eta|_g$ is at most the bound of Lemma~\ref{lemma:gradient} there.
\end{enumerate}
\end{theorem}

\subsection{The cover and the sign of \texorpdfstring{$p_0$}{p0}}
\label{subsection:cover-and-sign}

\begin{definition}[the grid]
\label{definition:grid}
Let $\ell$ be the double-precision number nearest to $0.03$.
The grid consists of the $943{,}296$ cubes
\[
    \prod_{k=1}^4\bigl[i_k\ell,(i_k+1)\ell\bigr],
    \qquad
    0\leq i_1,i_3,i_4\leq33,
    \quad
    0\leq i_2\leq23,
\]
in the coordinates $(\Re X_2,\Im X_2,\Re X_3,\Im X_3)$, with integers $i_k$.
A sub-cube of depth $j$ of a cube is a cube of side $\ell/2^j$ in it with corners on the same lattice.
\end{definition}

\begin{lemma}[the grid holds $\mathcal F$]
\label{lemma:grid}
The cubes of the grid cover $\mathcal F$, and the corners and centres of the cubes and of their sub-cubes of depth at most $3$ have exact double-double values.
\end{lemma}

\begin{lemma}[the cubes that miss $\mathcal F$]
\label{lemma:miss}
If interval arithmetic \cite{mpmath} proves that one of the inequalities
\[
    \Im X_2\leq\Re X_2,
    \quad
    \Re X_2,\Im X_2,\Re X_3,\Im X_3\geq0,
    \quad
    |X_3|^2\leq|X_2|^2,
    \quad
    |X_1|^8\leq1,
    \quad
    |X_2|^8\leq|X_1|^8,
\]
with $|X_1|^8=|1+X_2^4+X_3^4|^2$, fails on all of a cube, then the cube misses $\mathcal F$.
\end{lemma}

\begin{lemma}[the box of a cube]
\label{lemma:box-radius}
Let $K$ be a cube of the grid or a sub-cube, $c$ its centre, and $p$, $E$ the frame of $c$.
If ball arithmetic gives $|1+X_2^4+X_3^4|>0$ on $K$ and finite numbers $\delta$ and $h_{\max}$, then:
\begin{enumerate}
\item
there is exactly one continuous fourth root $X_1$ of $-1-X_2^4-X_3^4$ on $K$ with the value $X_1(c)$ of the frame at $c$; it is holomorphic on a neighbourhood of $K$, and the other roots are its images under phases;
\item
for each $t\in K$ the point $(X_1(t),t)$ is $y(w,h)$ with $w\in\mathbb D^2_\delta$ and $|h|\leq h_{\max}$.
\end{enumerate}
\end{lemma}

The numbers $\delta$ and $h_{\max}$ come from a bound of the second derivatives of $X_1$ on $K$ and Taylor's formula, in exact rationals from directed bounds, rounded upwards.
For the boxes of the proof, $0.0109\leq\delta\leq0.0252$ and $h_{\max}\leq\num{6.0e-3}$.

\begin{lemma}[the partition, a computation]
\label{lemma:partition}
A separate program proves, with the exact corners of the cubes:
\begin{enumerate}
\item
each cube of the grid that is not a cell gives a proof by Lemma~\ref{lemma:miss} that it misses $\mathcal F$;
\item
the cubes of the boxes of the proof, together with sub-cubes that miss $\mathcal F$ by Lemma~\ref{lemma:miss}, form a partition of each cell;
\item
$\delta$ and $h_{\max}$ of each box, computed again by Lemma~\ref{lemma:box-radius} from its cube, are at most the values that its computation used;
\item
each bound of $|\eta-1|$ lies in $(0,\num{5e-20}]$, each contraction factor $\lambda(R)$ is below $1$, each $\lambda_{\mathrm{low}}$ is at least $0.9985$, and each bound of the gradient and of the metric is finite.
\end{enumerate}
\end{lemma}

\begin{lemma}[the sign of $p_0$, a computation]
\label{lemma:sign}
Ball arithmetic gives $p_0(z^*)=426944.62554611\pm\num{4.6e-9}$ at $z^*=(1,e^{i\pi/4},0,0)$.
If $p_0\neq0$ on $M$, then $p_0>0$ on $M$.
\end{lemma}

The second statement holds because $p_0$ is real, $G$-invariant and homogeneous, and a path in $M$ joins each point to an image of $z^*$ under $G$, as in the proof of \cref{subsection:diameter-bound}.
This lemma is necessary, because the coefficient matrix of $p_0$ is not positive definite: only $2$ of its $64$ blocks are positive definite, and most of the others are singular up to their rounding.

\begin{proof}[Proof of Proposition~\ref{proposition:certified-residual}]
Let $x\in M$.
By step 1 of \cref{subsection:certified-errors-bounds} and Lemma~\ref{lemma:potential-symmetry}, an element of $G$ or a conjugation maps $x$ to a point over $\mathcal F$ and keeps $\eta$, $|\nabla\eta|_g$ and the sign of $p_0$.
By Lemmas~\ref{lemma:grid}, \ref{lemma:miss} and~\ref{lemma:partition}, this point lies over the cube of a box of the proof.
By Lemma~\ref{lemma:box-radius} and Lemma~\ref{lemma:partition}(3), a phase of $X_1$ maps it to a point $y(w,h)$ of that box with $w\in\mathbb D^2_\delta$ and $|h|\leq h_{\max}$.
So Theorem~\ref{theorem:one-box}(1) gives $p_0\neq0$ on $M$, and Lemma~\ref{lemma:sign} gives $p_0>0$ on $M$; hence $s_p>0$ for every box.
Then Theorem~\ref{theorem:one-box}(2) and~(3) give the bounds, and their largest values over the boxes are the numbers of Proposition~\ref{proposition:certified-residual}.
\end{proof}

\subsection{Bounds for the metric}
\label{subsection:certified-metric-bounds}

The same Taylor models give the computer-assisted bounds of \cref{section:bounding-quantities}.
Let $g_{\mathrm{FS}}$ be the matrix of $\partial\bar\partial\log(1+|y|^2)$.
A \emph{chart} $a=(i,j)$ is a pair of indices $i\neq j$; its coordinates are $t_a=(X_k/X_i,X_l/X_i)$ with $\{k,l\}=\{0,1,2,3\}\setminus\{i,j\}$, and its domain is $D_a=\{|X_i|>0.99|X_x|\text{ for }x\neq i,\ |X_j|>0.98|X_x|\text{ for }x\notin\{i,j\}\}$.
Up to a linear isometry of its coordinates, each parent chart of \cref{section:bounding-quantities} is a chart $a$ with its domain in $D_a$, such as $a=(0,1)$ for $U^*$.
Let $G_a$ be the real matrix of the metric of volume $1$ in the real coordinates of $t_a$, and $k_a$ the Frobenius norm of the real tensor $DG_a$ normalised by $G_a^{1/2}$ at the point; then $\tilde K$, $\lambda$ and $\Lambda$ of \cref{equation:eigenvalues-metric,equation:3rd-der-potential} are bounded by the supremum of $k_a$ and the extreme eigenvalues of $G_a$ over all charts $a$ and $x\in D_a$.

\begin{lemma}[the reduction of the atlas]
\label{lemma:atlas-reduction}
For each chart $a$ and $x\in D_a$ there are an element $\gamma$ of $G$, possibly composed with the conjugation, and a chart $b$ such that $\gamma x$ lies over $\mathcal F$ and in $D_b$, $k_b(\gamma x)=k_a(x)$, and $G_b(\gamma x)$ is an orthogonal conjugate of $G_a(x)$.
\end{lemma}

\begin{lemma}[the charts of a cube]
\label{lemma:chart-mask}
If interval arithmetic proves on a cube that $|X_i|\leq0.99|X_x|$ for some $x\neq i$, or $|X_j|\leq0.98|X_x|$ for some $x\notin\{i,j\}$, then no point of $M$ over the cube lies in $D_{(i,j)}$.
\end{lemma}

\begin{lemma}[eigenvalues of Hermitian $2\times2$ models]
\label{lemma:eigenvalues}
\begin{enumerate}
\item
For Hermitian positive definite $A$, $B$ and an invertible constant matrix $P$, the eigenvalues of $A^{-1}B$ lie in
\[
    \bigl[\lambda_{\min}(B')/\lambda_{\max}(A'),\ \lambda_{\max}(B')/\lambda_{\min}(A')\bigr],
    \qquad
    A'=P^TA\bar P,
    \quad
    B'=P^TB\bar P.
\]
\item
For Taylor models of a $2\times2$ matrix that is Hermitian on the slice, lower bounds of the real parts of the diagonal entries and upper bounds of the moduli of the entries and of the difference of the diagonal entries give, by the closed formula for the eigenvalues and by Gershgorin's discs, a lower bound of $\lambda_{\min}$ and an upper bound of $\lambda_{\max}$ at each point of the slice.
\end{enumerate}
\end{lemma}

\begin{lemma}[the metric against Fubini--Study, and Ricci]
\label{lemma:metric-fs-ricci}
On a box, with $P$ from $g_{\mathrm{FS}}$ at the centre, Lemma~\ref{lemma:eigenvalues} applied to the models of $g_{\mathrm{FS}}$ and $\Phi_{w\bar w}$ bounds the eigenvalues of $(k_0g_{\mathrm{FS}})^{-1}\Phi_{w\bar w}$, that is, of $\omega$ relative to $\omega_0$.
Since $\Ric(\omega)=-i\partial\bar\partial\log\eta$, each eigenvalue of $\Ric$ relative to $g$ has modulus at most
\[
    \frac{2}{c_0\,\lambda_{\min}(\Phi_{w\bar w})}
    \left(
    \frac{|\partial_w\partial_{\bar w}\eta|}{\inf\eta}
    +
    \frac{|\partial_w\eta|^2}{(\inf\eta)^2}
    \right)
\]
on the box, with Frobenius norms.
\end{lemma}

\begin{lemma}[the normalised derivative of the metric]
\label{lemma:k-bound}
In a chart $a$, the models of $t_a$, of its Jacobian $J_a$ and of $J_a^{-1}$ give models of $g_a=J_a^{-T}\Phi_{w\bar w}\bar J_a^{-1}$ and of its derivatives $\Gamma_a$, after a congruence with a constant matrix from the centre.
At each point of the box, $k_a$ is at most each of these two numbers:
\begin{enumerate}
\item
$\sqrt{8/c_0}$ times the supremum of $|\Gamma_a|_F$ over the infimum of $\lambda_{\min}(g_a)^{3/2}$ on the box;
\item
$\sqrt{8/c_0}$ times the square root of the bound of one Taylor model of the squared Frobenius norm of the normalised tensor, formed with the model of $g_a^{-1}$ at each point.
\end{enumerate}
The eigenvalues of $G_a=c_0g_a$ lie in the bounds of Lemma~\ref{lemma:eigenvalues} for $g_a$, times $c_0$.
\end{lemma}

On the boxes of Proposition~\ref{proposition:certified-residual} the first bound of $k_a$ is loose in the strongly distorted charts near $|X_1|\approx|X_2|\approx|X_3|\approx0.8$: there $k_a$ is $48$ to $55$ at the box centre, but the bound on the box reaches $439$.
The bounds of the metric use the models of $\Phi_{w\bar w}$ only up to the order $10$, so a separate computation bounds the metric alone, with Taylor models of the order $20$ on the cells and $12$ on their sub-cubes, and splits a cube, at most $3$ times, while the bound of $k_a$ exceeds $60$.

\begin{lemma}[the separate computation for the metric, a computation]
\label{lemma:metric-pass}
Its $261{,}373$ boxes satisfy (1) to (3) of Lemma~\ref{lemma:partition}, each contraction factor is below $1$, and each bound of the metric is finite.
Over these boxes, for the charts that Lemma~\ref{lemma:chart-mask} does not exclude, $k_a<59.9992$, the eigenvalues of $G_a$ lie in $(0.016109,0.64047)$, and the eigenvalues of $\omega$ relative to $\omega_0$ lie in $(0.596699,2.176518)$.
\end{lemma}

\begin{lemma}[the transition maps, a computation]
\label{lemma:transitions}
In ball arithmetic, a branch-and-bound over boxes of $D_{(0,1)}$ shows that the transition maps $F$ from the chart $(0,1)$ to the charts of the $11$ other pairs satisfy $\sup|D^2F|\leq\num{44.25}$ and $\sup|D^3F|\leq\num{882.2}$ in real norms.
\end{lemma}

\begin{lemma}[the transitions between two charts]
\label{lemma:chart-pairs}
For charts $a,a'$, a point $x$ in both domains, and $n=2,3$, let $N_n(a,a')$ be the norm of the normalised derivative $D^nF_{aa'}$ of \cref{equation:parent-transition-second-derivative-bound,equation:parent-transition-third-derivative-bound} at $x$.
If $x$ also lies in the domain of a chart $b$, then
\[
    N_2(a,a')\leq N_2(a,b)+N_2(a',b),
    \qquad
    N_3(a,a')\leq N_3(a,b)+N_3(a',b)+3N_2(a,b)N_2(a',b)+3N_2(a',b)^2 .
\]
\end{lemma}

These lemmas give the constants of \cref{section:bounding-quantities}.
Lemmas~\ref{lemma:metric-fs-ricci} and~\ref{lemma:metric-pass} give $\num{0.5966}\,\omega_0\leq\omega\leq\num{\numMetricComparison}\,\omega_0$.
On the boxes of Proposition~\ref{proposition:certified-residual}, the largest bound of Lemma~\ref{lemma:metric-fs-ricci} for $\Ric$ is $\num{1.334e-12}$, which gives $\kappa=\num{\numRicciBound}$ in \cref{eqn:ricci-bound}.
Lemmas~\ref{lemma:atlas-reduction} and~\ref{lemma:metric-pass} give $\tilde K=\num{\numKtilde}$, $\lambda=\num{\numLambdaLower}$ and $\Lambda=\num{\numLambdaUpper}$.
With Lemma~\ref{lemma:transitions}, the normalised second and third derivatives of the transitions with the chart $(0,1)$ are at most $\Lambda^{1/2}\lambda^{-1}\sup|D^2F|$ and $\Lambda^{1/2}\lambda^{-3/2}\sup|D^3F|$, and Lemmas~\ref{lemma:atlas-reduction} and~\ref{lemma:chart-pairs} with $b=(0,1)$ give \cref{equation:parent-transition-second-derivative-bound,equation:parent-transition-third-derivative-bound} for all pairs of charts with
\[
    \tilde c_1=2\Lambda^{1/2}\lambda^{-1}\sup|D^2F|,
    \qquad
    \tilde c_2=2\Lambda^{1/2}\lambda^{-3/2}\sup|D^3F|+6\bigl(\Lambda^{1/2}\lambda^{-1}\sup|D^2F|\bigr)^2,
\]
which are at most $\num{\numCOneTilde}$ and $\num{\numCTwoTilde}$.

\subsection{The computation}
\label{subsection:certified-computation}

\Cref{tab:certified-parameters} gives the parameters of the computation.
Each parameter is a free choice: the proof holds for every choice, and only the size of the bounds and the cost depend on it.

\begin{table}[ht]
\centering
\small
\begin{tabular}{@{}lp{0.62\textwidth}@{}}
\toprule
item & value\\
\midrule
grid & cubes of side $0.03$; $113{,}357$ cells may meet $\mathcal F$\\
reference centres & $332$, the centres of the cubes of side $0.18$ that contain cells\\
models of $\eta-1$ & order $28$ and $\vartheta=0.35$ if $r<0.3$ at the centre; order $32$ and $\vartheta=0.12$ otherwise\\
target & $\num{5e-20}$; a box above it uses the order $32$, then splits into the sub-cubes of half the side, at most twice\\
series & enough terms for a rest below $10^{-34}$\\
the metric alone & order $20$ on the cells and $12$ on the sub-cubes; split while the bound of $\tilde K$ is above $60$, at most $3$ times\\
ball arithmetic & $512$ bits for the frames, $256$ bits for the local forms\\
\bottomrule
\end{tabular}
\caption{The parameters of the computer-assisted proof.}
\label{tab:certified-parameters}
\end{table}

The computation of Proposition~\ref{proposition:certified-residual} evaluated $149{,}941$ boxes, and $134{,}552$ of them are in the proof: $9{,}793$ cells at the order $28$, $97{,}064$ cells at the order $32$, and $27{,}695$ sub-cubes at the order $32$ from $3{,}538$ split cells.
The other sub-cubes of the split cells miss $\mathcal F$.
The median bound of $|\eta-1|$ is $\num{1.8e-21}$, and $99.9\%$ of the bounds are below $\num{4.8e-20}$.
The bounds of $|\partial_w\eta|$ are at most $\num{2.43e-17}$, and $\lambda_{\min}(\Phi_{w\bar w})\geq\num{0.9985}$ on every box.
The largest bound of $|\nabla\eta|_g$ is at $r\approx0.92$.

The computation ran in four parts, each on one NVIDIA H200 GPU, in $2.9$ to $3.4$ hours per part, $12.5$ GPU-hours in total, or $0.3$ seconds per box on average.
This includes the bounds for the metric on the same boxes.
The local forms of the $332$ reference centres take $1.5$ hours on $48$ CPU cores.
The separate computation for the metric took $1.5$ to $1.7$ hours on each of four GPUs.
The exact check takes less than a minute on $16$ CPU cores.

\section{AI statement}

LLMs were used in the preparation of the manuscript and accompanying code.
The authors assume responsibility for all content.

\bibliographystyle{amsplain}
\bibliography{library}

\end{document}

%% file: preamble-fabian.tex
\usepackage[english]{babel}
\usepackage[utf8]{inputenc}
\usepackage[T1]{fontenc}
\usepackage[english,cleanlook]{isodate}
\usepackage{lmodern}
\usepackage[babel=true]{microtype}
\usepackage{amsmath}
\usepackage{mathtools}
\usepackage{amssymb}
\usepackage{amsthm}
\usepackage{mathrsfs}
\usepackage{enumitem}
\usepackage{hyperref}
\usepackage{esint}

\usepackage[margin=3.2cm]{geometry}

\newcommand{\CC}{\mathbb{C}}
\renewcommand{\P}{\mathbb{P}}

\renewcommand{\Im}{\operatorname{Im}}
\renewcommand{\Re}{\operatorname{Re}}

\newcommand{\Sym}{\operatorname{Sym}}

\newcommand{\Stab}{\operatorname{Stab}}

\newcommand{\Ric}{\operatorname{Ric}}

\let\temp\|
\renewcommand{\|}[1]{\temp #1 \temp}

\newcommand{\vol}{\operatorname{vol}}
\newcommand{\diam}{\operatorname{diam}}

\renewcommand{\tilde}[1]{\widetilde{#1}}

\usepackage[capitalise]{cleveref}

\numberwithin{equation}{section}

\newtheorem{proposition}[equation]{Proposition}
\crefname{proposition}{Proposition}{Propositions}

\newtheorem{lemma}[equation]{Lemma}
\crefname{lemma}{Lemma}{Lemmas}

\newtheorem{corollary}[equation]{Corollary}
\crefname{corollary}{Corollary}{Corollaries}

\newtheorem*{corollary*}{Corollary}
\crefname{corollary*}{Corollary}{Corollaries}

\newtheorem{theorem}[equation]{Theorem}
\crefname{theorem}{Theorem}{Theorems}

\crefname{exercise}{Exercise}{Exercises}

\crefname{task}{Task}{Task}

\crefname{conjecture}{Conjecture}{Conjectures}

\crefname{algorithm}{Algorithm}{algorithms}

\newtheorem*{theorem*}{Theorem}
\crefname{theorem*}{Theorem}{Theorems}

\crefname{claim}{Claim}{Claims}

\theoremstyle{remark}

\crefname{question}{Question}{Questions}

\newtheorem{definition}[equation]{Definition}
\crefname{definition}{Definition}{Definitions}

\crefname{example}{Example}{Examples}

\newtheorem{remark}[equation]{Remark}
\crefname{remark}{Remark}{Remarks}

\crefname{assumption}{Assumption}{Assumptions}

\crefname{problem}{Problem}{Problems}